\documentclass[10pt]{smfart}
\usepackage{a4wide}
\usepackage{amsmath,amssymb,smfthm}
\usepackage[francais]{babel}
\usepackage[latin1]{inputenc}
\usepackage[all]{xy}
\usepackage{graphicx}
\usepackage{epic,eepic}
\usepackage{epsfig}
\usepackage{color}

\def\R{\mathbb{R}}

\def\C{\mathbb{C}}

\def\Fq{\mathbb{F}_q}

\def\L{\mathcal{L}}
\def\La{\Lambda}

\def\O{\mathcal{O}}

\def\*{^\times }
\def\dpt{\displaystyle}

\def\Gm{\mathbb{G}_m}         
\def\l{\lambda}
\def\g{\gamma}

\def\a{\alpha}
\def\b{\beta}

\def\s{\sigma}
\def\ph{\varphi}

\def\e{\epsilon}
\def\ssi{\Leftrightarrow}
\def\lssi{\Longleftrightarrow}
\def\impl{\Rightarrow}
\def\limpl{\Longrightarrow}
\def\drt{\rightarrow}
\def\ldrt{\longrightarrow}
\def\Q{\mathbb{Q}}

\def\Qp{\mathbb{Q}_p}

\def\Zp{\mathbb{Z}_p}
\def\Z{\mathbb{Z}}
\def\N{\mathbb{N}}

\def\Gal{\text{Gal}}
\def\={\! = \!}

\def\spec{\text{Spec}}
\def\spf{\text{Spf}}

\def\limp{\underset{\longleftarrow}{\text{ lim }}\;}

\def\iso{\xrightarrow{\;\sim\;}}

\def\End{\text{End}}
\def\Aut{\text{Aut}}
\def\GL{\hbox{GL}}

\def\xrig{\xrightarrow}

\def\M{\mathcal{M}}

\def\X{\mathfrak{X}}

\def\bc{\backslash}
\def\AA{\mathcal{A}}

\def\<{<\hspace{-1mm}}
\def\>{\hspace{-1mm}>}

\def\Lie{\text{Lie}}

\def\dem{{\it Démonstration. }}

\def\Fil{\text{Fil}}

\def\LT{\mathcal{L}\mathcal{T}}
\def\D{{\mathcal{D}r}}

\author{Laurent Fargues}
\address{CNRS-IHES-université Paris-Sud Orsay}
\email{laurent.fargues@math.u-psud.fr}
\date{}

\begin{document}

\newtheorem{Fait}{Fait}

\title[Application de Hodge-Tate  et
immeuble du groupe linéaire]{Application de Hodge-Tate duale d'un groupe de Lubin-Tate, immeuble
de Bruhat-Tits du groupe linéaire et filtrations de
ramification}
\maketitle

\begin{figure}[htbp]
   \begin{center}
\begin{picture}(100,80)
\path(0,0)(100,0)
\path(0,0)(50,86)
\path(50,86)(100,0)
\path(56,19)(38,27)
\path(56,19)(68,27)
\path(38,27)(39,13)
\path(54,39)(41,45)
\path(38,27)(54,39)
\path(54,39)(56,19)
\path(57,9)(39,13)
\path(57,9)(69,13)
\path(29,31)(27,19)
\path(62,45)(53,53)
\path(29,31)(41,45)
\path(62,45)(68,27)
\path(39,13)(27,19)
\path(69,13)(77,19)
\path(27,19)(28,9)
\path(53,53)(44,57)
\path(41,45)(53,53)
\path(68,27)(69,13)
\path(58,4)(40,6)
\path(58,4)(70,6)
\path(25,33)(20,22)
\path(66,47)(58,57)
\path(25,33)(35,47)
\path(66,47)(73,31)
\path(40,6)(28,9)
\path(70,6)(78,9)
\path(20,22)(19,13)
\path(58,57)(52,63)
\path(35,47)(44,57)
\path(73,31)(77,19)
\path(28,9)(19,13)
\path(78,9)(84,13)
\path(19,13)(19,6)
\path(52,63)(45,66)
\path(44,57)(52,63)
\path(77,19)(78,9)
\path(58,2)(40,3)
\path(58,2)(70,3)
\path(22,34)(17,23)
\path(68,49)(61,59)
\path(22,34)(32,49)
\path(68,49)(76,33)
\path(40,3)(28,4)
\path(70,3)(79,4)
\path(17,23)(14,15)
\path(61,59)(56,66)
\path(32,49)(39,59)
\path(76,33)(81,22)
\path(28,4)(19,6)
\path(79,4)(84,6)
\path(14,15)(13,9)
\path(56,66)(51,70)
\path(39,59)(45,66)
\path(81,22)(84,13)
\path(19,6)(13,9)
\path(84,6)(88,9)
\path(13,9)(13,4)
\path(51,70)(47,72)
\path(45,66)(51,70)
\path(84,13)(84,6)
\path(58,1)(41,1)
\path(58,1)(70,1)
\path(21,35)(16,24)
\path(69,50)(63,60)
\path(21,35)(30,50)
\path(69,50)(77,34)
\path(41,1)(28,2)
\path(70,1)(79,2)
\path(16,24)(12,16)
\path(63,60)(58,67)
\path(30,50)(37,60)
\path(77,34)(83,23)
\path(28,2)(20,3)
\path(79,2)(85,3)
\path(12,16)(10,11)
\path(58,67)(54,72)
\path(37,60)(42,67)
\path(83,23)(86,15)
\path(20,3)(13,4)
\path(85,3)(89,4)
\path(10,11)(9,6)
\path(54,72)(51,74)
\path(42,67)(47,72)
\path(86,15)(88,9)
\path(13,4)(9,6)
\path(89,4)(92,6)
\path(9,6)(9,3)
\path(51,74)(48,76)
\path(47,72)(51,74)
\path(88,9)(89,4)
\path(58,0)(41,0)
\path(58,0)(70,0)
\path(21,35)(15,24)
\path(70,50)(63,60)
\path(21,35)(30,50)
\path(70,50)(78,35)
\path(41,0)(28,1)
\path(70,0)(79,1)
\path(15,24)(11,17)
\path(63,60)(59,68)
\path(30,50)(36,60)
\path(78,35)(84,24)
\path(28,1)(20,1)
\path(79,1)(85,1)
\path(11,17)(8,11)
\path(59,68)(55,73)
\path(36,60)(41,68)
\path(84,24)(88,16)
\path(20,1)(14,2)
\path(85,1)(89,2)
\path(8,11)(7,7)
\path(55,73)(53,76)
\path(41,68)(44,73)
\path(88,16)(90,11)
\path(14,2)(9,3)
\path(89,2)(92,3)
\path(7,7)(6,4)
\path(53,76)(50,78)
\path(44,73)(48,76)
\path(90,11)(92,6)
\path(9,3)(6,4)
\path(92,3)(94,4)
\path(6,4)(6,2)
\path(50,78)(48,79)
\path(48,76)(50,78)
\path(92,6)(92,3)
\path(58,0)(41,0)
\path(58,0)(70,0)
\path(20,35)(14,25)
\path(70,50)(64,61)
\path(20,35)(29,50)
\path(70,50)(79,35)
\path(41,0)(29,0)
\path(70,0)(79,0)
\path(14,25)(10,17)
\path(64,61)(59,68)
\path(29,50)(35,61)
\path(79,35)(84,24)
\path(29,0)(20,0)
\path(79,0)(85,0)
\path(10,17)(7,12)
\path(59,68)(56,73)
\path(35,61)(40,68)
\path(84,24)(89,17)
\path(20,0)(14,1)
\path(85,0)(89,1)
\path(7,12)(6,8)
\path(56,73)(54,77)
\path(40,68)(43,73)
\path(89,17)(91,11)
\path(14,1)(9,1)
\path(89,1)(92,1)
\path(6,8)(5,5)
\path(54,77)(52,79)
\path(43,73)(46,77)
\path(91,11)(93,7)
\path(9,1)(6,2)
\path(92,1)(94,2)
\path(5,5)(4,3)
\path(52,79)(50,80)
\path(46,77)(48,79)
\path(93,7)(94,4)
\path(6,2)(4,3)
\path(94,2)(96,3)
\path(4,3)(4,1)
\path(50,80)(49,81)
\path(48,79)(50,80)
\path(94,4)(94,2)
\path(58,0)(41,0)
\path(58,0)(70,0)
\path(20,35)(14,25)
\path(70,50)(64,61)
\path(20,35)(29,50)
\path(70,50)(79,35)
\path(41,0)(29,0)
\path(70,0)(79,0)
\path(14,25)(10,17)
\path(64,61)(60,68)
\path(29,50)(35,61)
\path(79,35)(85,25)
\path(29,0)(20,0)
\path(79,0)(85,0)
\path(10,17)(7,12)
\path(60,68)(56,73)
\path(35,61)(40,68)
\path(85,25)(89,17)
\path(20,0)(14,0)
\path(85,0)(89,0)
\path(7,12)(5,8)
\path(56,73)(54,77)
\path(40,68)(43,73)
\path(89,17)(92,12)
\path(14,0)(10,0)
\path(89,0)(92,0)
\path(5,8)(4,5)
\path(54,77)(52,79)
\path(43,73)(45,77)
\path(92,12)(94,8)
\path(10,0)(7,1)
\path(92,0)(94,1)
\path(4,5)(3,3)
\path(52,79)(51,81)
\path(45,77)(47,79)
\path(94,8)(95,5)
\path(7,1)(4,1)
\path(94,1)(96,1)
\path(3,3)(3,2)
\path(51,81)(50,82)
\path(47,79)(49,81)
\path(95,5)(96,3)
\path(4,1)(3,2)
\path(96,1)(97,2)
\path(3,2)(3,1)
\path(50,82)(49,83)
\path(49,81)(50,82)
\path(96,3)(96,1)
\path(58,0)(41,0)
\path(58,0)(70,0)
\path(20,35)(14,25)
\path(70,50)(64,61)
\path(20,35)(29,50)
\path(70,50)(79,35)
\path(41,0)(29,0)
\path(70,0)(79,0)
\path(14,25)(10,17)
\path(64,61)(60,68)
\path(29,50)(35,61)
\path(79,35)(85,25)
\path(29,0)(20,0)
\path(79,0)(85,0)
\path(10,17)(7,12)
\path(60,68)(57,73)
\path(35,61)(39,68)
\path(85,25)(89,17)
\path(20,0)(14,0)
\path(85,0)(89,0)
\path(7,12)(5,8)
\path(57,73)(54,77)
\path(39,68)(43,73)
\path(89,17)(92,12)
\path(14,0)(10,0)
\path(89,0)(92,0)
\path(5,8)(3,6)
\path(54,77)(53,80)
\path(43,73)(45,77)
\path(92,12)(94,8)
\path(10,0)(7,0)
\path(92,0)(94,0)
\path(3,6)(3,4)
\path(53,80)(52,81)
\path(45,77)(46,80)
\path(94,8)(95,5)
\path(7,0)(4,0)
\path(94,0)(96,0)
\path(3,4)(2,2)
\path(52,81)(51,83)
\path(46,80)(48,81)
\path(95,5)(96,3)
\path(4,0)(3,1)
\path(96,0)(97,1)
\path(2,2)(2,1)
\path(51,83)(50,83)
\path(48,81)(49,83)
\path(96,3)(97,2)
\path(3,1)(2,1)
\path(97,1)(98,1)
\path(2,1)(2,0)
\path(50,83)(49,84)
\path(49,83)(50,83)
\path(97,2)(97,1)
\path(58,0)(41,0)
\path(58,0)(70,0)
\path(20,35)(14,25)
\path(70,50)(64,61)
\path(20,35)(29,50)
\path(70,50)(79,35)
\path(41,0)(29,0)
\path(70,0)(79,0)
\path(14,25)(10,17)
\path(64,61)(60,68)
\path(29,50)(35,61)
\path(79,35)(85,25)
\path(29,0)(20,0)
\path(79,0)(85,0)
\path(10,17)(7,12)
\path(60,68)(57,74)
\path(35,61)(39,68)
\path(85,25)(89,17)
\path(20,0)(14,0)
\path(85,0)(89,0)
\path(7,12)(5,8)
\path(57,74)(54,77)
\path(39,68)(42,74)
\path(89,17)(92,12)
\path(14,0)(10,0)
\path(89,0)(92,0)
\path(5,8)(3,6)
\path(54,77)(53,80)
\path(42,74)(45,77)
\path(92,12)(94,8)
\path(10,0)(7,0)
\path(92,0)(94,0)
\path(3,6)(2,4)
\path(53,80)(52,82)
\path(45,77)(46,80)
\path(94,8)(96,6)
\path(7,0)(4,0)
\path(94,0)(96,0)
\path(2,4)(2,2)
\path(52,82)(51,83)
\path(46,80)(47,82)
\path(96,6)(97,4)
\path(4,0)(3,0)
\path(96,0)(97,0)
\path(2,2)(1,1)
\path(51,83)(50,84)
\path(47,82)(48,83)
\path(97,4)(97,2)
\path(3,0)(2,0)
\path(97,0)(98,0)
\path(1,1)(1,1)
\path(50,84)(50,84)
\path(48,83)(49,84)
\path(97,2)(98,1)
\path(2,0)(1,1)
\path(98,0)(98,1)
\path(1,1)(1,0)
\path(50,84)(49,84)
\path(49,84)(50,84)
\path(98,1)(98,0)
\path(58,0)(41,0)
\path(58,0)(70,0)
\path(20,35)(14,25)
\path(70,50)(64,61)
\path(20,35)(29,50)
\path(70,50)(79,35)
\path(41,0)(29,0)
\path(70,0)(79,0)
\path(14,25)(10,17)
\path(64,61)(60,68)
\path(29,50)(35,61)
\path(79,35)(85,25)
\path(29,0)(20,0)
\path(79,0)(85,0)
\path(10,17)(7,12)
\path(60,68)(57,74)
\path(35,61)(39,68)
\path(85,25)(89,17)
\path(20,0)(14,0)
\path(85,0)(89,0)
\path(7,12)(5,8)
\path(57,74)(55,77)
\path(39,68)(42,74)
\path(89,17)(92,12)
\path(14,0)(10,0)
\path(89,0)(92,0)
\path(5,8)(3,6)
\path(55,77)(53,80)
\path(42,74)(45,77)
\path(92,12)(94,8)
\path(10,0)(7,0)
\path(92,0)(94,0)
\path(3,6)(2,4)
\path(53,80)(52,82)
\path(45,77)(46,80)
\path(94,8)(96,6)
\path(7,0)(5,0)
\path(94,0)(96,0)
\path(2,4)(1,2)
\path(52,82)(51,83)
\path(46,80)(47,82)
\path(96,6)(97,4)
\path(5,0)(3,0)
\path(96,0)(97,0)
\path(1,2)(1,2)
\path(51,83)(51,84)
\path(47,82)(48,83)
\path(97,4)(97,2)
\path(3,0)(2,0)
\path(97,0)(98,0)
\path(1,2)(1,1)
\path(51,84)(50,84)
\path(48,83)(49,84)
\path(97,2)(98,1)
\path(2,0)(1,0)
\path(98,0)(98,0)
\path(1,1)(1,0)
\path(50,84)(50,85)
\path(49,84)(49,84)
\path(98,1)(98,1)
\path(1,0)(1,0)
\path(98,0)(99,0)
\path(1,0)(1,0)
\path(50,85)(49,85)
\path(49,84)(50,85)
\path(98,1)(98,0)
\path(54,39)(56,19)
\path(54,39)(56,19)
\path(56,19)(38,27)
\path(38,27)(54,39)
\path(56,19)(38,27)
\path(38,27)(54,39)
\path(56,19)(57,9)
\path(56,19)(57,9)
\path(38,27)(29,31)
\path(54,39)(62,45)
\path(38,27)(29,31)
\path(54,39)(62,45)
\path(38,27)(39,13)
\path(68,27)(69,13)
\path(39,13)(27,19)
\path(41,45)(53,53)
\path(54,39)(41,45)
\path(56,19)(68,27)
\path(57,9)(58,4)
\path(57,9)(58,4)
\path(29,31)(25,33)
\path(62,45)(66,47)
\path(29,31)(25,33)
\path(62,45)(66,47)
\path(39,13)(40,6)
\path(69,13)(70,6)
\path(27,19)(20,22)
\path(53,53)(58,57)
\path(41,45)(35,47)
\path(68,27)(73,31)
\path(27,19)(28,9)
\path(77,19)(78,9)
\path(28,9)(19,13)
\path(44,57)(52,63)
\path(53,53)(44,57)
\path(69,13)(77,19)
\path(58,4)(58,2)
\path(58,4)(58,2)
\path(25,33)(22,34)
\path(66,47)(68,49)
\path(25,33)(22,34)
\path(66,47)(68,49)
\path(40,6)(40,3)
\path(70,6)(70,3)
\path(20,22)(17,23)
\path(58,57)(61,59)
\path(35,47)(32,49)
\path(73,31)(76,33)
\path(28,9)(28,4)
\path(78,9)(79,4)
\path(19,13)(14,15)
\path(52,63)(56,66)
\path(44,57)(39,59)
\path(77,19)(81,22)
\path(19,13)(19,6)
\path(84,13)(84,6)
\path(19,6)(13,9)
\path(45,66)(51,70)
\path(52,63)(45,66)
\path(78,9)(84,13)
\path(58,2)(58,1)
\path(58,2)(58,1)
\path(22,34)(21,35)
\path(68,49)(69,50)
\path(22,34)(21,35)
\path(68,49)(69,50)
\path(40,3)(41,1)
\path(70,3)(70,1)
\path(17,23)(16,24)
\path(61,59)(63,60)
\path(32,49)(30,50)
\path(76,33)(77,34)
\path(28,4)(28,2)
\path(79,4)(79,2)
\path(14,15)(12,16)
\path(56,66)(58,67)
\path(39,59)(37,60)
\path(81,22)(83,23)
\path(19,6)(20,3)
\path(84,6)(85,3)
\path(13,9)(10,11)
\path(51,70)(54,72)
\path(45,66)(42,67)
\path(84,13)(86,15)
\path(13,9)(13,4)
\path(88,9)(89,4)
\path(13,4)(9,6)
\path(47,72)(51,74)
\path(51,70)(47,72)
\path(84,6)(88,9)
\path(58,1)(58,0)
\path(58,1)(58,0)
\path(21,35)(21,35)
\path(69,50)(70,50)
\path(21,35)(21,35)
\path(69,50)(70,50)
\path(41,1)(41,0)
\path(70,1)(70,0)
\path(16,24)(15,24)
\path(63,60)(63,60)
\path(30,50)(30,50)
\path(77,34)(78,35)
\path(28,2)(28,1)
\path(79,2)(79,1)
\path(12,16)(11,17)
\path(58,67)(59,68)
\path(37,60)(36,60)
\path(83,23)(84,24)
\path(20,3)(20,1)
\path(85,3)(85,1)
\path(10,11)(8,11)
\path(54,72)(55,73)
\path(42,67)(41,68)
\path(86,15)(88,16)
\path(13,4)(14,2)
\path(89,4)(89,2)
\path(9,6)(7,7)
\path(51,74)(53,76)
\path(47,72)(44,73)
\path(88,9)(90,11)
\path(9,6)(9,3)
\path(92,6)(92,3)
\path(9,3)(6,4)
\path(48,76)(50,78)
\path(51,74)(48,76)
\path(89,4)(92,6)
\path(58,0)(58,0)
\path(58,0)(58,0)
\path(21,35)(20,35)
\path(70,50)(70,50)
\path(21,35)(20,35)
\path(70,50)(70,50)
\path(41,0)(41,0)
\path(70,0)(70,0)
\path(15,24)(14,25)
\path(63,60)(64,61)
\path(30,50)(29,50)
\path(78,35)(79,35)
\path(28,1)(29,0)
\path(79,1)(79,0)
\path(11,17)(10,17)
\path(59,68)(59,68)
\path(36,60)(35,61)
\path(84,24)(84,24)
\path(20,1)(20,0)
\path(85,1)(85,0)
\path(8,11)(7,12)
\path(55,73)(56,73)
\path(41,68)(40,68)
\path(88,16)(89,17)
\path(14,2)(14,1)
\path(89,2)(89,1)
\path(7,7)(6,8)
\path(53,76)(54,77)
\path(44,73)(43,73)
\path(90,11)(91,11)
\path(9,3)(9,1)
\path(92,3)(92,1)
\path(6,4)(5,5)
\path(50,78)(52,79)
\path(48,76)(46,77)
\path(92,6)(93,7)
\path(6,4)(6,2)
\path(94,4)(94,2)
\path(6,2)(4,3)
\path(48,79)(50,80)
\path(50,78)(48,79)
\path(92,3)(94,4)
\path(58,0)(58,0)
\path(58,0)(58,0)
\path(20,35)(20,35)
\path(70,50)(70,50)
\path(20,35)(20,35)
\path(70,50)(70,50)
\path(41,0)(41,0)
\path(70,0)(70,0)
\path(14,25)(14,25)
\path(64,61)(64,61)
\path(29,50)(29,50)
\path(79,35)(79,35)
\path(29,0)(29,0)
\path(79,0)(79,0)
\path(10,17)(10,17)
\path(59,68)(60,68)
\path(35,61)(35,61)
\path(84,24)(85,25)
\path(20,0)(20,0)
\path(85,0)(85,0)
\path(7,12)(7,12)
\path(56,73)(56,73)
\path(40,68)(40,68)
\path(89,17)(89,17)
\path(14,1)(14,0)
\path(89,1)(89,0)
\path(6,8)(5,8)
\path(54,77)(54,77)
\path(43,73)(43,73)
\path(91,11)(92,12)
\path(9,1)(10,0)
\path(92,1)(92,0)
\path(5,5)(4,5)
\path(52,79)(52,79)
\path(46,77)(45,77)
\path(93,7)(94,8)
\path(6,2)(7,1)
\path(94,2)(94,1)
\path(4,3)(3,3)
\path(50,80)(51,81)
\path(48,79)(47,79)
\path(94,4)(95,5)
\path(4,3)(4,1)
\path(96,3)(96,1)
\path(4,1)(3,2)
\path(49,81)(50,82)
\path(50,80)(49,81)
\path(94,2)(96,3)
\path(58,0)(58,0)
\path(58,0)(58,0)
\path(20,35)(20,35)
\path(70,50)(70,50)
\path(20,35)(20,35)
\path(70,50)(70,50)
\path(41,0)(41,0)
\path(70,0)(70,0)
\path(14,25)(14,25)
\path(64,61)(64,61)
\path(29,50)(29,50)
\path(79,35)(79,35)
\path(29,0)(29,0)
\path(79,0)(79,0)
\path(10,17)(10,17)
\path(60,68)(60,68)
\path(35,61)(35,61)
\path(85,25)(85,25)
\path(20,0)(20,0)
\path(85,0)(85,0)
\path(7,12)(7,12)
\path(56,73)(57,73)
\path(40,68)(39,68)
\path(89,17)(89,17)
\path(14,0)(14,0)
\path(89,0)(89,0)
\path(5,8)(5,8)
\path(54,77)(54,77)
\path(43,73)(43,73)
\path(92,12)(92,12)
\path(10,0)(10,0)
\path(92,0)(92,0)
\path(4,5)(3,6)
\path(52,79)(53,80)
\path(45,77)(45,77)
\path(94,8)(94,8)
\path(7,1)(7,0)
\path(94,1)(94,0)
\path(3,3)(3,4)
\path(51,81)(52,81)
\path(47,79)(46,80)
\path(95,5)(95,5)
\path(4,1)(4,0)
\path(96,1)(96,0)
\path(3,2)(2,2)
\path(50,82)(51,83)
\path(49,81)(48,81)
\path(96,3)(96,3)
\path(3,2)(3,1)
\path(97,2)(97,1)
\path(3,1)(2,1)
\path(49,83)(50,83)
\path(50,82)(49,83)
\path(96,1)(97,2)
\path(58,0)(58,0)
\path(58,0)(58,0)
\path(20,35)(20,35)
\path(70,50)(70,50)
\path(20,35)(20,35)
\path(70,50)(70,50)
\path(41,0)(41,0)
\path(70,0)(70,0)
\path(14,25)(14,25)
\path(64,61)(64,61)
\path(29,50)(29,50)
\path(79,35)(79,35)
\path(29,0)(29,0)
\path(79,0)(79,0)
\path(10,17)(10,17)
\path(60,68)(60,68)
\path(35,61)(35,61)
\path(85,25)(85,25)
\path(20,0)(20,0)
\path(85,0)(85,0)
\path(7,12)(7,12)
\path(57,73)(57,74)
\path(39,68)(39,68)
\path(89,17)(89,17)
\path(14,0)(14,0)
\path(89,0)(89,0)
\path(5,8)(5,8)
\path(54,77)(54,77)
\path(43,73)(42,74)
\path(92,12)(92,12)
\path(10,0)(10,0)
\path(92,0)(92,0)
\path(3,6)(3,6)
\path(53,80)(53,80)
\path(45,77)(45,77)
\path(94,8)(94,8)
\path(7,0)(7,0)
\path(94,0)(94,0)
\path(3,4)(2,4)
\path(52,81)(52,82)
\path(46,80)(46,80)
\path(95,5)(96,6)
\path(4,0)(4,0)
\path(96,0)(96,0)
\path(2,2)(2,2)
\path(51,83)(51,83)
\path(48,81)(47,82)
\path(96,3)(97,4)
\path(3,1)(3,0)
\path(97,1)(97,0)
\path(2,1)(1,1)
\path(50,83)(50,84)
\path(49,83)(48,83)
\path(97,2)(97,2)
\path(2,1)(2,0)
\path(98,1)(98,0)
\path(2,0)(1,1)
\path(49,84)(50,84)
\path(50,83)(49,84)
\path(97,1)(98,1)
\path(58,0)(58,0)
\path(58,0)(58,0)
\path(20,35)(20,35)
\path(70,50)(70,50)
\path(20,35)(20,35)
\path(70,50)(70,50)
\path(41,0)(41,0)
\path(70,0)(70,0)
\path(14,25)(14,25)
\path(64,61)(64,61)
\path(29,50)(29,50)
\path(79,35)(79,35)
\path(29,0)(29,0)
\path(79,0)(79,0)
\path(10,17)(10,17)
\path(60,68)(60,68)
\path(35,61)(35,61)
\path(85,25)(85,25)
\path(20,0)(20,0)
\path(85,0)(85,0)
\path(7,12)(7,12)
\path(57,74)(57,74)
\path(39,68)(39,68)
\path(89,17)(89,17)
\path(14,0)(14,0)
\path(89,0)(89,0)
\path(5,8)(5,8)
\path(54,77)(55,77)
\path(42,74)(42,74)
\path(92,12)(92,12)
\path(10,0)(10,0)
\path(92,0)(92,0)
\path(3,6)(3,6)
\path(53,80)(53,80)
\path(45,77)(45,77)
\path(94,8)(94,8)
\path(7,0)(7,0)
\path(94,0)(94,0)
\path(2,4)(2,4)
\path(52,82)(52,82)
\path(46,80)(46,80)
\path(96,6)(96,6)
\path(4,0)(5,0)
\path(96,0)(96,0)
\path(2,2)(1,2)
\path(51,83)(51,83)
\path(47,82)(47,82)
\path(97,4)(97,4)
\path(3,0)(3,0)
\path(97,0)(97,0)
\path(1,1)(1,2)
\path(50,84)(51,84)
\path(48,83)(48,83)
\path(97,2)(97,2)
\path(2,0)(2,0)
\path(98,0)(98,0)
\path(1,1)(1,1)
\path(50,84)(50,84)
\path(49,84)(49,84)
\path(98,1)(98,1)
\path(1,1)(1,0)
\path(98,1)(98,0)
\path(1,0)(1,0)
\path(49,84)(50,85)
\path(50,84)(49,84)
\path(98,0)(98,1)
\path(54,39)(38,27)
\path(54,39)(68,27)
\path(56,19)(39,13)
\path(38,27)(41,45)
\path(56,19)(54,39)
\path(38,27)(56,19)
\path(56,19)(39,13)
\path(56,19)(69,13)
\path(38,27)(27,19)
\path(54,39)(53,53)
\path(38,27)(41,45)
\path(54,39)(68,27)
\path(38,27)(27,19)
\path(68,27)(77,19)
\path(39,13)(28,9)
\path(41,45)(44,57)
\path(54,39)(53,53)
\path(56,19)(69,13)
\path(57,9)(40,6)
\path(57,9)(70,6)
\path(29,31)(20,22)
\path(62,45)(58,57)
\path(29,31)(35,47)
\path(62,45)(73,31)
\path(39,13)(28,9)
\path(69,13)(78,9)
\path(27,19)(19,13)
\path(53,53)(52,63)
\path(41,45)(44,57)
\path(68,27)(77,19)
\path(27,19)(19,13)
\path(77,19)(84,13)
\path(28,9)(19,6)
\path(44,57)(45,66)
\path(53,53)(52,63)
\path(69,13)(78,9)
\path(58,4)(40,3)
\path(58,4)(70,3)
\path(25,33)(17,23)
\path(66,47)(61,59)
\path(25,33)(32,49)
\path(66,47)(76,33)
\path(40,6)(28,4)
\path(70,6)(79,4)
\path(20,22)(14,15)
\path(58,57)(56,66)
\path(35,47)(39,59)
\path(73,31)(81,22)
\path(28,9)(19,6)
\path(78,9)(84,6)
\path(19,13)(13,9)
\path(52,63)(51,70)
\path(44,57)(45,66)
\path(77,19)(84,13)
\path(19,13)(13,9)
\path(84,13)(88,9)
\path(19,6)(13,4)
\path(45,66)(47,72)
\path(52,63)(51,70)
\path(78,9)(84,6)
\path(58,2)(41,1)
\path(58,2)(70,1)
\path(22,34)(16,24)
\path(68,49)(63,60)
\path(22,34)(30,50)
\path(68,49)(77,34)
\path(40,3)(28,2)
\path(70,3)(79,2)
\path(17,23)(12,16)
\path(61,59)(58,67)
\path(32,49)(37,60)
\path(76,33)(83,23)
\path(28,4)(20,3)
\path(79,4)(85,3)
\path(14,15)(10,11)
\path(56,66)(54,72)
\path(39,59)(42,67)
\path(81,22)(86,15)
\path(19,6)(13,4)
\path(84,6)(89,4)
\path(13,9)(9,6)
\path(51,70)(51,74)
\path(45,66)(47,72)
\path(84,13)(88,9)
\path(13,9)(9,6)
\path(88,9)(92,6)
\path(13,4)(9,3)
\path(47,72)(48,76)
\path(51,70)(51,74)
\path(84,6)(89,4)
\path(58,1)(41,0)
\path(58,1)(70,0)
\path(21,35)(15,24)
\path(69,50)(63,60)
\path(21,35)(30,50)
\path(69,50)(78,35)
\path(41,1)(28,1)
\path(70,1)(79,1)
\path(16,24)(11,17)
\path(63,60)(59,68)
\path(30,50)(36,60)
\path(77,34)(84,24)
\path(28,2)(20,1)
\path(79,2)(85,1)
\path(12,16)(8,11)
\path(58,67)(55,73)
\path(37,60)(41,68)
\path(83,23)(88,16)
\path(20,3)(14,2)
\path(85,3)(89,2)
\path(10,11)(7,7)
\path(54,72)(53,76)
\path(42,67)(44,73)
\path(86,15)(90,11)
\path(13,4)(9,3)
\path(89,4)(92,3)
\path(9,6)(6,4)
\path(51,74)(50,78)
\path(47,72)(48,76)
\path(88,9)(92,6)
\path(9,6)(6,4)
\path(92,6)(94,4)
\path(9,3)(6,2)
\path(48,76)(48,79)
\path(51,74)(50,78)
\path(89,4)(92,3)
\path(58,0)(41,0)
\path(58,0)(70,0)
\path(21,35)(14,25)
\path(70,50)(64,61)
\path(21,35)(29,50)
\path(70,50)(79,35)
\path(41,0)(29,0)
\path(70,0)(79,0)
\path(15,24)(10,17)
\path(63,60)(59,68)
\path(30,50)(35,61)
\path(78,35)(84,24)
\path(28,1)(20,0)
\path(79,1)(85,0)
\path(11,17)(7,12)
\path(59,68)(56,73)
\path(36,60)(40,68)
\path(84,24)(89,17)
\path(20,1)(14,1)
\path(85,1)(89,1)
\path(8,11)(6,8)
\path(55,73)(54,77)
\path(41,68)(43,73)
\path(88,16)(91,11)
\path(14,2)(9,1)
\path(89,2)(92,1)
\path(7,7)(5,5)
\path(53,76)(52,79)
\path(44,73)(46,77)
\path(90,11)(93,7)
\path(9,3)(6,2)
\path(92,3)(94,2)
\path(6,4)(4,3)
\path(50,78)(50,80)
\path(48,76)(48,79)
\path(92,6)(94,4)
\path(6,4)(4,3)
\path(94,4)(96,3)
\path(6,2)(4,1)
\path(48,79)(49,81)
\path(50,78)(50,80)
\path(92,3)(94,2)
\path(58,0)(41,0)
\path(58,0)(70,0)
\path(20,35)(14,25)
\path(70,50)(64,61)
\path(20,35)(29,50)
\path(70,50)(79,35)
\path(41,0)(29,0)
\path(70,0)(79,0)
\path(14,25)(10,17)
\path(64,61)(60,68)
\path(29,50)(35,61)
\path(79,35)(85,25)
\path(29,0)(20,0)
\path(79,0)(85,0)
\path(10,17)(7,12)
\path(59,68)(56,73)
\path(35,61)(40,68)
\path(84,24)(89,17)
\path(20,0)(14,0)
\path(85,0)(89,0)
\path(7,12)(5,8)
\path(56,73)(54,77)
\path(40,68)(43,73)
\path(89,17)(92,12)
\path(14,1)(10,0)
\path(89,1)(92,0)
\path(6,8)(4,5)
\path(54,77)(52,79)
\path(43,73)(45,77)
\path(91,11)(94,8)
\path(9,1)(7,1)
\path(92,1)(94,1)
\path(5,5)(3,3)
\path(52,79)(51,81)
\path(46,77)(47,79)
\path(93,7)(95,5)
\path(6,2)(4,1)
\path(94,2)(96,1)
\path(4,3)(3,2)
\path(50,80)(50,82)
\path(48,79)(49,81)
\path(94,4)(96,3)
\path(4,3)(3,2)
\path(96,3)(97,2)
\path(4,1)(3,1)
\path(49,81)(49,83)
\path(50,80)(50,82)
\path(94,2)(96,1)
\path(58,0)(41,0)
\path(58,0)(70,0)
\path(20,35)(14,25)
\path(70,50)(64,61)
\path(20,35)(29,50)
\path(70,50)(79,35)
\path(41,0)(29,0)
\path(70,0)(79,0)
\path(14,25)(10,17)
\path(64,61)(60,68)
\path(29,50)(35,61)
\path(79,35)(85,25)
\path(29,0)(20,0)
\path(79,0)(85,0)
\path(10,17)(7,12)
\path(60,68)(57,73)
\path(35,61)(39,68)
\path(85,25)(89,17)
\path(20,0)(14,0)
\path(85,0)(89,0)
\path(7,12)(5,8)
\path(56,73)(54,77)
\path(40,68)(43,73)
\path(89,17)(92,12)
\path(14,0)(10,0)
\path(89,0)(92,0)
\path(5,8)(3,6)
\path(54,77)(53,80)
\path(43,73)(45,77)
\path(92,12)(94,8)
\path(10,0)(7,0)
\path(92,0)(94,0)
\path(4,5)(3,4)
\path(52,79)(52,81)
\path(45,77)(46,80)
\path(94,8)(95,5)
\path(7,1)(4,0)
\path(94,1)(96,0)
\path(3,3)(2,2)
\path(51,81)(51,83)
\path(47,79)(48,81)
\path(95,5)(96,3)
\path(4,1)(3,1)
\path(96,1)(97,1)
\path(3,2)(2,1)
\path(50,82)(50,83)
\path(49,81)(49,83)
\path(96,3)(97,2)
\path(3,2)(2,1)
\path(97,2)(98,1)
\path(3,1)(2,0)
\path(49,83)(49,84)
\path(50,82)(50,83)
\path(96,1)(97,1)
\path(58,0)(41,0)
\path(58,0)(70,0)
\path(20,35)(14,25)
\path(70,50)(64,61)
\path(20,35)(29,50)
\path(70,50)(79,35)
\path(41,0)(29,0)
\path(70,0)(79,0)
\path(14,25)(10,17)
\path(64,61)(60,68)
\path(29,50)(35,61)
\path(79,35)(85,25)
\path(29,0)(20,0)
\path(79,0)(85,0)
\path(10,17)(7,12)
\path(60,68)(57,74)
\path(35,61)(39,68)
\path(85,25)(89,17)
\path(20,0)(14,0)
\path(85,0)(89,0)
\path(7,12)(5,8)
\path(57,73)(54,77)
\path(39,68)(42,74)
\path(89,17)(92,12)
\path(14,0)(10,0)
\path(89,0)(92,0)
\path(5,8)(3,6)
\path(54,77)(53,80)
\path(43,73)(45,77)
\path(92,12)(94,8)
\path(10,0)(7,0)
\path(92,0)(94,0)
\path(3,6)(2,4)
\path(53,80)(52,82)
\path(45,77)(46,80)
\path(94,8)(96,6)
\path(7,0)(4,0)
\path(94,0)(96,0)
\path(3,4)(2,2)
\path(52,81)(51,83)
\path(46,80)(47,82)
\path(95,5)(97,4)
\path(4,0)(3,0)
\path(96,0)(97,0)
\path(2,2)(1,1)
\path(51,83)(50,84)
\path(48,81)(48,83)
\path(96,3)(97,2)
\path(3,1)(2,0)
\path(97,1)(98,0)
\path(2,1)(1,1)
\path(50,83)(50,84)
\path(49,83)(49,84)
\path(97,2)(98,1)
\path(2,1)(1,1)
\path(98,1)(98,1)
\path(2,0)(1,0)
\path(49,84)(49,84)
\path(50,83)(50,84)
\path(97,1)(98,0)
\path(58,0)(41,0)
\path(58,0)(70,0)
\path(20,35)(14,25)
\path(70,50)(64,61)
\path(20,35)(29,50)
\path(70,50)(79,35)
\path(41,0)(29,0)
\path(70,0)(79,0)
\path(14,25)(10,17)
\path(64,61)(60,68)
\path(29,50)(35,61)
\path(79,35)(85,25)
\path(29,0)(20,0)
\path(79,0)(85,0)
\path(10,17)(7,12)
\path(60,68)(57,74)
\path(35,61)(39,68)
\path(85,25)(89,17)
\path(20,0)(14,0)
\path(85,0)(89,0)
\path(7,12)(5,8)
\path(57,74)(55,77)
\path(39,68)(42,74)
\path(89,17)(92,12)
\path(14,0)(10,0)
\path(89,0)(92,0)
\path(5,8)(3,6)
\path(54,77)(53,80)
\path(42,74)(45,77)
\path(92,12)(94,8)
\path(10,0)(7,0)
\path(92,0)(94,0)
\path(3,6)(2,4)
\path(53,80)(52,82)
\path(45,77)(46,80)
\path(94,8)(96,6)
\path(7,0)(5,0)
\path(94,0)(96,0)
\path(2,4)(1,2)
\path(52,82)(51,83)
\path(46,80)(47,82)
\path(96,6)(97,4)
\path(4,0)(3,0)
\path(96,0)(97,0)
\path(2,2)(1,2)
\path(51,83)(51,84)
\path(47,82)(48,83)
\path(97,4)(97,2)
\path(3,0)(2,0)
\path(97,0)(98,0)
\path(1,1)(1,1)
\path(50,84)(50,84)
\path(48,83)(49,84)
\path(97,2)(98,1)
\path(2,0)(1,0)
\path(98,0)(98,0)
\path(1,1)(1,0)
\path(50,84)(50,85)
\path(49,84)(49,84)
\path(98,1)(98,1)
\path(1,1)(1,0)
\path(98,1)(99,0)
\path(1,0)(1,0)
\path(49,84)(49,85)
\path(50,84)(50,85)
\path(98,0)(98,0)
\end{picture}
\end{center}
\end{figure}

\begin{abstract}
L'un des buts de cet article est de décrire l'isomorphisme entre les
tours de Lubin-Tate et de Drinfeld au niveau de leurs squelettes
après quotient par
$\GL_n (\O_F)\times \O_D^\times$ ou bien $I\times \O_D^\times$ où
$\O_D$ est l'ordre maximal dans l'algèbre à division d'invariant
$\frac{1}{n}$ sur $F$ et $I$ un sous-groupe d'Iwahori de $\GL_n$. Nous 
donnons des applications à l'étude des sous-groupes canoniques sur les espaces
de Lubin-Tate, 
 la description
des orbites de Hecke sphériques dans ces espaces, les
domaines fondamentaux pour les correspondances de Hecke et
l'application des périodes de Gross-Hopkins. Nous-y étudions également
en détail les filtrations de ramification (inférieure et supérieure) et l'application de
Hodge-Tate d'un groupe formel $p$-divisible de dimension un. 
\end{abstract}

\begin{altabstract}
One of the goals of this article  is to describe the isomorphism
between Lubin-Tate and Drinfeld towers at the level of their skeletons
after taking quotient by  $\GL_n (\O_F)\times \O_D^\times$ or $I\times
\O_D^\times$ where $\O_D$ is the maximal order in the division algebra
with invariant $\frac{1}{n}$ over $F$ and $I$ a Iwahori subgroup of
$\GL_n$. We give applications to the theory of canonical subgroups on
Lubin-Tate spaces, the description of spherical Hecke orbits in those
spaces, fundamental domains for Hecke correspondences and the
Gross-Hopkins period mapping. We also study in details the
ramification filtrations (upper and lower) and the Hodge-Tate map of a one dimensional
formal $p$-divisible group. 
\end{altabstract}

\section*{Introduction}

L'un des buts de cet article est de décrire l'isomorphisme entre les
squelettes des tours de Lubin-Tate et de Drinfeld après quotient par
$\GL_n (\O_F)\times \O_D^\times$ ou bien $I\times \O_D^\times$ où
$\O_D$ est l'ordre maximal dans l'algèbre à division d'invariant
$\frac{1}{n}$ et $I$ un sous-groupe d'Iwahori de $\GL_n$. 

L'isomorphisme entre les tours de Lubin-Tate et de Drinfeld au niveau
des points (\cite{Points}) 
 est un
isomorphisme $\GL_n (F)\times D^\times$-équivariant en niveau infini 
$$\xymatrix@C=15mm{
\LT_\infty \ar[r]^\sim \ar[d] &
\D_\infty \ar[d] \\
\LT_\infty / \GL_n (\O_F) =\mathring{\mathbb{B}}^{n-1}  &
\D_\infty/\O_D^\times =  \Omega
}
$$ 
 qui induit  un homéomorphisme des espaces de Berkovich
associés  
$$
|\LT_\infty |\iso |\D_\infty|
$$
Par passage au quotient par $\GL_n (\O_F)\times \O_D^\times$ il
devrait donc induire une application 
$$
\O_D^\times \bc |\mathring{\mathbb{B}}^{n-1}| \ldrt \GL_n (F)\bc |\Omega |
$$
On  va décrire  cette 
application au niveau des squelettes de ces deux espaces.
\\

\begin{figure}[htbp]
   \begin{center}
      \includegraphics{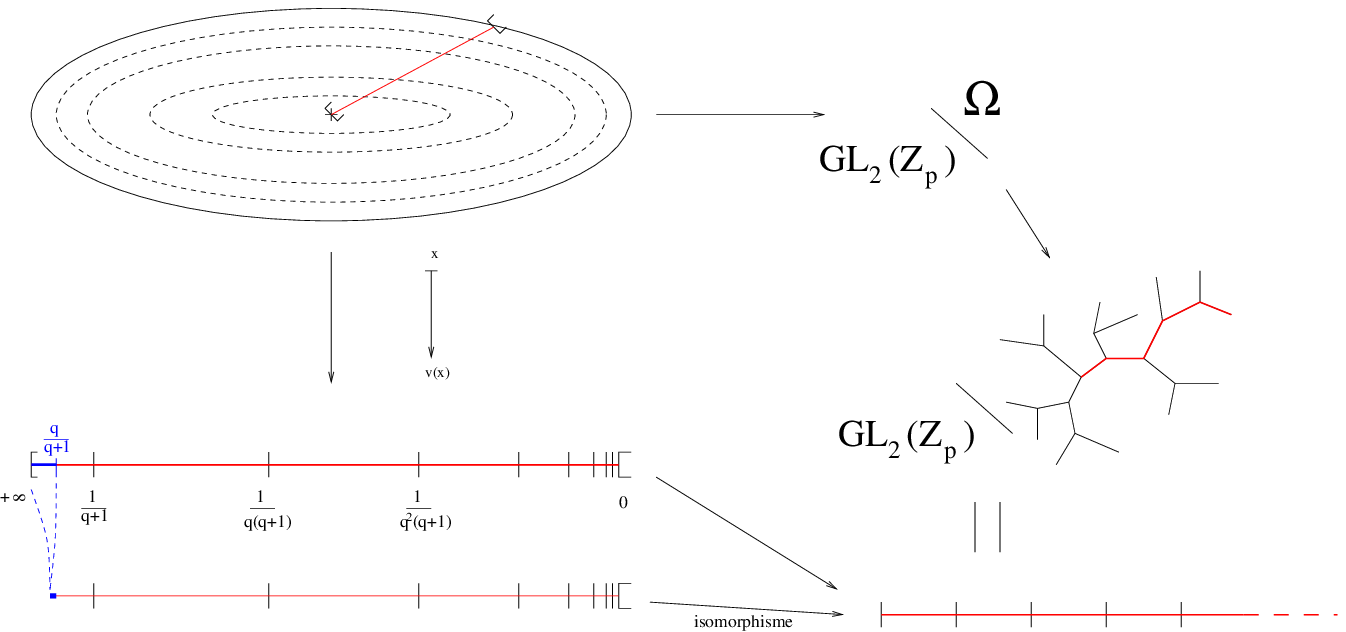}
   \end{center}
   \caption{\footnotesize Le cas de $\GL_2$
     }
   \label{ujikolpo}
\end{figure}

Plutôt que de tenter de décrire le résultat général explicitons ce que
cela signifie sur la figure \ref{ujikolpo}
dans le cas de $\GL_2$ pour le quotient par $\GL_2
(\Zp)\times \O_D^\times$ : 
\begin{itemize}
\item
L'espace de Lubin-Tate sans niveau (la tour de Lubin-Tate quotientée
par $\GL_2 (\Zp)$) est une boule ouverte $p$-adique au sens de
Berkovich $\mathring{\mathbb{B}^1}$. Dans ce cas là appelons
squelette de $\mathring{\mathbb{B}^1}$ un rayon de cette boule
$]0,+\infty]$. Il y a une rétraction (la flèche verticale de gauche) 
$|\mathring{\mathbb{B}^1} |\ldrt ]0,+\infty]$ donnée par la
valuation de la coordonnée dans la boule.
\item L'espace de Drinfeld sans niveau (la tour de Drinfeld après
  quotient par $\O_D^\times$) est l'espace $\Omega$ de Drinfeld ayant
  pour $\mathbb{C}_p$-points $\C_p\setminus \Qp$. Son squelette est
  l'arbre de Bruhat-Tits $\mathcal{I}$ de $\GL_2$. Il y a une
  rétraction (la flèche verticale de droite) $|\Omega|\ldrt
  \mathcal{I}$ qui après quotient par $\GL_2 (\Zp)$ fournit une
  rétraction $\GL_2 (\Zp) \bc |\Omega |\ldrt \GL_2 (\Zp) \bc \mathcal{I}$. 
\item Si $D$ désigne une demi-droite simpliciale d'origine la classe
  du réseau $\Zp^2$ dans l'arbre $\mathcal{I}$ alors $D$ est un
  domaine fondamental pour l'action de $\GL_2 (\Zp)$ sur $\mathcal{I}$
  et 
$D\iso \GL_2 (\Zp)\bc \mathcal{I}$.
\item L'isomorphisme entre les deux tours induit une application
  $]0,+\infty ]\ldrt D$
\item On peut décrire complètement la structure simpliciale 
sur $]0,+\infty]$ déduite de
  celle sur $D$ par l'application précédente.
\item Après quotient par une ``petite'' partie de $]0,+\infty]$
  l'application induit un isomorphisme (flèche du bas). 
\end{itemize}

Le résultat est du même type pour  $\GL_n$, bien qu'un peu plus compliqué à
énoncer. 

 Indépendamment de l'isomorphisme
entre les deux tours la structure simpliciale 
que nous explicitons sur l'espace de
Lubin-Tate a de nombreuses applications  comme l'étude des
sous-groupes canoniques, la détermination de domaines fondamentaux
pour les correspondances de Hecke et l'étude du morphisme des
périodes. 
\\
L'un des autres buts de cet article est d'étudier en détails la filtration donnée par la valuation des points de torsion sur un groupe formel $p$-divisible de dimension un.
\\

Décrivons succinctement le contenu de chacune des parties de l'article :
\begin{itemize}
\item Dans le premier chapitre nous donnons une formule pour la valuation $p$-adique des périodes de Hodge-Tate du dual de Cartier 
d'un groupe $p$-divisible formel sur un anneau de valuation (pas
forcément discrète) pour une valuation de hauteur $1$. 
En fait, nous considérons plus généralement le cas d'un $\O$-module formel $\pi$-divisible où $\O$ est l'anneau des entiers d'une extension de degré fini de $\Qp$. Dans ce cas la bonne notion de dualité remplaçant la dualité de Cartier est celle définie par Faltings (\cite{Faltings7}). Le lecteur ne connaissant pas la théorie de \cite{Faltings7} pourra supposer $\O=\Zp$.
\item Dans le second chapitre on étudie la filtration donnée par la
  valuation sur les points de torsion d'un group formel
  $\pi$-divisible de dimension $1$. 
Cette filtration fournit une famille de réseaux dans le module de Tate rationnel $V_p$. L'image dans l'immeuble de $\text{PGL} (V_p)$ de cet ensemble est un ensemble fini de sommets. De plus les éléments de valuation suffisamment petite décrivent un simplexe $S$ de cet immeuble. L'un des principaux résultats est que cet ensemble est contenu dans un appartement et peut être reconstruit géométriquement dans l'immeuble à partir du simplexe $S$ et du sommet donné par le réseau $T_p\subset V_p$.
\\
Nous donnons  également une description combinatoire du simplexe $S$ à partir du polygone de Newton de la multiplication par $\pi$ sur une loi de groupe formel associée. 
\item Le troisième chapitre est inspiré par les travaux d'Abbes-Saito et Abbes-Mokrane (\cite{AbbesMokrane}). La filtration sur les points de torsion étudiée dans le deuxième chapitre se comporte bien par restriction à un sous-groupe : si $H$ est un groupe formel $p$-divisible de dimension $1$ et $G_2 \subset G_1\subset H$ des sous-groupes plats finis alors
$
\forall \l\; \{x\in G_1\;|\; v(x)\geq \l \}\cap G_2 = \{ x\in G_2\;|\; v(x)\geq \l \}
$. Par contre 
cette filtration dite de ``ramification inférieure''  
 ne se comporte pas bien par isogénies. C'est le cas de la filtration définie en toute généralités dans \cite{AbbesMokrane}. Dans le cas que nous étudions des sous-groupes plats finis d'un groupe formel $p$-divisible de dimension $1$ l'algèbre de ces groupes est monogène et la filtration de ramification supérieure de \cite{AbbesMokrane} est obtenue par réindexation de la filtration de ramification inférieure via une fonction de Herbrand. Cela est expliqué dans l'appendice \ref{Ramitofm}. La terminologie ``inférieure/supérieur'' provient par analogie avec la théorie des groupes de ramification des groupes de Galois des corps locaux : les groupes de ramification inférieure se comportent bien par restriction à un sous-groupe de Galois tandis que ceux de ramification supérieure se comportent bien vis à vis d'un quotient. 

Nous étudions cette filtration de ramification supérieure ainsi que son image dans l'immeuble de la même façon que dans le chapitre deux. 
\item  Dans le quatrième chapitre on étudie le point de la réalisation géométrique de l'immeuble défini par l'application de Hodge-Tate du dual d'un $\O$-module formel $\pi$-divisible de dimension un. Ce point est la classe d'équivalence de la norme sur le module de Tate rationnel donnée par la valuation de l'application de Hodge-Tate étudiée dans le premier chapitre. On donne des formules intégrales pour cette norme en fonction des filtrations étudiées aux chapitres 2 et 3. Cette formule est particulièrement simple lorsque formulée en termes de la filtration de ramification supérieure (proposition \ref{kujipopom}). 
\\
L'un des principaux corollaires de ces formules est que ce point dans l'immeuble est contenu dans la réalisation géométrique $|S|$ du simplexe $S$ défini au chapitre 2. 
\item 
Dans le chapitre 5 on définit et étudie une structure simpliciale sur le squelette de l'espace de Lubin-Tate sans niveau. Le bon objet n'est pas en fait ce squelette mais plutôt un quotient de celui-ci, l'espace des polygones de Newton. On décrit complètement une structure simpliciale sur cet espace des polygones de Newton ainsi que l'action de certains opérateurs de Hecke sur cet espace simplicial. 

La définition de cette structure simpliciale est inspirée des résultats des chapitres 2 et 4. 
\item Dans le chapitre 6 on montre que la bijection entre les points des tours de Lubin-Tate et de Drinfeld induit un isomorphisme entre l'espace des polygones de Newton muni de la structure simpliciale définie au chapitre 5 et le quotient de l'immeuble de $\text{PGL}_n$ par un sous-groupe compact maximal.
\item Le chapitre 7 est consacré aux applications de la structure simpliciale sur l'espace des polygones de Newton et de l'action des opérateurs de Hecke sur celle-ci. 
Certains raisonnements sur l'espace de Lubin-Tate s'interprètent naturellement sur un appartement de l'immeuble. 

Par exemple on démontre que l'existence de sous-groupes canoniques en un sens généralisé est équivalent à ce que le point dans l'immeuble soit contenu dans un certain demi-appartement.
Cela démontre par un simple raisonnement géométrique que le ``bord'' de l'espace de Lubin-Tate est recouvert par des ouverts admissibles sur lesquels il ya des sous-groupes canoniques puisque c'est le cas dans l'immeuble. L'application quotient par un sous-groupe canonique se comprend également très facilement sur l'immeuble, de même que l'action des correspondances de Hecke.
\\
Couplé aux résultats des chapitres précédents cela donne une condition nécessaire et suffisante pour l'existence de sous-groupes canoniques généralisés en termes de l'application de Hodge-Tate du groupe $p$-divisible formel de dimension $1$, comme dans \cite{AbbesMokrane}. 

On généralise également le domaine fondamental de Gross-Hopkins grâce à cette étude sur l'immeuble : n'importe quel domaine fondamental polyèdral dans le simplexe standard $\text{Convexe} (e_0,\dots, e_n)\subset \R^{n+1}$ sous l'action du groupe des rotations engendrées par $e_0\mapsto e_1,e_1\mapsto e_2,\dots, e_n\mapsto e_0$ fournit un  domaine fondamental 
pour l'action des correspondances de Hecke dans l'espace de Lubin-Tate. 
\item Dans le chapitre 8 on généralise les résultats précédents au cas
  de l'espace de Lubin-Tate avec structure de niveau Iwahori. Dans ce
  cas là l'espace est une couronne $p$-adique généralisée et son squelette un simplexe ``ouvert''. On définit et étudie alors comme auparavant une structure simpliciale sur ce simplexe et montre que 
via l'isomorphisme entre les tours de Lubin-Tate et de Drinfeld 
ce simplexe ouvert est isomorphe au quotient de l'immeuble par un sous-groupe d'Iwahori.
Enfin on peut comprendre facilement grâce à cette étude le morphisme de l'espace de Lubin-Tate avec niveau Iwahori vers celui sans niveau au niveau des squelettes. 
Ce chapitre ne contient aucune démonstration, les démonstrations étant semblables à celles du cas de l'espace de Lubin-Tate sans niveau elles sont  laissées au lecteur. 
\item Enfin l'appendice A contient des rappels sur l'immeuble de Bruhat-Tits de  $\text{PGL}_n$.
\end{itemize}

\vspace{5mm}
Cet article est indépendant de \cite{Cellulaire} et \cite{iso4}. Seule
la compréhension d'une partie de \cite{Points} peut être utile. 
\\

Enfin, certains des aspects de cet article apparaissent déja dans les
travaux de Yu
\cite{Yu}. Cet article peut donc être en quelques sortes  considéré
comme une suite de \cite{Yu}, suite qui permet de comprendre pourquoi les
calculs effectués dans \cite{Yu} font appraître l'appartement d'un
immeuble de Bruhat-Tits.
\\

{\it Remerciements : L'auteur tient à remercier Alain Genestier et Vincent Lafforgue pour de nombreuses discussions sur le sujet. C'est en particulier Alain Genestier qui a suggéré 
d'introduire le simplexe  de la définition \ref{quatoki}, simplexe qui a suggéré à l'auteur d'étudier plus en détails les filtrations de ramification. 
Ils ont également suggéré à l'auteur l'étude du cas Iwahori faite au chapitre 8. 
}

\section{Une formule pour la valuation $p$-adique de l'application de
  Hodge-Tate du dual d'un groupe de Lubin-Tate}

Soit $F|\Qp$ une extension de degré fini et $\O=\O_F$ son anneau des
entiers. Soit $K|F$ un corps valué complet pour une valuation $v$ à
valeurs dans $\R$ étendant celle de $F$. Soit $H$ un $\O$-module
formel de dimension $1$ et de hauteur finie $n$ sur $\O_K$.

Supposons d'abord que $\O=\Zp$. Le but de cette section est de donner
une formule pour l'application composée 
$$
T_p (H^D) \xrig{\; \a_{H^D}\;} \omega_H\otimes
\O_{\widehat{\overline{K}}} \simeq \O_{\widehat{\overline{K}}}
\xrig{\; v \;} \R_+ \cup \{\infty \}
$$
où $\a_{H^D}$ est l'application de Hodge-Tate de $H^D$ : si $x\in T_p
(H^D)$, $x: \Qp/\Zp \ldrt H^D_{\O_{\widehat{\overline{K}}}}$ et $x^D:
H_{\O_{\widehat{\overline{K}}}} \ldrt \mu_{p^\infty/\O_{\widehat{\overline{K}}}}$ alors 
$$
\a_{H^D} (x) = (x^D)^*\frac{dT}{T}
$$

Pour $\O$ plus général que $\Zp$ nous donnons une formule pour la
composée 
$$
T_p (H^\vee) \xrig{\;\a_{H^\vee}^\O\;}  \omega_H \otimes
\O_{\widehat{\overline{K}}} \simeq \O_{\widehat{\overline{K}}}
\xrig{\; v \; } \R_+ \cup \{\infty \}
$$
où $H^\vee$ est le dual strict au sens de Faltings (\cite{Faltings7}). Si
$\LT$ désigne un groupe de Lubin-Tate de $\O$-hauteur $1$ alors à
$x\in T_p (H^\vee)$ est associé un morphisme $x^\vee  :
H_{\O_{\widehat{\overline{K}}}} \ldrt
\LT_{/\O_{\widehat{\overline{K}}}}$ qui définit donc $\a_{H^\vee}^\O (x) =
(x^\vee)^* \beta$ après choix d'un générateur $\beta$ de $\omega_{\LT}$.
\\

Il est clair que pour le problème auquel on s'intéresse on peut
supposer que $K=\widehat{\overline{K}}$, ce que nous ferons dans la suite.

\subsection{Périodes de Hodge-Tate de certains schémas en groupes de
  type $(p,\dots,p)$}
\subsubsection{Le cas $\O=\Zp$}

Soit $G$ un schéma en groupes fini localement libre d'ordre $p$ sur une base
affine $\spec (R)$ au dessus de $\spec(\Zp)$. 
D'après  \cite{OortTate} ou plus généralement \cite{Ray1} il
 existe alors $\gamma,\delta\in R$ tels que
$\gamma \delta = w$ où $w\in \Zp$ est une constante universelle de
valuation $p$-adique $1$ tels que 
\begin{eqnarray*}
G &\simeq & \spec (R[T]/(T^p-\delta T)) \\
G^D &\simeq & \spec (R[U]/(U^p-\gamma U))  
\end{eqnarray*}
Alors, 
\begin{eqnarray*}
\omega_G &\simeq &R/\delta R .dT \\
\omega_{G^D} &\simeq & R/\gamma R.dU 
\end{eqnarray*}
et 
\begin{eqnarray*}
\a_{G^D} : G^D & \ldrt & \omega_{G} \\
u & \longmapsto & (u \text{ mod } \delta ) . dT
\end{eqnarray*}
Si $R = \O_K$ avec $K$ comme précédemment alors $v(\gamma )+ v(\delta)
=1$ et 
$$
\forall u\in G^D(\O_K) \;\; v(u) = \frac{v(\gamma)}{p-1}
$$
et donc on connaît le sous-module $\O_K.\text{Im}\, \a_{G^D}$
de $\omega_{G^D}$ dès que l'on connait
$v(\text{Ann } \omega_G )= v (\delta)$ ou bien $v(\text{Ann }
\omega_{G^D})= v (\gamma)$.

\subsubsection{Le cas  $\O$  général}

Soit $\pi$ une uniformisante de $\O$ et $q=p^r$ le cardinal de son
corps résiduel. On note $v$ la valuation normalisée de $F$.

Soit $R$ une $\O$-algèbre. Soit $G$ un schéma en groupes fini et localement libre sur $\spec
(R)$. Supposons le muni d'une action de $\O/\pi\O$ et de type
$(p,\dots,p)$ relativement à cette action. L'anneau $R$ étant une
$\O$-algèbre il y a un caractère 
$$
\chi : (\O/\pi\O)^\times \xrig{\text{Teichmüller}} \O^\times \ldrt R^\times
$$
Alors, d'après \cite{Ray1} il existe  $(\gamma_i,\delta_i)_{i\in
  \Z/r\Z}\in R^{\Z/r\Z}$ tels que $\gamma_i\delta_i =w\in \Zp$ est de valuation
$p$-adique $1$, localement sur $\spec (R)\;\;$ $G\simeq \spec (A)$ avec
$$
A= R[T_i]_{i\in \Z/r\Z} / (T_i^p - \delta_i T_{i+1})
$$
et l'action de $(\O/\pi\O)^\times$ sur $A$ induite par l'action de
$\O/\pi\O$ sur $G$ se fasse sur $T_i$ à travers le caractère $\chi^{p^i}$. 
On a alors
$$
\omega_G \simeq \bigoplus_{i\in \Z/r\Z} R/\delta_{i-1} R . dT_i
$$
Supposons maintenant de plus que l'action de $\O$ sur $\omega_G$ soit
l'action naturelle induite par la structure de $\O$-algèbre de
$R$. Alors,
$$
\forall i\neq r-1 \;\;\delta_i\in R^\times
$$
et donc, si 
$$
\delta= \delta_{0}^{\, p^{r-1}}\delta_1^{\, p^{r-2}}\dots \delta_{r-2}^{\, p}
\delta_{r-1}
$$
on a 
$$
A \simeq R[T]/(T^q-\delta T)
$$
Le complexe de co-Lie de $G$ s'identifie alors à 
$$
l_G \simeq  [ R \xrig{\; \times \delta \; } R  ] 
$$
Supposons maintenant que $G$ est muni d'une action stricte de $\O$ au
sens de \cite{Faltings7} relevant l'action de $\O/\pi \O$ sur $G$. 
D'après \cite{Faltings7} l'ensemble de ces relèvements est un torseur
sous $H^{-1} (\End (l_G)) \simeq Ann_R (\delta)$. Supposons 
maintenant que $R$ est sans $p$-torsion. Cela implique $ Ann_R (\delta)=(0)$. 
 D'après ce qui précède il
existe donc une unique telle $\O$-action stricte : c'est celle définie
dans le chapitre 3 de \cite{Faltings7} sur le groupe noté $G_{u,v}$. 
On a donc identifié $G$ muni de son action stricte de $\O$ et d'après
le chapitre 3 de \cite{Faltings7} le dual strict d'un tel groupe est
connu. Rappelons en effet qu'alors il existe $\gamma \in R$ tel que 
$\gamma \delta =w'\in \O$ où $w'$ est une uniformisante de $F$ et que 
le dual stricte s'identifie à 
$$
G^\vee \simeq \spec (R[U]/(U^q-\gamma U))
$$
Soit $\LT$ un groupe formel de Lubin-Tate de $\O$-hauteur $1$. Alors, d'après 
\cite{Faltings7}, pour un choix de coordonnée formelle $V$ sur $\LT$
l'accouplement 
$$
G\times G^\vee \ldrt \LT[\pi]
$$
est donnée par
$$
V\longmapsto T\otimes U 
$$
Soit alors 
$$
\a_{G^\vee}^\O : G^\vee \ldrt \omega_G
$$
l'application de Hodge-Tate relative à $\O$.
Avec l'identification 
$$
\omega_G\simeq R/\delta R . dT
$$
cette application s'identifie donc à 
$$
u \longmapsto (u\text{ mod } \delta).dT
$$
Lorsque $R=\O_K$ avec $K$ comme précédemment on en déduit que l'on
connaît $\O_K.\text{Im} \a^\O_{G^\vee}$ dès que l'on connaît
$v(\text{Ann }\omega_G) = v(\delta)$ ou bien $v(\text{Ann }
\omega_{G^\vee} )= v(\delta)$.

\subsection{Calcul de la valuation $p$-adique de $\a_{H^\vee}^\O (x)$}\label{corian}

\subsubsection{Notations}

Soit $H$ un  $\O$-module $\pi$-divisible formel de dimension $1$ et de $\O$-hauteur
$n$ sur $\spec(O_K)$. Supposons également que sa fibre spéciale est formelle.
 Nous allons calculer $v(\a_{H^\vee}^\O (x))$ pour $x\in T_p
(H^D)$. Nous noterons $\widehat{H}$ le groupe formel associé sur $\spf
(\O_K)$. Alors $H(\O_K) = \bigcup_{k\geq 1} H[\pi^n] (\O_K) \subset
\widehat{H} (\O_K)$. Il y a une
``valuation''
$$
v: \widehat{H} (\O_K) \ldrt \R_{>0}
$$
qui définit une filtration dite de ramification inférieure sur
$\widehat{H} (\O_K)$ et donc sur les points de torsion (cf. section \ref{Filtraminf}).
 Cette ``valuation'' est définie de la façon suivante : fixons un isomorphisme de $\spf(\O_K)$-schémas formels pointés
$$
\widehat{H}\iso \spf (\O_K[[T]])
$$
où $\widehat{H}$ est pointé par sa section unité et $\spf (\O_K[[T]])$ par la section $T=0$. Cet isomorphisme induit une bijection $\widehat{H} (\O_K) \simeq \{x\in\O_K\;|\; v(x)>0 \}$. Si 
via cette bijection $y\in 
\widehat{H} (\O_K)$ correspond à $x\in \O_K$ on pose alors $v(y)=v(x)$. On vérifie aussitôt que cette définition ne dépend pas de l'isomorphisme de schémas formels pointés choisi. 

On utilisera systématiquement le
jeu entre la fibre générique et les modèles entiers en écrivant pour
$G$ un groupe fini localement libre sur $\O_K$ 
$$
G(\O_K) = G(K)
$$
et
$$
T_p (H) = \underset{k}{\limp} H[\pi^k] (K) = \underset{k}{\limp} H[\pi^k] (\O_K)
$$
\'Etant donné que $K$ est algébriquement clos on considérera toujours
les fibres génériques des groupes finis sur $\O_K$ comme des groupes
abstraits.

Soit $G$ un groupe $p$-divisible sur $\spec (\O_K)$ et $D$ un
sous-groupe fini de la fibre générique de $G$. On notera $D^{adh}$
l'adhérence schématique de $D$ dans $G[p^k]$ pour $k>>0$ (et cela ne
dépend pas de $k$). Dans la suite il n'y aura jamais d'ambiguïté pour
un $D$ donné sur le groupe $G$ dans lequel on prend l'adhérence
schématique, c'est pourquoi $G$ n'intervient pas dans la notation.

\subsubsection{Premiers calculs}

Pour $G$ un groupe fini localement libre muni d'une action stricte de
$\O$ 
le morphisme de faisceaux
fppf $\a_{G^\vee}^\O : G^\vee
\ldrt \omega_G$ est naturel en $G$, tout morphisme strict $f: G_1\ldrt G_2$
induit un diagramme commutatif 
$$
\xymatrix{
G_2^\vee \ar[r]^{ \a_{G_2^\vee}^\O} \ar[d]^{f^\vee} & \omega_{G_2} \ar[d] \\
G_1^\vee \ar[r]^{\a_{G_1^\vee}^\O} & \omega_{G_1}
}
$$
En particulier $\forall k\in \N^*$ l'inclusion $H[\pi^k]\hookrightarrow
H[\pi^{k+1}]$ induit un diagramme commutatif de morphismes de schémas en groupes
$$
\xymatrix@C=17mm{
H^\vee [\pi^{k+1}] \ar[r]^{\a_{H^\vee[p^{k+1}]}^\O} \ar@{->>}[d] & \omega_{H[\pi^{k+1}]}
\ar@{->>}[d] & \omega_H/\pi^{k+1} \omega_H  \ar[l]_\sim \ar@{->>}[d] \\
H^\vee [\pi^k] \ar[r]^{\a_{H^\vee[p^{k}]}^\O}  & \omega_{H[\pi^{k}]}
 & \omega_H/\pi^{k} \omega_H  \ar[l]_\sim
}
$$
et un diagramme de morphismes de groupes
$$
\xymatrix@C=12mm{
T_p (H^\vee) \ar[r]^{\a_{H^\vee}^\O} \ar@{->>}[d] & \omega_H \ar@{->>}[d] \\
T_p (H^\vee)/\pi^k T_p (H^\vee) \ar[r] \ar[d]^\simeq & \omega_H /\pi^k \omega_H
\ar[d]^\simeq \\
H[\pi^k]^\vee (\O_K) \ar[r]^{\a_{H[\pi^k]^\vee}^\O} & \omega_{H[\pi^k]}
}
$$
où $T_p (H^\vee)$ est le groupe des $(x_k)_{k\geq1}$, $x_k \in H[\pi^k]^\vee (K) =
H[\pi^k]^\vee (\O_K)$, $\pi x_{k+1} = x_k$. 
Ainsi si $x=(x_k)_{k\geq 1}$ pour calculer $\a_{H^\vee}^\O (x)$ il suffit de
calculer $\a_{H[\pi^k]^\vee}^\O (x_k)$ pour tout $k$ qui s'identifie à
$\a_{H^\vee}^\O(x)$ mod $\pi^k$.
\\

Soit donc $x\in T_p (H^\vee)$ dont on veut calculer $v(\a_{H^D}^\O (x))$. On
peut supposer que $x\notin \pi T_p (H^\vee)$ c'est à dire que le morphisme
associé $T_p (H) \ldrt \O_F (1)$ est surjectif où $F(1)$ désigne le
caractère de Lubin-Tate. 
 On fera donc cette
hypothèse. On constate que la valuation  de $\a_{H^D}^\O (x)$
ne dépend que du sous-module engendré $\O.x \subset T_p (H^\vee)$ qui
est facteur direct dans $T_p (H^\vee)$. Via la dualité parfaite
$$
T_p (H)\times T_p (H^\vee) \ldrt \O_F (1)
$$
de tels sous-modules correspondent aux sous-$\O$-modules $M \subset T_p
(H)$  facteur direct de rang $n-1$, $M=\left (\O.x \right )^\perp$.

Cela reste valable modulo $\pi^k$. Si $x\in H [\pi^k]^\vee (K)\setminus H
[\pi^{k-1}]^\vee (K)$, modulo une unité $\a_{H[\pi^k]^\vee} (x)$ ne dépend que du
sous-module engendré $C=<x>$ et de tels sous-modules sont en bijection
avec les sous-modules $C^\perp\subset H[\pi^k](K)$ facteurs directs de
rang $n-1$ sur $\O/\pi^k O$.

\begin{lemm}
L'opération d'adhérence schématique commute à la dualité de Cartier-Faltings :
si $C\subset H[\pi^k]^\vee (K)$ est un sous-groupe alors 
$$
\left (C^{adh} \right )^\vee \simeq H[\pi^k] / \left (C^\perp \right )^{adh}
$$
\end{lemm}
\dem
De la suite exacte 
$$
0\ldrt C^{adh} \ldrt H[\pi^k]^\vee \ldrt H[\pi^k]^\vee/C^{adh}\ldrt 0
$$
on déduit d'après le théorème 8 de \cite{Faltings7} la suite exacte
$$
0 \ldrt \left (  H[\pi^k]^\vee/C^{adh}\right )^\vee \ldrt H[\pi^k] \ldrt \left (
  C^{adh}\right )^\vee \ldrt 0
$$
Le sous-groupe fini localement libre de gauche coïncide en fibre
générique avec $C^\perp$. Il est donc égal à $\left ( C^\perp \right
)^{adh}$.
\qed

\begin{prop}\label{rrgkjm}
Soit $D \subset H$ un sous-groupe fini localement libre sur
$\O_K$. Alors $\omega_D \simeq \O_K/\gamma \O_K$ où 
$$
v (\gamma) = \sum_{\l\in D\setminus 0} v(\l)
$$
Le complexe de co-Lie de $D$ est isomorphe au complexe
$$
\left [ \O_K \xrig{\; \gamma \; } \O_K \right ]
$$
\end{prop} 
\dem
Avec un choix de bonnes ``coordonnées formelles'' (on entend par là un
isomorphisme de schémas formels pointés entre $\widehat{H}$ et  $\spf (\O_K [[T]])$) à la
source et au but l'isogénie de groupes formels 
$$
\widehat{H} \ldrt \widehat{H}/C
$$
s'écrit $\dpt{T \prod_{x\in D\setminus \{0 \}} (T-x)}$.
\qed

\begin{rema}
Dans cette dernière proposition l'assertion concernant la valuation de
$\gamma$ 
 est l'analogue de la proposition 4 du
chapitre IV de \cite{Serre1} reliant valuation de la différente et les
groupes de ramifications inférieurs d'une extension de corps locaux.
\end{rema}

\begin{coro}\label{drhg}
Soient $D_1 \subset D_2$ des groupes fini localement libres sur $\O_K$
sous-groupes  de $H$. Alors la suite 
$$
0 \ldrt \omega_{D_2/D_1} \ldrt \omega_{D_2} \ldrt \omega_{D_1} \ldrt 0
$$
est exacte.
\end{coro}
\dem
D'après la proposition précédente le groupe de cohomologie $H^{-1}$ du
complexe de co-Lie de nos groupes est nul puisque $\O_K$ est sans $p$-torsion.
\qed
\\

Soit donc maintenant $C=\O.y\subset H[\pi^k]^\vee (K)$ facteur direct de rang 
$1$ et notons $C^\perp \subset H[\pi^k](K)$ son orthogonal. 
Considérons le diagramme 
$$
\xymatrix{
C^{adh}  \ar[r]^{\a_1=\a_{C^{adh}}^\O} \ar@{->>}[d]^{q_1} & \omega_{(C^{adh})^\vee}
\ar@{^(->}[r] \ar@{->>}[d]^{q_2} &
\omega_H/\pi^k \omega_H \\
C^{adh}/C[\pi^{k-1}]^{adh} \ar[r]^{\a_2} & \omega_{(C^{adh}/C[\pi^{k-1}]^{adh})^\vee}
}
$$
L'isomorphisme $\left ( C^{adh}\right )^\vee \simeq 
  H[\pi^k]/\left ( C^\perp \right )^{adh}$ implique que si $\omega_{(C^{adh})^\vee}
\simeq \O_K/\gamma \O_K$ alors 
$$
v(\gamma ) = k- \sum_{z\in C^\perp \setminus \{ 0 \} } v(z)
$$
De même si $\omega_{(C[\pi^{k-1}]^{adh})^\vee} \simeq \O_K/\gamma' \O_K$ alors 
$$
v(\gamma') = k-1-\sum_{z\in C^\perp [\pi^{k-1}] \setminus \{ 0 \}} v(z)
$$
Soit maintenant $\gamma''$ tel que $\omega_{(C^{adh}/C[\pi^{k-1}]^{adh})^\vee}
\simeq \O_K/\gamma'' \O_K$. On  déduit donc du corollaire \ref{drhg} que
$$
v(\gamma'')=1-\sum_{z\in C^\perp \setminus C^\perp [\pi^{k-1}]} v(z)
$$
Nous allons maintenant utiliser les résultats de la section
\ref{corian}. Le groupe $C^{adh}/C[\pi^{k-1}]^{adh}$ est de type $(p,\dots,p)$ 
et sont dual strict vérifie les hypothèses de la section \ref{corian}.
 Avec les notations du diagramme précédent, $q_1 (x)$ engendre
les points à valeurs dans $K$ de ce groupe comme $\O/\pi \O$-module. On en
déduit que $\a_2 (q_1 (x))= \beta \text{ mod } \gamma'' \O_K$ où 
$$
v(\beta ) = \frac{\sum_{z\in C^\perp \setminus C^\perp [\pi^{k-1}]} v(z)}{q-1}
$$
et donc, si 
$$
v(\beta) < v(\gamma")
$$
c'est à dire si  
$$
\sum_{z\in C^\perp \setminus C^\perp [\pi^{k-1}]} v(z)<1-\frac{1}{q}
$$
alors $\a_1 (y) \neq 0 \in \O_K/ \gamma \O_K$ et
$$
v(\a_1 (y)) = \frac{\sum_{z\in C^\perp \setminus C^\perp [\pi^{k-1}]} v(z)}{q-1}
$$
Si maintenant notre élément $y\in H[\pi^k]^\vee (K)$ provient d'un $x\in T_p (H^\vee)
\setminus \pi T_p (H^\vee)$, puisque $\omega_{(C^{adh})^D} \hookrightarrow
\omega_H/\pi^k \omega_H$ est une injection on a  
\begin{eqnarray*}
v(\a_{H^\vee}^\O (x)) &=& \frac{\sum_{z\in C^\perp \setminus C^\perp [\pi^{k-1}]}
  v(z)}{q-1} + k - v(\gamma) \\
&=& \frac{\sum_{z\in C^\perp \setminus C^\perp [\pi^{k-1}]} v(z)}{q-1}
+ \sum_{z \in C^\perp \setminus \{ 0 \}} v(z)
\end{eqnarray*}

\begin{defi}
Soit $A$ un ensemble de points de torsion de $H$. On note 
$$
v(A) = \sum_{z\in A\setminus \{ 0\}} v(z)
$$
\end{defi}

Résumons ce que nous avons démontré jusqu'à maintenant.

\begin{prop}\label{balkla}
Soit $x\in T_p (H^\vee)\setminus \pi T_p (H^\vee)$ et $M=
(\O.x)^\perp\subset T_p(H)$ facteur direct de rang $n-1$. 
 Notons pour tout entier $k\geq 1$ $\;M[\pi^k]$ le
sous-groupe des points de $\pi^k$-torsion associé dans $H[\pi^k](K)$. Si l'entier
$k$ est tel que
$$
v(M[\pi^k]\setminus M[\pi^{k-1}])<1-\frac{1}{q}
$$
alors
\begin{eqnarray*}
v(\a_{H^\vee}^\O)(x) &= &\frac{v(M[\pi^k]\setminus M[\pi^{k-1}])}{q-1} +
v(M[\pi^k]) \\
&=& \frac{1}{q-1} \left ( q v(M[\pi^k]) - v(M[\pi^{k-1}])\right )
\end{eqnarray*}
\end{prop}

Reste à voir qu'il existe un tel entier $k$, ce que nous allons
faire.

\subsection{La formule finale} \label{fbmqsd}

\begin{figure}[htbp]\label{kolde}
   \begin{center}
      \includegraphics{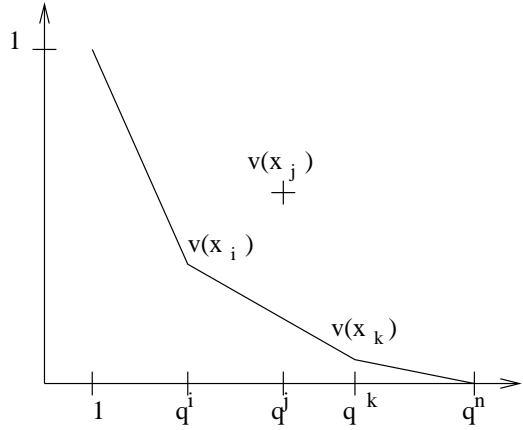}
   \end{center}
   \caption{\footnotesize Le polygone de Newton de
    la multiplication par $\pi$}
   \label{NewGH}
\end{figure}

Rappelons maintenant (\cite{HopkinsGross}, \cite{Cellulaire} chapitre 1) qu'il existe une loi de groupe formelle associée à
$H$ telle
que le polygone de Newton de la multiplication par
$\pi$ sur cette loi soit l'enveloppe convexe des points 
$$
(0,\infty), (1,1), (v(x_1),q),\dots, (v(x_i),q^i),\dots,
(v(x_{n-1}),q^{n-1}), (0,q^n) 
$$
où $x_1,\dots,x_{n-1} \in K$ et $\forall i\; v(x_i)>0$ (cf. figure
\ref{kolde}).
Rappelons également la recette suivante : si $y\in \widehat{H} (\O_K)$
la valuation des points s'envoyant sur $y$ par la multiplication par
$\pi$ sur $\widehat{H} (\O_K)$ s'obtient en prenant l'enveloppe convexe
du polygone précédent et du point $(0, v(y))$. 

\begin{lemm}\label{vaenzero}
On a 
$$
\underset{k\drt +\infty}{\lim}
\underset{z\in H[\pi^k]\setminus H[\pi^{k-1}] }{\text{sup}} v(z) =0
$$
\end{lemm}
\dem Remarquant que 
$$
H[\pi^{k+1}]\setminus H[\pi^k] = \{\; y \in \widehat{H} (\O_K)\;|\; \pi y \in
H[\pi^k]\setminus H[\pi^{k-1}]\;\} 
$$
on en déduit facilement le résultat en utilisant les rappels
précédents sur le polygone de Newton. 
\qed

Ainsi les valuations des points de torsion de $\widehat{H} (\O_K)$
forment un ensemble discret dans $]0,+\infty [$ s'accumulant en $0$.

\begin{prop}
Soit $M \subset T_p (H)$ un sous-$\O$-module facteur direct de rang $n-1$. Alors, 
$$
\text{pour } k>>0 \;\; \forall y\in M[\pi^k]\setminus M[\pi^{k-1}] \;\; \forall
z\in H[\pi] (K) \;\;\; \frac{v(y)}{q^n} \leq v(z)
$$
De plus si $k$ est un entier vérifiant cette condition alors
$$
v(M[\pi^{k+1}]\setminus M[\pi^k])= \frac{1}{q} v ( M[\pi^k]\setminus M[\pi^{k-1}])
$$
\end{prop}
\dem
La première assertion résulte du lemme précédent. Quant à la seconde,
il y a une application surjective 
$$
\text{multiplication par } \pi : 
M[\pi^{k+1}]\setminus M[\pi^k] \twoheadrightarrow 
M[\pi^k]\setminus M[\pi^{k-1}] 
$$
et si $k$ vérifie la condition énoncée alors d'après les rappels sur
le polygone de Newton de la multiplication par $\pi$ 
$$
\forall z \in H[\pi^{k+1}]\setminus H[\pi^k] \;\; v(z) = \frac{1}{q^n}
v(\pi .z)
$$
Remarquant que par l'application multiplication par $\pi$ précédente un 
élément a $q^{n-1}$ antécédents on en déduit le résultat.
\qed

\begin{coro}
La condition de la proposition \ref{balkla} est toujours vérifiée pour
$k>>0$. La fonction 
$$
k\longmapsto \frac{1}{q-1} \left ( q v(M[\pi^k])- v(M[\pi^{k-1}])\right )
$$
est constante  dès que $k\geq k_0$ où $k_0$ est tel que $\forall y \in M[\pi^{k_0}]\setminus M[p^{k_0-1}]
\; \forall z\in H[\pi] \; v(y)\leq q^n v(z)$. Elle 
fournit donc une formule pour la valuation de $\a_{H^\vee}^\O (x)$ lorsque
$k$ vérifie cette condition.
\end{coro}

\begin{lemm}
Soit $M \subset T_p (H)$ facteur direct de rang $n-1$.
Soit $M\otimes \Qp/\Zp \subset H (\O_K)$ l'ensemble des points de torsion associé. 
 Alors la
fonction 
$$
k\longmapsto \frac{1}{q-1} \left ( q v(M[\pi^k])- v(M[\pi^{k-1}])\right )
$$
est égale à $v \left ( M\otimes \Qp/\Zp\right)$ pour $k>>0$.
\end{lemm}
\dem
Soit $k_0$ un entier comme dans le corollaire précédent.
\begin{eqnarray*}
v \left ( M\otimes \Qp/\Zp \right)  &=& v(M[\pi^{k_0}]) +
\sum_{k\geq k_0+1}  v\left ( M[\pi^{k}] \setminus M[\pi^{k-1}]\right ) \\
&=&  v(M[\pi^{k_0}]) +
\sum_{k\geq k_0+1} \left ( v(M[\pi^{k}]) - v(M[\pi^{k-1}]) \right ) \\
&=&  v(M[\pi^{k_0}]) + (v(M[\pi^{k_0+1}])-v (M[\pi^{k_0}])) \sum_{i\geq 0}
q^{-i} \\
&=&  \frac{1}{q-1} \left ( p v(M[\pi^{k_0+1}])- v(M[\pi^{k_0}])\right )
\end{eqnarray*}
\qed

Résumons tous les résultats obtenus dans le théorème suivant. 

\begin{theo}\label{tklopui}
Soit $x\in T_p (H^\vee)$ non-nul et notons $\ph : H\ldrt \LT$ le
morphisme associé.
Alors 
$$
v(\a_{H^\vee}^\O(x)) = \sum_{z\in \ker \ph_{| H[p^\infty](K)\setminus \{ 0 \}}} v(z) 
$$
Il existe de plus un entier $k_0$ tel que 
$$
\forall y\in \ker \ph [\pi^{k_0}](K) \setminus \ker \ph [\pi^{k_0-1}](K)\;\; \forall y\in
\ker \ph[\pi](K) \;\; v(y)\geq q^n v(z)
$$
Alors, si $k\geq k_0$ 
$$
v(\a_{H^\vee}^\O (x) )= \frac{1}{q-1} \left (q v(\ker \ph [\pi^k])-v
(\ker \ph [\pi^{k-1}])\right )
$$
\end{theo}
\dem
Il existe un entier $k\geq 0$ tel que $x\in \pi^k T_p
(H^\vee)\setminus \pi^{k+1} T_p (H^\vee)$. Quitte à considérer la
factorisation
$$
\ph : H \twoheadrightarrow H/H[\pi^k] \ldrt \LT
$$
on peut se ramener au cas où $x\notin \pi T_p (H^\vee)$ étudié
précédemment. En effet, si $C$ est un sous-groupe fini de $H$ contenant
$H[\pi^k]$ alors d'après la proposition \ref{rrgkjm} et le corollaire
\ref{drhg}  
$$
\sum_{x\in C\setminus \{0 \}} v(x) = k+ \sum_{x\in (C/H[\pi^k]) \setminus
  \{0 \}} v(x)
$$
où $C/H[\pi^k] \subset H/H[\pi^k]$ et la fonction valuation sur
$C/H[\pi^k]$ est celle déduite de la fonction valuation sur le
$\O$-module formel $H/H[\pi^k]$.
\qed

\begin{rema}
Soit $f:\widehat{H}_1\ldrt \widehat{H}_2$ une isogénie de groupes
formels de dimension $1$ sur $\O_K$. Après choix d'isomorphismes de schémas formels pointés $\widehat{H}_1\simeq \spf
(\O_K [[T]])$, resp. $\widehat{H}_2\simeq \spf (\O_K [[T]])$, d'après le théorème de
factorisation de Weierstrass 
$$
f^* T = \prod_{y\in \ker f} (T-y) \times u
$$
où $u\in \O_K[[T]]^\times$. Alors, l'application induite $f^*:
\omega_{H_2}\ldrt \omega_{H_1}$ s'identifie à 
$$
\O_K \xrig{\;\times f'(0) = u(0)\prod_{y\in \ker f \setminus \{ 0 \}} y \;}  \O_K
$$
et $v(f'(0))= \sum_{y\in \ker f \setminus \{ 0 \}} v(y)$. 
Pour un $\ph : H\ldrt \LT$ comme dans le théorème précédent, si $n>1$,
après un choix d'isomorphismes $\widehat{H}\simeq \spf (\O_K[[T]])$,
$\widehat{\LT}\simeq \spf (\O_K [[T]])$, $\ph^* T\in \O_K [[T]]$ ne
vérifie pas les conditions du théorème de factorisation de
Weierstrass. Mais la formule donnée dans le théorème précédent dit que
pour $v \left ( (\ph^*T)'(0)\right )$ tout se passe comme si 
$$
\ph^* T = \prod_{y\in \ker \ph\setminus \{ 0 \} } (T-y) \times \text{unité}
$$
bien que ce produit infini n'ait pas de sens. Plus précisément, pour
tout entier $k\geq 1$ on peut factoriser 
$$
\ph^* T = \prod_{y\in \ker \ph [\pi^k]\setminus \{ 0 \}} (T-y) \times
g_k (T)
$$
où $g_k (0) \notin \O_K^\times$ ce qui exprime le fait que le noyau du
morphisme $\ph_{|H[\pi^k]}: H[\pi^k]\ldrt \LT [\pi^k]$ n'est pas
plat.
 Néanmoins le théorème précédent est équivalent à ce que 
$$
\underset{k\ldrt +\infty}{\lim} v(g_k (0)) =0
$$
L'auteur de cet article a essayé de trouver une démonstration du théorème
précédent sous la forme que l'on vient de donner
à partir de manipulations élémentaires sur les séries
formelles mais n'a pas réussi. Bien sûr les considérations précédentes
montrent tout de même que $v(\a^\O_{H^\vee}  (x))\geq \sum_{y\in\ker
  \ph\setminus \{0 \}} v(y)$.   
\end{rema}

\section{Filtration de ramification inférieure}\label{Filtraminf}

Soit $K|F$ un corps valué complet pour une valuation $v$ à valeurs
dans $\R$ étendant celle de $F$. On supposera comme dans le chapitre
précédent que $K=\widehat{\overline{K}}$. 

Soit $H$ un $\O$-module $\pi$-divisible de dimension $1$ et de
$\O$-hauteur $n$ sur $\spec (\O_K)$ dont la fibre spéciale est
formelle. 

\subsection{Définition et premières propriétés}\label{knuypot}

Posons $V=V_p(H)$, un $F$-espace vectoriel de dimension $n$. On a
l'identification
$$
V_p (H) = \{\; (x_i)_{i\in\Z}\;|\; x_i \in \widehat{H} (\O_K)\; \pi
x_{i+1} = x_i \text{ et } x_i=0 \text{ pour } i<<0\;\}
$$
Avec cette identification soit 
\begin{eqnarray*}
v : V &\ldrt & ]0,+\infty ] \\
(x_i)_{i\in\Z} & \longmapsto & v(x_0)
\end{eqnarray*}
Via l'identification $V_p (H)\otimes \Qp/\Zp = H(\O_K)$ cette
``valuation'' $v$ s'écrit aussi comme l'application composée 
$$
 V_p (H) \ldrt V_p (H)\otimes\Qp/\Zp = H(\O_K) \subset \widehat{H}
(\O_K) \xrig{\; v \; } ]0,+\infty ] 
$$
Elle vérifie 
\begin{itemize}
\item $\forall x,y\in V\;\; v(x+y) \geq \inf \{ v(x),v(y)\}$
\item $\forall a\in \O\; \forall x\in V\;\; v(ax)\geq v(x) $ et
  $v(\pi x)> v(x)$
\item $\forall x\in V\;\; v(\pi^k x) =+\infty$ pour $k>>0$ \\
et $v(\pi^{k-1} x) = \frac{1}{q^n} v(\pi^k x)$ pour $k<<0$
\end{itemize} 

\begin{defi}
Pour tout $\l \in ]0,+\infty]$ posons 
$$
\Fil_\l V= \{\; x\in V\; |\; v(x)\geq \l \; \}
$$
qui définit donc une filtration décroissante de $V$.
\end{defi}

 Il résulte des
propriétés énoncées de $v$ que les $\Fil_\l V$ sont des
sous-$\O$-modules. De plus
\begin{itemize}
\item D'après le lemme \ref{vaenzero} les $\Fil_\l V$ sont des
  réseaux dans $V$ 
\item $\Fil_\infty V= T_p (H)$
\item Pour $\l \geq \frac{1}{q-1} \;\; \Fil_\l V= \Fil_\infty V$
\item $\exists \e >0$ tel que $\forall 0<\l \leq \e \;\; \pi^{-1}
  \Fil_\l = \Fil_{\frac{\l}{q^n}}$
\end{itemize}

La filtration est donc entièrement déterminée par un nombre fini de
réseaux, les autres étant déterminés par
périodicité. Ses sauts forment un ensemble discret 
s'accumulant en $0$ dans $\R_{>0}$. Plus précisément, elle a la forme
suivante
$$
\Fil_{\l_1}\subsetneq \dots \subsetneq \Fil_{\l_r} \subsetneq
\Fil_{\l_{r+1}}
\subsetneq \dots \subsetneq \Fil_{\l_{r+t}} \subsetneq \pi^{-1} 
\Fil_{\l_{r+1}}
\subsetneq \dots \subsetneq \pi^{-1} \Fil_{\l_{r+t}} \subsetneq 
\pi^{-2} \Fil_{\l_{r+1}} \subsetneq \dots
$$
où $\l_1>\dots> \l_{r} >\l_{r+1}>\dots >\l_{r+t}$. 

La ``valuation'' $v$ sur $V$ donne lieu à une ``valuation''
$$
v: V/\Fil_\infty \ldrt ]0,+\infty]
$$
Pour tout entier $k\geq 1$ il y a un isomorphisme naturel
$$
\pi^{-k} \Fil_\infty /\Fil_\infty \iso H[\pi^k] (\O_K)
$$
et via cet isomorphisme la ``valuation'' n'est rien d'autre que la
valuation sur les points de $\pi^k$-torsion. Ainsi 
$$
\forall k\geq 1\; \; \left (\Fil_\l \cap \pi^{-k} \Fil^\infty \right )/\Fil_\infty
\iso \{ \; x\in H[\pi^k] (\O_K)\;|\; v(x)\geq \l \;\}
$$
$$
\Fil_\l/\Fil_\infty = \{\; x\in H(\O_K)\;|\; v(x)\geq \l \;\}
$$

Notons enfin deux propriétés cruciales qui n'ont pas encore été énoncées de la
``valuation'' $v$ sur $V/\Fil_\infty V$ :
$$
(*)\;\; \forall x\in V/\Fil_\infty V \;\; \{ v(y)\; |\; y\in
V/\Fil_\infty V\; \pi y =x \;\} \text{ ne dépend que de } v(x)
$$
$$
(**) \;\; \forall x,x'\in V/\Fil_\infty V\;\; v(x)\leq v(x') \limpl 
\sup \{ v(y)\;|\; \pi y =x\; \} \leq \sup \{ v(y')\;|\; \pi y' =x'\;\}
$$
Celles-ci résultent des rappels faits au début de la section \ref{fbmqsd}
sur le polygone de Newton de la
multiplication par $\pi$ sur une loi de groupe formel associée à $H$.

\subsection{La filtration de ramification est contenue dans un
  appartement de l'immeuble}

Soit $\mathcal{I} (V)$ l'immeuble de $\text{PGL} (V)$ vu comme
ensemble simplicial
(cf. appendice \ref{kujitopo}). Considérons les classes d'homothéties
$$
\left ( [\Fil_\l V]\right )_{\l\in ]0,+\infty ]}
$$
Elles forment un ensemble fini de sommets dans $\mathcal{I} (V)$. 
On renvoie à l'appendice \ref{kujitopo} pour les définitions de bases concernant les
appartements et les quartiers dans l'immeuble $\mathcal{I} (V)$. 

\begin{theo}\label{Ramapp}
L'ensemble fini $\{\; [\Fil_{\,\l} V] \; |\; \l\in ]0,+\infty ]\;\}$ est
contenu dans un appartement de $\mathcal{I} (V)$. Plus précisément,
cet ensemble est contenu dans un quartier de sommet $[\Fil_\infty V ]$
dans un appartement. Il est également connexe.
\end{theo}
\dem
Considérons la filtration par des sous $\Fq$-espaces vectoriels 
 sur les points de $\pi$-torsion de $H$
$$
\Fil_\l \left ( \pi^{-1}\Fil_\infty V/\Fil_\infty V \right )\;\;\;\l\in
]0,+\infty ]
$$
Soit $ (e_i^{(1)} )_{1\leq i \leq n}$ une base de $\pi^{-1} \Fil_\infty
V /\Fil_\infty V$ scindant la filtration précédente, c'est à dire
telle que pour des entiers $1\leq \a_1 <\dots <\a_r = n$ cette
filtration soit donnée par 
$$
(0)\subsetneq <e_1^{(1)},\dots,e_{\a_1}^{(1)}> \subsetneq \dots
\subsetneq <e_1^{(1)},\dots, e^{(1)}_{\a_i}> \subsetneq \dots
\subsetneq  <e_1^{(1)},\dots, e^{(1)}_{n}>
$$
Pour tout entier $i$ tel que $1\leq i\leq n$ soit la suite
$(e_i^{(k)})_{k\geq 1}$ définie par récurrence de la façon suivante :
\begin{eqnarray*}
e_i^{(k)}& \in & \pi^{-k}\Fil_\infty V /\Fil_\infty V \\
\pi e_i^{(k+1)} &=& e_i^{(k)} \\
\text{ et } v(e_i^{(k+1)}) &=& \sup \{\; v(x) \; | \; \pi x
=e_i^{(k)}\; \}
\end{eqnarray*}
Les $(e_i^{(k)})_{1\leq i \leq n}$ forment donc une base du $\O/\pi^k
\O$-module libre $\pi^{-k}\Fil_\infty V/\Fil_\infty V$. La relation de
compatibilité $\pi e_i^{(k+1)} = e_i^{(k)}$ implique qu'ils
fournissent une base du module de Tate $T_p (H)$. En termes d'algèbre
linéaire cela se traduit de la façon suivante : si pour tout $i$ et
$k$ on fixe $\widetilde{e_i^{(k)}} \in V$ un relèvement de $e_i^{(k)}
\in V/\Fil_\infty V$ alors en posant 
$$
\forall i \;\; e_i = \underset{k\drt + \infty}{\lim} \pi^k
\widetilde{e_i^{(k)}} \in \Fil_\infty V
$$
On a 
$$
\Fil_\infty V = <e_1,\dots, e_n> \text{ et } e_i^{(k)} \equiv \pi^{-k}
e_i \text{ mod } \Fil_\infty 
$$
Montrons que $\forall \l\;\; \Fil_\l V$ est dans l'appartement associé
à la base $(e_1,\dots, e_n)$ c'est à dire s'écrit sous la forme 
$<\pi^{-k_1} e_1,\dots, \pi^{-k_n} e_n>$ où $k_i\in \N$.
D'après le lemme \ref{kopghut} de l'appendice \ref{kujitopo} il suffit de montrer que 
$\forall (a_i)_i \in F^n \;\; v (\sum_i a_i e_i )$ ne dépend que des
$(v(a_i))_{1\leq i \leq n}$. 
\\

Posons pour tout $k\leq 0$ $e_i^{(k)}=0\in V/\Fil_\infty V$. Tout $x
\in V/\Fil_\infty V$, $x\neq 0$, peut s'écrire sous la forme 
$$
x= \sum_{j=1}^l \sum_{i=\delta_{j}+1}^{\delta_{j+1}} a_i e_i^{(k(j))}
$$
où 
$$
0=\delta_1 <\delta_2 < \dots <\delta_{l+1} =n
$$
$$
k(1)>k(2)>\dots >k (l ) \geq 0
$$
$$
\forall i\; a_i \in \O_F \text{ et }
\forall j>1 \;\; a_{\delta_j}\in \O_F^\times
$$
Nous allons démontrer par récurrence sur l'entier $k\geq 1$
l'assertion suivante :
\begin{eqnarray*}
(A_k) \;\; &&\text{Soit } x \text{ écrit sous la forme précédente avec }
k(1)=k \\
&&\text{alors  }  v(x) = \inf\{ v( e_{\delta_{j+1}}^{(k(j))})\;| 1 \leq
j \leq l\;\} 
\end{eqnarray*}
Si $k=1$, $(A_1)$ est vérifiée par définition des
$(e_i^{(1)})_i$. Supposons donc $(A_k)$ vérifiée et soit $x$ comme
précédemment avec $k(1)=k+1$. Alors, d'après la propriété $(**)$ de la
``valuation'' $v$ sur $V/\Fil_\infty V$ donnée à la fin de la section \ref{knuypot}
$$
\forall j\in \{ 1,\dots,l\}\;\; v(e_{\delta_{j}+1}^{(k(j))})\geq 
v(e_{\delta_{j}+2}^{(k(j))}) \geq \dots \geq v( e_{\delta_{j+1}}^{(k(j))})
$$
et donc 
$$
v(x)\geq \inf \{ v( e_{\delta_{j+1}}^{(k(j))})\; |\; 1 \leq j\leq l\;\}
$$
De plus, par hypothèse de récurrence 
$$
v(\pi x) = \inf \{ v ( e_{\delta_{j}+1}^{(k(j)-1)})\;|\; 1 \leq j \leq l
\; \}
$$
Soit $j_0 \in \{ 1, \dots,l \}$ tel que 
$$
v( e_{\delta_{j_0+1}}^{(k(j_0)-1)}) =  \inf \{ v ( e_{\delta_{j}+1}^{(k(j)-1)})\;|\; 1 \leq j \leq l
\; \}
$$
Alors, toujours d'après la propriété $(**)$ 
$$
\inf \{ v ( e_{\delta_{j+1}}^{(k(j))})\;|\; 1 \leq j \leq l \; \} =v(e_{\delta_{j_0+1}}^{(k(j_0))})
$$
Or d'après la propriété $(*)$
$$
\sup \{ v(y)\;|\; \pi y =\pi x\;\} = \sup \{ v(z) \;|\; \pi z =
e_{\delta_{j_0+1}}^{(k(j_0)-1)}\; \} = v( e_{\delta_{j_0+1}}^{(k(j_0))})
$$
Donc 
$$
v(x) = v(e_{\delta_{j_0+1}}^{(k(j_0))})
$$
d'où $(A_{k+1})$. 

Le résultat sur le fait que la filtration de ramification est contenue
dans l'appartement s'en déduit. 
\\

Vérifions maintenant que cette filtration est contenus dans le
quartier de sommet $[\Fil_\infty]$ égal à 
$$
\{[<\pi^{a_1} e_1,\dots, \pi^{a_n} e_n>] \; |\; a_1 \leq a_2\leq \dots
\leq a_n \;\}
$$
Pour cela il suffit de vérifier que 
$$
\forall i\leq j \; \forall a \in\Z\; v(\pi^a e_i) \geq v(\pi^a e_j)
$$
ce qui est clair une nouvelle fois grâce à la propriété (**).
\\

La connexité de l'ensemble de sommets associé résulte de ce qu'étant
donné que $\forall x\in V\setminus \{ 0\}$ $\; v(\pi x)>v(x)$ on a 
$$\forall \l \neq +\infty\; \exists\l'>\l \; \pi \Fil_\l \subset
\Fil_{\l'} \subsetneq \Fil_\l
$$ 
\qed

\begin{prop}\label{Muropk}
Soit $(e_1,\dots,e_n)$ la base construite dans le théorème
précédent. Soient $\l_1\geq \dots \geq \l_n$ les valuations des points
de $\pi$-torsion non-nuls de $H$, où il y a $q^i-q^{i-1}$ points de
valuation $\l_i$. Soit $I=\{ i\in \{ 1,\dots ,n-1 \}\;|\;
\l_i=\l_{i+1} \}$. A un $i\in I$ est associé un mur de
l'appartement associé à la base $(e_1,\dots,e_n)$ défini par 
$$
\{[<\pi^{a_1} e_1,\dots, \pi^{a_n} e_n>]\;|\; a_i = a_{i+1} \}
$$
L'ensemble fini
$[\Fil_\l V ]_{\l}$ est contenu dans l'intersection de ces murs
associés a $I$.
\end{prop}
\dem 
Cela se déduit aisément de la démonstration du théorème précédent.
\qed

\subsection{L'algorithme de calcul de la filtration de ramification}\label{ghuji}

De la démonstration du théorème \ref{Ramapp} on peut extraire
l'algorithme suivant. Soient 
$$
\l_1^{(1)}\geq \dots \geq \l_n^{(1)}
$$
les pentes du polygone de Newton noté $\mathcal{N}$ de la
multiplication par $\pi$ sur une loi de groupe formelle associée à $H$
où $\l_i^{(1)}$ désigne la pente comprise entre les abscisses
$q^{i-1}$ et $q^i$. Définissons pour tout $i$ la suite
  $(\l_i^{(k)})_{k \geq 1}$ par récurrence de la façon suivante :
$\l_i^{(k+1)}$ est la plus grande pente de l'enveloppe convexe de
$\mathcal{N}$ et du point $(0,\l_i^{(k)})$, c'est à dire la pente
comprise entre les abscisses $0$ et $1$.

On a donc $\l_i^{(k+1)}<\l_i^{(k)}$ et pour $k>>0\;\;
\dpt{\l_i^{(k+1)}= \frac{1}{q^n} \l_i^{(k)}}$. 

D'après la démonstration du théorème \ref{Ramapp} il existe une base
$(e_1,\dots,e_n)$ de $T_p (H) = \Fil_\infty V$ telle que 
$$
\forall i\;\forall k\geq 1 \;\; v(\pi^{-k}e_i)=\l_i^{(k)}
$$
Pour tout $\mu\in\R_{>0}$ et tout $i$ soit 
$$
k(i,\mu) = \sup \{k\geq 0\;|\; \l_i^{(k)}\geq \mu \; \}
$$
où l'on a posé $\forall i\; \l_i^{(0)}=+\infty$.
Alors 
$$
k(1,\mu)\geq k(2,\mu) \geq \dots \geq k(n,\mu)
$$
et 
$$
\boxed{
\Fil_\mu V = \bigoplus_{i=1}^n \O_F.  \pi^{-k(i,\mu)} e_i
}
$$
De plus, pour une telle base 
$$
\forall (x_1,\dots,x_n)\in F^n\;\; \boxed{ v(\sum_{i=1}^n x_i e_i) = \inf
\{\l_i^{(-v(x_i))}\;|\; 1\leq i\leq n \} }
$$
où l'on a posé $\l_i^{(k)}=+\infty$ si $k\leq 0$.

\begin{defi}\label{kkddssee}
On appellera base adaptée du module de Tate de $H$ une base
$(e_1,\dots,e_n)$ de
$\Fil_\infty V$ telle que 
$$
\forall 1\leq i\leq n \;\forall k\geq 1 \;\; v(\pi^{-k} e_i) =\l_i^{(k)}
$$
\end{defi}

Ainsi pour une base adaptée $(e_1,\dots,e_n)$ la filtration de
ramification est située dans le quartier 
$$
\{[<\pi^{-a_1} e_1,\dots,\pi^{-a_n} e_n>] \;|\; a_1\geq \dots \geq a_n \}
$$

\subsection{L'application $\l \longmapsto \l^{(k)}$ et la fonction de Herbrand}\label{lapelctose}

On renvoie à l'appendice \ref{Ramitofm} pour les définitions et propriétés de base concernant les fonctions de Herbrand.

Soit $\mathcal{N} : [0 , q^n] \ldrt \R$ le polygone de Newton de
pentes $\l_1\geq \dots \geq \l_n$. 
Pour $\l\in ]0,+\infty [$ on note $(\l^{(k)})_{k\geq 1}$ la suite
définie par $\l^{(1)} = \l$ et $\l^{(k+1)}$ est la plus grande pente
de l'enveloppe convexe de $\mathcal{N}$ et  de $(0,\l^{(k)})$. 

Pour une fonction convexe $\ph$
considérons sa duale (cf. appendice \ref{Ramitofm})
$$
\ph^* (t) = \sup \{ u\;|\; \ph (\bullet) \geq -t \bullet + u \}
$$
Notons $\eta (\bullet) = \mathcal{N}^{\,*} (\bullet)$ et $\Psi = \eta^{-1}$ la fonction de
Herbrand de $H[\pi]$ (cf. appendice \ref{Ramitofm}). Alors 
$$
\forall \l\;\; \boxed{\l^{(k)} = \Psi^{\circ (k-1)} (\l)}
$$
et $ \Psi^{\circ (k-1)}$ est la fonction de Herbrand de
  $H[\pi^{k-1}]$. On a $ \Psi^{\circ (k-1)} =( \eta^{\circ
    (k-1)})^{-1}$ et 
$$
 \eta^{\circ
    (k-1)} (s) = \int_0^s |\{ x\in H[\pi^{k-1}]\;|\; v(x)\geq t \} |
  \, dt
$$

\subsection{La facette de l'immeuble associée à la filtration de
  ramification}\label{kkjjuuii}

\subsubsection{Première définition}

\begin{defi}\label{quatoki}
Soit $\e\in ]0,+\infty[$ tel que $\e\leq \l_n$, c'est à dire $\forall
x\in \pi^{-1} \Fil_\infty V\; v(x)\geq \e$. On note $S$ le simplexe de
sommets 
$$
\{[\Fil_\l V]\;|\; 0<\l \leq \e\;\}
$$
dans l'immeuble $\mathcal{I} (V)$. 
\end{defi}

Cette définition a bien un sens car pour $0<\l\leq \e$ on a
$\Fil_{\l/q^n} V = \pi^{-1} \Fil_\l V$ qui définit donc une suite
périodique de réseaux. 

\subsubsection{Deuxième définition comme intersection de demi-appartements}

Soit $(e_1,\dots, e_n)$ une base adaptée du module de Tate de $H$ (définition \ref{kkddssee}). Soit
$\mathcal{A}$ l'appartement de $\mathcal{I} (V)$ de sommets 
$$
\{[\bigoplus_{i=1}^n \O_F.\pi^{-a_i} e_i] \;|\; a_i \in\Z \}
$$
Posons 
$$
\forall i<j\;\; \a_{ij} ( \bigoplus_{i=1}^n \O_F.\pi^{-a_i} e_i) = a_i - a_j
$$
les fonctions racines sur l'appartement. 

 D'après la section \ref{rutoipm} de
l'appendice les simplexes maximaux dans $\mathcal{A}$ sont en
bijection avec les collections d'entiers relatifs
$$
(b_{ij})_{1\leq i < j \leq n}
$$
telles que 
$$
b_{ij} \in \Z \text{ et } \forall i<j<k \;\; b_{ik}\in \{ b_{ij}+b_{jk},
b_{ij} + b_{jk} +1 \}
$$
\`A $(b_{ij})_{i<j}$ est associé le simplexe dont les sommets sont
$$
\{ x\in \mathcal{A} \;|\; \forall i<j \; \a_{ij} (x)\in \{ b_{ij}, b_{ij}+1 \}\;\}
$$
La proposition qui suit dit que la manière dont s'ordonnent les
valuations des points de torsion en fonction de leur ordre, les $(\l_i^{(k)})_{1\leq i\leq n,
  k\geq 1}$,  déterminent un simplexe dans l'appartement qui est une
facette de $S$. 

\begin{prop}\label{desmuap}
Soient $(b_{ij})_{1\leq i <j\leq n}$, $b_{ij}\in \N$ tels que 
$$
\forall i<j \;\; \l_i^{(b_{ij}+1)} \geq \l_j > \l_i^{(b_{ij}+2)}
$$
Alors, $(b_{ij})_{1\leq i<j\leq n}$ définit un simplexe maximal dans
$\mathcal{A}$ dont $S$ est une facette égale à l'intersection des murs
$$
\{x\;|\; \a_{ij} (x)=b_{ij}\} \;\text{ où } i<j \text{ et } \l_i^{(b_{ij}+1)}=\l_j
$$
avec ce simplexe maximal.
\end{prop}
\dem
Soient $i<j<k$. On a les inégalités 
$$
\l_j^{(b_{jk}+1)} \geq \l_k >\l_j^{(b_{jk}+2)}
$$
et 
$$
\l_i^{(b_{ij}+1)} \geq \l_j > \l_i^{(b_{ij}+2)}
$$
Mais d'après la remarque \ref{voutikj} (énoncée après cette démonstration) on déduit 
des deux inégalités précédentes 
$$
\l_i^{(b_{ij}+b_{jk}+1)} \geq \l_k > \l_j^{(b_{ij}+b_{jk} + 3)}
$$
Ce qui implique que $b_{ik}\in \{ b_{ij} + b_{jk},  b_{ij} + b_{jk}+1
\}$.  On en déduit donc d'après la section \ref{rutoipm} de
l'appendice que la donnée des $(b_{ij})_{i<j}$ définit bien un
simplexe maximal dans $\mathcal{A}$. 

Soit maintenant $\mu\in ]0,+\infty [$ suffisamment petit.
Supposons par exemple que $\forall 1\leq i <j\leq n \;\;
\mu <\l_i^{(b_{ij}+2)}$. Utilisons la formule donnée dans la section
\ref{ghuji} pour la filtration de ramification : 
$$
\text{si } \forall i\; k(i,\mu) \text{ est tel que }
\;\l_i^{(k(i,\mu))} \geq \mu >\l_i^{(k(i,\mu)+1)}
$$
alors 
$$
\Fil_\mu V = \bigoplus_{i=1}^n \pi^{- k(i,\mu)} \O_F .e_i 
$$
Soient maintenant $1\leq i< j\leq n$. De l'inégalité 
$$
\l_i^{(b_{ij}+1)} \geq \l_j > \l_i^{(b_{ij}+2)}
$$
on déduit en appliquant de nouveau la propriété de la remarque \ref{voutikj}
et l'inégalité $k(i,\mu)\geq b_{ij}+2$ 
que 
$$
\mu > \l_i^{(k(i,\mu)+1)} \geq \l_j^{(k(i,\mu) - b_{ij}+1)}
$$
et 
$$
\l_j^{(k(i,\mu)-b_{ij}-1)} > \l_i^{(k(i,\mu))} \geq \mu
$$
ce qui implique que 
$$
k(j,\mu) \in \{ k(i,\mu)-b_{ij}-1, k(i,\mu)-b_{ij} \}
$$
c'est à dire $k(i,\mu) - k(j,\mu) \in \{ b_{ij},b_{ij}+1 \}$ et donc
$[\Fil_\mu V]$ est dans la chambre définie par $(b_{ij})_{i<j}$.

De même si $\l_i^{(b_{ij}+1)} = \l_j$ on en déduit que $k(i,\mu) -
k(j,\mu) =b_{ij}$ ce qui implique l'assertion sur le fait que
$[\Fil_\mu V]$ est dans l'intersection  des murs donnés dans la
proposition.

Réciproquement on doit vérifier qu'étant donnés des entiers strictement positifs $(k_i)_i$ tels que 
$$
\forall i<j \; k_i - k_j \in \{ b_{ij},b_{ij}+1 \} \text{ et si }
\l_i^{(b_{ij})+1} =\l_j \text{ alors } k_i -k_j= b_{ij} 
$$
il existe $\mu$ tel que $[\Fil_\mu V] = <\pi^{-k_1} e_1,\dots,
\pi^{-k_n} e_n>$ (nécessairement pour un tel $\mu$ $\;[\Fil_\mu V]\in S$ car $\forall i\; k_i>0$). 
 Cela équivaut à trouver  $\mu$ tel que 
$$
\forall i \;\; \l_i^{(k_i)} \geq \mu > \l_i^{(k_i+1)}
$$
Soit l'intervalle de $\R$ $I_i= ] \l_i^{(k_i+1)}, \l_i^{(k_i)} ]$. On vérifie à partir des
hypothèses, de la définition des $(b_{ij})_{i,j}$  et de la remarque \ref{voutikj} 
que 
$$
\forall i<j\;\; \l_i^{(k_i)}>\l_j^{(k_j+1)} \;\text{ et } \;
\l_j^{(k_j)}> \l_i^{(k_i+1)}
$$
Cela implique que $\forall  i\neq j\; I_i \cap I_j \neq \emptyset$. On
conclu avec le lemme qui suit.
\qed

\begin{lemm}
Soit $(I_i)_{1\leq i\leq n}$ une famille d'intervalles de $\R$. Alors
$$
\bigcap_{i=1}^n I_i \neq \emptyset \ssi \forall i\neq j\; I_i \cap I_j
\neq \emptyset
$$
\end{lemm}

\begin{rema}\label{voutikj}
 La propriété suivante est utilisée à maintes reprises
$$
\forall \a,\b\;\forall a,b,c\geq 1\;\; \l_{\a}^{(a)} \geq \l_\b^{(b)}
\limpl \l_\a^{(a+c)} \geq \l_\b^{(b+c)}
$$
et de même en remplaçant $\geq $ par $>$.
\end{rema}

\subsubsection{Lien entre la facette et la filtration de ramification}

\begin{prop}
On a l'égalité suivante 
$$
\forall \l\in ]0, +\infty [\;\;
\forall l\in \N \; l\geq 1 \;\; \Fil_\l V = \pi^{1-l} \Fil_{\l^{(l)}}
V + \Fil_\infty V
$$ 
\end{prop}
\dem
Utilisant l'équivalence 
$$
\forall k\geq 1 \forall \mu_1,\mu_2\in ]0,+\infty [
\;\; \mu_1\geq \mu_2 \lssi \mu_1^{(k)} \geq \mu_2^{(k)}
$$
on vérifie facilement avec les notations de la section \ref{ghuji}
que 
$$
\forall i\;\forall l\geq 1\;\; k(i,\l)   = \sup \{ k(i,\l^{(l)})+1-l,0 \}
$$
\qed

\begin{defi}
Soient deux sommets $x,y\in \mathcal{I} (V)$. On pose
$$
x\vee y = \{ [\La_1+\La_2]\;|\; [\La_1]=x,\; [\La_2]=y \;\}
$$
un ensemble fini de sommets de $\mathcal{I} (V)$ contenu dans l'enclos délimité par $x$ et $y$. 
\end{defi}

\begin{prop}\label{rughv}
Le simplexe $S$ ainsi que le sommet $[T_p (H)]=[\Fil_\infty V]$ déterminent
complètement les sommets de l'immeuble associés à la filtration de
ramification inférieure. Cet ensemble, contenu  dans l'enclos délimité
par $S$ et $[T_p (H)]$, est  égal à 
$$
\bigcup_{x\in S} (x\vee [\Fil_\infty V])
$$
\end{prop}
\dem
Il suffit d'appliquer la proposition précédente couplée au faits suivants 
$$
\forall \l\in ]0,+\infty [\;\; \underset{l\drt +\infty}{\lim }
\l^{(l)}  = 0 
$$
et 
$
\forall l\geq 1
$
l'application $\l\longmapsto \l^{(l)}$ est une bijection de $]0,+\infty[$ dans lui-même. 
\qed

\begin{figure}[htbp]
   \begin{center}
      \includegraphics{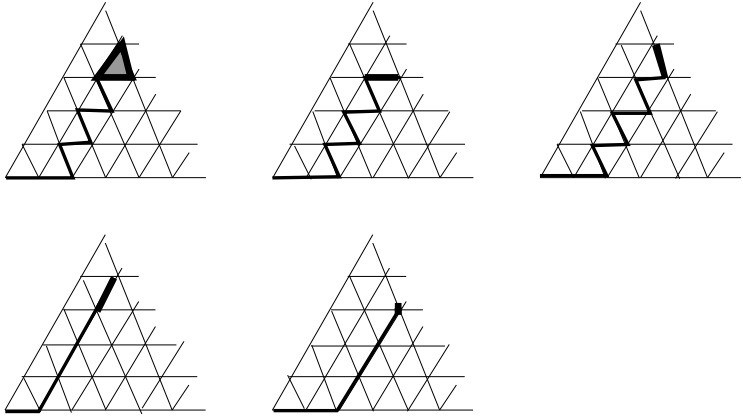}
   \end{center}
   \caption{\footnotesize Les quatre configurations possibles pour la
     filtration de ramification inférieure dans le cas de $\GL_3$
     }
   \label{config_ram}
\end{figure}

\section{Filtration de ramification supérieure}

On continu d'utiliser les notations de la section précédente.

La filtration étudiée précédemment se comporte bien vis à vis de la
restriction à un sous-schéma en groupes, par exemple $\forall \l\;
\forall n\leq m\;\; H[\pi^n](K)_\l = H[\pi^m] (K)_\l \cap H[\pi^n]
(K)$. 
Il s'agit de l'analogue de la filtration de ramification inférieure
du groupe de Galois d'une extension de corps locaux (\cite{Serre1}). 

La filtration que nous allons étudier dans cette section est adéquat
aux isogénies, et donc plus adaptée à l'action de $\GL_n (F)$ par les correspondances de Hecke sur les
espaces de Lubin-Tate.

\begin{defi}
On note pour tout $n\in\N^*$ et $\l\in ]0,+\infty]\;\;$ $H[\pi^n]^\l$
la filtration d'Abbes-Saito de $H[\pi^n]$ (cf. appendice
\ref{Ramitofm}). 
\end{defi}

Il s'agit d'une filtration décroissante qui d'après la proposition
\ref{zuytok} de l'appendice \ref{Ramitofm} 
 vérifie $\forall n\leq m
\;\; \pi^{m-n} H[\pi^m]^\l = H[\pi^n]^\l$.

\begin{defi}
On note $\psi=\eta^{-1}$ la fonction de Herbrand de $H[\pi]$
(cf. appendice \ref{Ramitofm}). 
\end{defi}

On a donc 
$$H[\pi^k]^\l (K) = H[\pi^k](K)_{\psi^k (\l)}$$
Rappelons (cf. section \ref{lapelctose}) que la fonction $\psi^k$ est également
celle notée $\l \longmapsto \l^{(k+1)}$ dans les sections précédentes.

\begin{defi}
Soit pour $\l\in ]0,+\infty [ $ 
$$
\Fil^{\,\l} V_p (H) = \underset{k\geq 1}{\limp} H[\pi^k]^\l (K) \subset
T_p (H)
$$
\end{defi}

On a donc en termes de la filtration de ramification inférieure et de
la fonction de Herbrand 
$$
\Fil^{\, \l} V = \underset{k\geq 1}{\limp} \left (  \Fil_{ \l^{(k+1)} } V \cap
 \pi^{-k} \Fil_\infty V \right )/
\Fil_\infty V 
$$
qui définit une filtration décroissante telle que $\Fil^\l V = T_p
(H)$ pour $\l$ suffisamment petit. 

\begin{lemm}
Soit $\l\in ]0,+\infty]$. Alors pour $k>>0$ 
$$
H[\pi^{k+1}]^\l = \{ x\in H[\pi^{k+1}] \;|\; \pi x \in H[\pi^k]^\l \}
$$
\end{lemm} 
\dem 
Il suffit de choisir $k$ tel que $\l^{(k+1)} = \frac{1}{q^n}
\l^{(k)}$. 

\qed

D'après le lemme précédent les $\Fil^{\,\l} V$ forment donc des
réseaux dans $V$ et on a 
$$
\forall \l\;\; \Fil^{\, \l} V = \pi^k \Fil_{\l^{(k+1)}} V \cap T_p (H)
\;\text{ pour } k>>0 
$$
On en déduit la propositions suivante

\begin{prop}
La filtration $(\Fil^{\,\l} V)_{\l \in ]0,+\infty]}$ forme une chaîne
périodique de
réseaux dans $T_p (H)$ pour $\l>>0$ 
 telle que le simplexe associé dans l'immeuble de
$\text{PGL} (V)$ soit le simplexe noté $S$ dans la définition
\ref{quatoki} de la section
précédente. 
\end{prop}

On a une description explicite de la filtration de ramification
supérieure analogue à celle donnée dans la section \ref{ghuji} pour la
filtration de ramification inférieure :

\begin{prop}
Soit $(e_i)_{1\leq i \leq n}$ une base adaptée du module de Tate $T_p
(H)$. Posons
$$
\forall i\;\;\forall \mu\in ]0,+\infty [\;\;\; l(i,\mu) = \inf \{
l\geq 0\;|\; \l_i \geq \mu^{(l+2)}\;\}
$$
Alors
$$
\Fil^{\,\mu} V = \bigoplus_{i=1}^n \O_F. \pi^{l(i,\mu)} e_i
$$
\end{prop}
\dem
Il suffit d'écrire pour $\l\in \N$ 
\begin{eqnarray*}
\pi^l e_i \in \Fil^{\,\mu} V &\lssi& \text{pour } k>>0 \;\pi^{l-k} e_i
\in \Fil_{\mu^{(k+1)}} V \\ &\lssi& \text{pour } k>>0 \;\l_i^{(k-l)}
\geq \mu^{(k+1)} \\
&\lssi & \l_i \geq \mu^{(l+2)}
\end{eqnarray*}
\qed

\begin{defi}
Soient $x,y$ deux sommets de $\mathcal{I} (V)$. On pose 
$$
x\wedge y = \{[\La_1\cap \La_2 ]\;|\; [\La_1]=x,\; [\La_2]=y\;\}
$$
\end{defi}

Voici la proposition cousine de la proposition \ref{rughv}

\begin{prop}
Le simplexe $S$ ainsi que le sommet $[\Fil_\infty V]$ déterminent
complètement les sommets de l'immeuble associés à la filtration de
ramification supérieure. Cet ensemble, contenu dans l'enclos
délimité  par $[\Fil_\infty V]$ et $S$, est égal à
$$
\bigcup_{x\in S} (x\wedge [\Fil_\infty V])
$$
\end{prop}

\section{Le point de l'immeuble associé à l'application de Hodge-Tate
  duale}\label{pohtmlsu}

\subsection{Une formule générale}

Soit $H$ comme dans les chapitres précédents.
Soit 
$$
\a_{H^D}^\O : T_p (H^\vee)\ldrt \omega_H \simeq \O_K
$$
et 
$$
\a_{H^\vee}^\O(-1) : T_p (H)^* \ldrt \omega_H (-1) \simeq \O_K
$$
où $F(1)$ désigne le module de Tate d'un groupe de Lubin-Tate de
$\O$-hauteur $1$. 

Considérons la norme ``additive'' associée
\begin{eqnarray*}
\|.\| : T_p (H)^* & \ldrt & \R\cup \{+\infty \} \\
\ph & \longmapsto & v (\a_{H^D} (-1) (\ph))
\end{eqnarray*}
Soit $(\Fil_\l V_p (H))_{\l\in]0,+\infty ]}$ la filtration de
ramification inférieure telle que définie dans la section \ref{knuypot}.

\begin{prop}
Soit $\mu\in ]0,+\infty]$. 
Soit $\ph \in T_p (H)^* \setminus \{ 0\}$ tel que $\ph (\Fil_{\mu}
V_p (H)) = \O_F$. Soit  $M=\ker \ph
\otimes \Qp \subset V_p (H)$. Soit $\Fil_\l M = \Fil_\l V_p (H) \cap
M$ la filtration induite sur $M$. Alors 
$$
\|\ph\| =
\int^{\mu}_{0} \left (|\Fil_\mu V_p (H)/\Fil_\infty V_p (H)|. |\Fil_\l M /\Fil_{\mu} M | -1
\right ) d\l +   v (\Fil_{\mu}V_p (H) /\Fil_\infty V_p (H) )
$$ 
\end{prop}
\dem
Soit  $\widetilde{\ph} =\ph\otimes Id : T_p (H) \otimes F/\O_F \ldrt F/\O_F$. 
D'après le théorème \ref{tklopui}
$$
\| \ph \| = \sum_{x\in \ker \widetilde{\ph} \setminus \{ 0 \}} v(x) 
$$
De plus  il y a une suite exacte 
$$
0 \ldrt \Fil_{\mu}V /\Fil_\infty V \ldrt \ker \widetilde{\ph} \ldrt
M/ \Fil_{\mu} M \ldrt 0
$$
Soient $(\l_k)_{k\geq 0}$, $\l_{k+1}<\l_k$ les sauts de la filtration
$(\Fil_\l V_p (H))_{\l \geq \mu}$ (avec $\l_0=\mu$). Alors
\begin{eqnarray*}
 \sum_{x\in \ker \widetilde{\ph} \setminus (\Fil_\mu V / \Fil_\infty
   V) } v(x)  &=&
 \sum_{k=1}^{+\infty} 
\left ( |(\Fil_{\l_k} M+\Fil_\mu V ) /\Fil_{\infty} V | - |(\Fil_{\l_{k-1}}
  M+ \Fil_\mu V ) /\Fil_{\infty}  V | \right ) \l_k \\
&=&  \sum_{k=1}^{+\infty} 
\left (( |(\Fil_{\l_k} M +\Fil_\mu V ) /\Fil_{\infty} V |-1) - (|\Fil_{\l_{k-1}}
  M+\Fil_\mu V /\Fil_{\infty} V |-1) \right ) \l_k \\
&=& \sum_{k=1}^{+\infty} \left ( |\Fil_\mu V/\Fil_\infty V|. 
|(\Fil_{\l_k} M + \Fil_\mu V ) /\Fil_\mu V | -1
\right )(\l_k - \l_{k+1}) \\
&=& \int^{\mu}_0 \left ( |\Fil_\mu V/\Fil_\infty V |. 
|\Fil_\l M/\Fil_\mu M | -1 \right ) d
\l
\end{eqnarray*}
\qed

\begin{coro}
Soit $\l_n$ la plus petite valuation des points de
$H[\pi](\O_K)$. Soit $\mu\in ]0,+\infty [$. Soit $\e$ tel que 
$\e \leq \mu$ et $0<\e\leq \l_n$. 
Il existe alors une constante $C$
telle que $\forall \ph \in V_p (H)^*$ tel que $\ph (\Fil_\mu V_p (H)) =\O_F$ si
$M=\ker \ph\otimes \Qp$ on ait  
$$
\| \ph \|= \Delta_\mu \left[ \int^\mu_\e |\Fil_\l M/\Fil_\infty M | d\l + \frac{q}{q-1}
\int^\e_{\e/q^n} |\Fil_\l M/\Fil_\infty M | d \l \right ] + C
$$
où $\Delta_\mu= |\Fil_\mu V_p (H)/\Fil_\infty V_p (H)|$. 
\end{coro}
\dem
Il suffit d'utiliser le fait que $\forall \l \leq \e$ on a
$\Fil_{\l/q^n} M = \pi^{-1} \Fil_\l M$ et que donc
\begin{eqnarray*}
\int^\e_0 |\Fil_\l M/\Fil_\infty M | d\l &=& \sum_{k=0}^{+\infty}
\int^{\e/q^{nk}}_{\e/q^{n (k+1)}} |\Fil_\l M/\Fil_\infty M| d \l \\
&=& \sum_{k=0}^{+\infty} \frac{1}{q^{nk}} \int^\e_{\e/q^n}
\underbrace{|\Fil_{\l/q^{nk}} M /\Fil_\infty M |}_{|\pi^{-k} \Fil_\l
  M/\Fil_\infty M| = q^{(n-1)k} |\Fil_\l M /\Fil_\infty M| } d \l \\
&=& \frac{q}{q-1} \int^\e_{\e/q^n} |\Fil_\l M/\Fil_\infty M| d\l 
\end{eqnarray*}
\qed

\begin{defi}
Soit $\La$ un réseau dans $V_p (H)$. On note $\La^\vee = \{\ph\in V_p
(H)^*\;|\; \ph (\La) \subset \O_F \}$ le réseau dual. On note alors 
$\left ( \Fil_\l V_p (H)^* \right )_{\l \in ]0,+\infty ] }$ la
filtration croissante de $V_p(H)^*$ égale à $\left ( (\Fil_\l V_p
  (H))^\vee\right )_{\l\in ]0,+\infty ]}$. On note $\|.\|_\l$ la norme
``additive'' sur $V_p (H)^*$
associée au réseau $\Fil_\l V_p (H)^*$. 
\end{defi}

\begin{theo}\label{Fokkomp}
Soit la fonction
\begin{eqnarray*}
 f: ] 0,+\infty ] & \ldrt & \R \\
\l & \longmapsto & |\Fil_\l V_p (H) / \Fil_\infty V_p  (H) |
\end{eqnarray*}
Soient $\mu$ et $\e$  comme dans le corollaire
précédent. Il existe alors une constante $C$ telle que 
\begin{eqnarray*}
\forall \ph \in V_p (H)^*\setminus \{ 0 \} \;\;\; \|\ph\| &=& \int^\mu_0
f(\l) \, q^{\| \ph
    \|_\l - \| \ph \|_\mu} d \l + \| \ph \|_\mu + C \\
&=& \int^\mu_\e f(\l) \, q^{\|\ph \|_\l - \| \ph \|_\mu} d\l +
\frac{q}{q-1} \int^\e_{\e/q^n} f(\l)\, q^{\| \ph \|_\l - \| \ph
    \|_\mu} d \l + \| \ph\|_\mu +  C
\end{eqnarray*}
\end{theo}
\dem 
Lorsqu'on remplace $\ph$ par $\pi^k \ph$ les membres des formules
données sont tous translatés de $k$. On peut donc supposer que $\ph
(\Fil_\mu V_p (H))=\O_F$ c'est à dire $\|\ph\|_\mu =0$. 

Il suffit alors de remarquer que si $M=\ker \ph\otimes \Qp$ on a une
suite exacte
$$
0\ldrt \Fil_\l M/\Fil_\mu M \ldrt \Fil_\l V_p (H)/\Fil_\mu V_p
(H) \ldrt \ph \left ( \Fil_\l V_p (H) \right )/\O_F \ldrt 0
$$
et que $\ph \left ( \Fil_\l V_p (H) \right ) = \pi^{\|\ph \|_\l }
\O_F$. Cela implique que 
$$
|\Fil_\mu V_p (H)/\Fil_\infty V_p (H)| . 
| \Fil_\l M/\Fil_\mu M| = f(\l)\, q^{\|\ph\|_\l}
$$
d'où le résultat par application  du corollaire
précédent.
\qed

\begin{coro}\label{Kopolm}
Soit $\mu \in ]0,+\infty [$ tel que $\mu \leq \l_n$. Il existe alors
une constante $C$ telle que 
$$
\forall \ph \in V_p (H)^*\setminus \{ 0\}\;\; \boxed{
\|\ph\| = \frac{q}{q-1} \int^\mu_{\mu/q^n} f(\l)\, q^{\| \ph \|_\l - \| \ph
    \|_\mu} d \l + \| \ph\|_\mu +  C}
$$
\end{coro}

\subsection{Applications}

\begin{prop}\label{shibouzouk}
Soit $[\|.\|] \in |\mathcal{I} (V_p (H)^*)|$ le point de l'immeuble
associé à l'application de Hodge-Tate duale de $H$. Alors, $[\|.\|]$
appartient à la réalisation géométrique du simplexe $S$ associé à la
filtration de ramification de $H$.  
\end{prop}
\dem
Il résulte du corollaire \ref{Kopolm} que la  norme $\| .\|$
ne dépend que des normes associées à des réseaux dont la classe
d'homothétie est dans $S$. D'où le résultat.
\qed

\begin{prop}
Soient $\l_1 \geq \dots \geq \l_n$ les pentes du polygone de Newton
associé à $H$.
 Soit $(e_1,\dots, e_n)$ une base adaptée de $T_p
(H)$ et $(e_1^*,\dots ,e_n^*)$ la base duale associée. Soit 
$
\forall i<j\;\; 
\a_{ij}$ la racine $\text{diag} (t_1,\dots,t_n) \mapsto 
t_i - t_j$ dans la base précédente de $V_p
(H)^*$. 

 Si pour $i<j$ $\;\l_i^{(b_{ij}+1)}= \l_j$ alors
 le
point de Hodge-Tate dual appartient  au mur $\a_{ij} = - b_{ij}$
\end{prop}
\dem
C'est une conséquence de la description donnée du simplexe $S$ comme
intersection de demi-appartements et de murs (proposition
\ref{desmuap}). 
\qed

\begin{exem}\label{jjdfii}
Si $\l_1= \dots = \l_n$ c'est à dire le polygone de Newton de la
multiplication par $\pi$ est plat alors $[\|.\|] = [T_p (H)^*]$. 
\end{exem}

\begin{exem}
Si $a_1\geq \dots \geq a_n \geq 0$ sont des entiers et le polygone de Newton de $H$ vérifie 
$$
\forall i<j\;\; \l_i^{(a_i -a_j+1 )} =\l_j
$$
alors le point de Hodge-Tate de $H$ est un sommet dans
l'appartement. Dans une base comme dans la proposition précédente ce
sommet est $ [<\pi^{-a_1} e_1^*,\dots, \pi^{-a_n} e_n^*>]$.
\end{exem}

\subsection{Expression en termes de la filtration de ramification supérieure}

\begin{prop}\label{kujipopom}
Soit $\|.\|^\l$ la norme sur $V_p (H)^*$ 
associée au réseau $\left ( \Fil^{\,\l} V_p (H)\right )^\vee$. 
Pour $\mu>>0$ il existe une constante $C$ telle que 
$$
\forall \ph \in V_p (H)^* \setminus \{ 0\}\;\;
\| \ph \| = \frac{q}{q-1} \int^\mu_{\mu-1} q^{\|\ph \|^\l - \|\ph
  \|^\mu} d \l + \|\ph\|^\mu + C
$$
\end{prop}
\dem
Effectuer un changement de variables dans l'intégrale couplé à la
formule exprimant la fonction $\eta$ de Herbrand comme une intégrale.
\qed

\section{Une structure simpliciale sur le squelette de l'espace de
  Lubin-Tate}

\subsection{Rappels sur l'espace de Lubin-Tate}\label{ujikoln}

Nous reprenons les notations de \cite{Cellulaire}. Soit 
$$
\X\simeq \spf (\O_{\widehat{F^{nr}}} [[x_1,\dots,x_{n-1}]])
$$
l'espace de Lubin-Tate des déformations d'un $\O$-module $\pi$-divisible formel de
dimension $1$ et de hauteur $n$. On suppose les coordonnées
$(x_1,\dots,x_{n-1})$ normalisées de telle manière que sur une loi de
groupe formelle universelle le polygone de Newton de la multiplication
par $\pi$ soit l'enveloppe convexe des points
$(1,1)$,$(q^i,v(x_i))_{1\leq i \leq n-1}$ et $(q^n,0)$ (cf. la section
1.5 de \cite{Cellulaire}). 

Soit $\X^{an}$ l'espace de Berkovich associé, $\X^{an}\simeq
\mathring{\mathbb{B}}^{n-1}$. On appelle squelette de $\X^{an}$
l'espace topologique 
$$
S (\X^{an}) = ]0,+\infty ]^{n-1} \subset |\X^{an}|
$$
où l'inclusion est donnée par l'application qui à
$(r_1,\dots,r_{n-1})\in  ]0,+\infty ]^{n-1}$ associe la norme
``additive'' 
$$
\left | \sum_{\underline{\a} = (\a_1,\dots,\a_{n-1})\in \N^{n-1}}
  a_{\underline{\a}} \,x_1^{\a_1}\dots x_n^{\a_n} \right | = \underset{\underline{\a}}{\inf} \{ v(a_{\underline{\a}})+ \a_1
r_1 +\dots +\a_{n-1} r_{n-1} \}
$$
Il y a une rétraction
\begin{eqnarray*}
r: |\X^{an}| & \ldrt & S(\X^{an}) \\
|.| & \longmapsto & (|x_i|)_{1\leq i \leq n-1}
\end{eqnarray*}
où $|.|$ désigne une semi-norme ``additive''.

\subsection{L'espace de Newton}

\begin{defi}
On note $Newt$ l'ensemble des fonctions convexes $f: [1,q^n] \ldrt ]0,+\infty
[$ linéaires sur les segments $[q^i,q^{i+1}]$ pour $1\leq i\leq n-1$,
ne s'annulant pas sur $[1,q^n[$, 
 vérifiant $f(1)=1$ et $f(q^n)=0$.
\end{defi}

Une telle fonction $f\in Newt$ est décroissante et vérifie $\text{Im}
f = [0 1]$. On muni $Newt$ de la topologie de la convergence uniforme
sur les fonctions continues de $[1, q^n] $ à valeurs dans $[0 1]$. On
vérifie aussitôt le lemme suivant :

\begin{lemm}
L'ensemble $Newt$ est convexe dans l'espace des fonctions de $[1,q^n]$
à valeurs dans $\R$. 
\end{lemm}

Ainsi on a une ``structure affine'' sur $Newt$ au sens où l'on a une
notion de barycentre :
$$
\forall I\text{ fini  }\; \forall (\g_i)_i\in \R_+^I\; \text{ vérifiant }
\sum_i \g_i =1 \;\forall (f_i)_i\in Newt^I\;\; \sum_i \g_i f_i \in
Newt 
$$ 

Une description concrète de $Newt$ que l'on utilisera est la suivante
:
$$
Newt = \{ (\l_1,\dots,\l_n)\in \R^n\;|\; \l_1\geq \dots \geq \l_n >0
\text{ et } \sum_{i=1}^n (q^i-q^{i-1}) \l_i =1 \}
$$
où $\l_i$ désigne la pente du polygone de Newton entre les abscisses
$q^{i-1}$ et $q^i$. Dans cette description le barycentre dans $Newt$
 consiste à
 prendre le barycentre des pentes $(\l_i)_i$, la structure affine est
 induite par celle de $\R^n$. 

Il y a une application continue surjective 
\begin{eqnarray*}
S (\X^{an}) & \ldrt & Newt \\
(v_1,\dots,v_{n-1}) & \longmapsto & \text{Conv} ( (1,1),(q,v_1),\dots,
(q^{n-1}, v_{n-1}),(q^n,0)) 
\end{eqnarray*}
où Conv désigne l'enveloppe convexe. Pour toute extension valuée complète
$K|\widehat{F^{nr}}$ l'application composée
$$
\X^{an} (K) \ldrt |\X^{an}| \xrig{\; \;r \; \;}S (\X^{an}) \ldrt Newt
$$
est celle qui à une déformation $(H,\rho)$ avec $H$ un $\O$-module
$\pi$-divisible sur $\O_K$ associe le polygone de la multiplication
par $\pi$ sur une loi de groupe formelle associée.

\begin{rema}
Soit $D$ l'algèbre à division d'invariant $1/n$ telle que $D^\times$
agisse  sur l'espace
de Lubin-Tate (cf. \cite{Cellulaire}). L'application $\X^{an} (K) \ldrt Newt$ est $D^\times$-invariante.
\end{rema}

\begin{defi}
Soit un polygone dans $Newt$ donné par ses pentes $\l_1\geq \dots\geq 
\l_n$. On notera pour tout $k\in \Z$ et tout $i$ 
$\;\l_i^{(k)}$ le nombre donné par $\l_i^{(k)}= +\infty$ si $k\leq 0$,
$\l_i^{(1)} =\l_i$ et $\l_i^{(k+1)}$ est la plus grande pente de
l'enveloppe convexe de $(0,\l_i^{(k)})$ et du polygone dont on est parti.
\end{defi}

\subsection{Les opérateurs $(\pi^{-a_1},\dots,\pi^{-a_n})$ sur
  l'espace de Newton}

Dans cette section nous définissons et étudions des opérateurs qui
sont des sections des correspondances de Hecke non-ramifiées, sections 
définies 
uniquement 
au niveau des polygones de Newton.

\subsubsection{Définition}

\begin{prop}\label{komopo}
Soit $\ph: H_1 \ldrt H_2$ une isogénie entre $\O$-modules formels de
dimension $1$. Le polygone de Newton de la multiplication par $\pi$
sur $H_2$ ne dépend que de celui de la multiplication par $\pi$ sur
$H_1$ ainsi que de celui de $\ker \ph$ (on entend par là la collection
 des $v(x), x\in \ker \ph \setminus \{ 0 \}$, comptées avec multiplicité).
\end{prop}
\dem 
Soit $C =\{ x\in H_1 \;|\; \pi x\in \ker \ph \}$. Alors $C/\ker \ph
\simeq H_2 [\pi]$. Après un choix de bonnes coordonnées formelles la
composition d'isogénies de groupes formels
$$
H_1 \xrig{\; \a\; } H_1/\ker \ph \xrig{\; \b \;} H_1/C
$$
s'écrit
$$
\b^* T = g(T),\;\;\; \a^* (T) = f(T)
$$
où $f$ et $g$ sont deux polynômes unitaires dans $\O_K [T]$. 
Le polygone de Newton de $f$ est connu par hypothèse. Celui de $f\circ g$
également d'après les rappels du début de la section \ref{fbmqsd}. 

Le résultat est donc une conséquence de ce que si $f$ et $g$ sont deux polynômes
unitaires dans $\O_K [T]$ alors $Newt (f\circ g)^*=Newt (f)^*\circ Newt (g)^*$ (cf. le lemme \ref{korTo} de l'appendice
\ref{Ramitofm}) et que ces fonctions sont bijectives, donc $Newt (f\circ g)$ et $Newt (f)$
déterminent $ Newt (g)^*$ et donc $Newt (g)$. 
\qed

\begin{coro}\label{kujokipm}
Soit $(a_1,\dots, a_n)\in \N^n$. Soit $H$ un $\O$-module
$\pi$-divisible formel de
dimension $1$ et hauteur $n$ sur $\O_K$ avec $K$ algébriquement clos. Soit
$(e_1,\dots,e_n)$ une base adaptée  de son module de Tate (cf. définition \ref{kkddssee}). 
 Soit $C\subset H$ l'adhérence schématique du
sous-groupe fini $<\pi^{-a_1} \bar{e}_1,\dots, \pi^{-a_n} \bar{e}_n>\subset T_p
(H)\otimes F/\O_F$ de la fibre
générique de $H$ et $H'=H/C$. 
Alors le polygone de Newton associé à $H'$ ne dépend que de celui de
$H$ et des $(a_i)_{1\leq i\leq n}$.
\end{coro}
\dem
D'après les formules données dans la section \ref{ghuji} les
valuations des éléments de $C$ ne dépendent que des
$(\l_i^{(k)})_{1\leq i \leq n,k\geq 1}$ qui eux mêmes ne dépendent que
des $(\l_i)_{1\leq i \leq n-1}$. Le résultat est donc une conséquence
de la proposition précédente.
\qed
\\

On peut donc définir $\forall (a_1,\dots ,a_n)\in \N^n$ un opérateur 
$$
(\pi^{-a_1},\dots, \pi^{-a_n}): Newt \ldrt Newt
$$
tel que $\forall K$, via l'application $\X^{an} (K) \ldrt Newt$,
celui-ci soit donné comme dans le corollaire précédent par le quotient
par le sous-groupe adhérence schématique de $<\pi^{-a_1} \bar{e}_1, \dots, \pi^{-a_n} \bar{e}_n>$.

\subsubsection{Une formule explicite pour les isogénies}\label{abbhui}

\subsubsection{Un lemme stupide}

\begin{lemm}
Soit $R$ un anneau et $F(X,Y)\in R[[X,Y]]$ une loi de groupe formelle commutative de dimension $1$. Soit $X\underset{F}{-} Y =F(X,\iota (Y))\in R[[X,Y]]$ où $\iota\in R[[T]]$ est l'inversion formelle associée à $F$ (l'unique série telle que $\iota (0)=0$ et $F(T,\iota (T))=0$. Alors 
$$
\exists h\in R[[X,Y]]^\times\;\;\; X\underset{F}{-} Y = (X-Y)\times h
$$
\end{lemm}
\dem
Soit l'anneau $Y$-adique $A=R[[Y]]$ et $X\underset{F}{-} Y \in A[[X]]$. Appliquons 
le lemme de division de Weierstrass à cette série et au polynôme $X-Y\in A[X]$.
$$
\exists h\in A[[X]]\;\; \exists a\in A\;\;\; X\underset{F}{-} Y = (X-Y) h + a
$$
Mais si l'on fait $X=Y$ on obtient $a=0$. De plus $X\underset{F}{-} Y \equiv X-Y \text{ mod deg }2$. Donc $h(0)=1$ et $h$ est inversible.
\qed

\begin{rema}
Le lemme précédent est faux si l'on remplace le couple $(X\underset{F}{-} Y,X-Y)$ par $(X\underset{F}{+} Y, X+Y)$. 
\end{rema}

\begin{coro}
Soit $\ph : H_1\ldrt H_2$ une isogénie de $\O$-modules formels de dimension $1$ sur $\O_K$ ($K$ algébriquement clos). Alors 
$$
\forall x\in H_1 (\O_K)\;\; v(\ph (x))= \sum_{a\in \ker \ph } v(x-a)
$$
où $x-a$ désigne la différence entre $x$ et $a$ dans le groupe $H_1 (\O_K)$.
\end{coro}
\dem
On sait qu'après un bon choix de coordonnées formelles à la source et au but l'isogénie s'écrit 
$$
\prod_{a\in \ker \ph} (T-a)
$$
Donc $v (\ph (x))= \sum_{a\in \ker \ph} v (x-a)$. On conclut d'après le lemme précédent.
\qed

\subsubsection{La formule}

Soient $(a_1,\dots,a_n)\in \N^n$ et $H$ sur $\O_K$ dont le polygone de
Newton associé a pour pentes $\l_1\geq \dots \l_n$. Soit
$(e_1,\dots,e_n)$ une base adaptée de $T_p (H)$. Rappelons que d'après la section \ref{ghuji} 
$$ 
\forall (x_1,\dots,x_n)\in F^n\;\; v(\sum_{i=1}^n x_i e_i) = \inf \{
\l_i^{(-v(x_i))}\;|\; 1\leq i\leq n \}
$$
Soit $C$ l'adhérence schématique du sous-groupe $<\pi^{-a_1} \bar{e}_1,\dots
\pi^{-a_n} \bar{e}_n>\subset T_p (H)\otimes F/\O_F$ et $\ph :
H\twoheadrightarrow H/C$· Soit $x=\sum_{i=1}^n x_i e_i$ dont on veut
calculer $v(\ph (x))$. On peut supposer que pour un sous-ensemble
$I\subset \{ 1,\dots ,n \}\;\; x=\sum_{i\in I} x_i e_i$ où $\forall
i\in I\; v(x_i)<-a_i$. Alors
\begin{eqnarray*}
v(\ph(x)) &=& \sum_{y\in \ker \ph} v(x+y) \\
&=& \sum_{(t_1,\dots, t_n)\in \O/\pi^{a_1}\oplus \dots \oplus
  \O/\pi^{a_n}} v\left (\sum_{i\in I} (x_i+t_i \pi^{-a_i}) e_i +
  \sum_{i\notin I} t_i \pi^{-a_i} e_i \right )  \\
&=&  \sum_{(t_1,\dots, t_n)\in \O/\pi^{a_1}\oplus \dots \oplus
  \O/\pi^{a_n}} \inf  \{ \l_i^{(-v(x_i))}\;|\; i\in I\} \cup \{
\l_i^{(a_i-v(t_i))}\;|\; i\notin I \} \\
&=& q^{\sum_{i\in I} a_i} \sum_{(k_i)_{i\notin I}  \atop {0 \leq k_i \leq a_i} } \prod_{i\notin I}
 A(k_i) \; \inf \{ \l_i^{(-v(x_i))} \; | \; i\in I \} \cup \{
\l_i^{( k_i)}\; |\; i\notin I \} 
\end{eqnarray*}
où l'on a posé
$$
A (k) = \left \{ q^k - q^{k-1} \text{ si } k\geq 1 \atop 
                  1 \;\;\;\;\;\;\;\;\;\;\;\;\;\,\text{       si } k=0
\right.
$$

\subsubsection{Composition des opérateurs $(\pi^{-a_1},\dots, \pi^{-a_n})$}

\begin{theo}\label{kakuyata}
Soient $H$ un $\O$-module $\pi$-divisible formel de dimension $1$ et
hauteur $n$ sur $\O_K$ tel que
le polygone de Newton associé ait pour pentes $\l_1\geq \dots \geq
\l_n$. Soit $(e_1,\dots,e_n)$ une base adaptée de $T_p (H)$.
 Soit $C$ l'adhérence schématique
de $<\pi^{-a_1} \bar{e}_1,\dots,\pi^{-a_n} \bar{e}_n> \subset T_p (H)\otimes
F/\O_F$ et $\ph : H\ldrt H/C = H'$. Soit $\ph_* : T_p (H)
\hookrightarrow T_p (H')$ le morphisme induit. Soit $\s\in
\mathfrak{S}_n$ une permutation telle que 
$$
\l_{\s(1)}^{(a_{\s(1)})} \geq \l_{\s (2)}^{(a_{\s(2)})} \geq \dots \geq
  \l_{\s (n)}^{(a_{\s (n)})}
$$
Alors, si l'on pose 
$$\forall i\in \{1,\dots, n\}\;\; \e_i = \ph_* ( \pi^{-a_{\s(i)} } e_{\s (i)})
$$
$(\e_i)_i$ est une base adaptée de $T_p (H')$. 
\end{theo}
\dem
Reprenons la formule donnée dans la section précédente 
\begin{eqnarray*}
v(\ph (x)) &=& q^{\sum_{i\in I} a_i} \sum_{(k_i)_{i\notin I}  \atop {0 \leq k_i \leq a_i} } (\prod_{i\notin I}
 A(k_i)) \; \inf \{ \l_i^{(-v(x_i))} \; | \; i\in I \} \cup \{
\l_i^{( k_i)}\; |\; i\notin I \} 
 \\
&=& \sum_{(k_i)_{1\leq i \leq n}  \atop {0 \leq k_i \leq a_i} }
(\prod_{1\leq i \leq n}
 A(k_i)) \; \inf \{ \l_i^{(-v(x_i))} \; | \; i\in I \} \cup \{
\l_i^{( k_i)}\; |\; 1\leq i \leq n \} 
\\
&=& \underset{i\in I}{\inf}\;  \sum_{(k_j)_{1\leq j \leq n}  \atop {0
    \leq k_j \leq a_j} } ( \prod_{1\leq j \leq n}
A(k_j )) \; \inf \{ \l_i^{(-v(x_i))}  \} \cup \{
\l_j^{( k_j)}\; |\; 1 \leq j \leq n \} \\
&=& \underset{i\in I}{\inf}\; v (\ph ( x_i e_i))
\end{eqnarray*}
et de plus $\forall i\;  v (\ph ( x_i e_i))$ ne dépend que de
$v(x_i)$. On a
$$
\l_{\s(1)}^{(a_{\s(1)}+1)} \geq \l_{\s (2)}^{(a_{\s(2)}+1)} \geq \dots \geq \l_{\s (n)}^{(a_{\s(n)}+1)}
$$
et donc grâce à la formule ci-dessus 
$$
v(\ph ( \pi^{-a_{\s(1)}-1} e_1)) \geq \dots \geq v( \ph (\pi^{-a_{\s
    (n)}-1} e_{\s (n)}))
$$
ce qui implique que cette suite de valuations est 
la suite des pentes du polygone de Newton associé à  $H'$. Reste à voir que 
$$
\forall i\;\forall k\geq 1 \;\forall y\in H'[\pi] \;\; 
v(\ph (\pi^{-a_i-k-1} e_i)) \geq v(\ph (\pi^{-a_i-k-1} e_i)+y)
$$
Mais cela résulte de la formule $v(\ph (x))=  \underset{i\in I}{\inf}\;
v (\ph ( x_i e_i))$.
\qed

\begin{coro}\label{compoheklp}
Soient $(a_1,\dots, a_n), (b_1,\dots, b_n)\in \N^n$. 
\begin{itemize}
\item
Pour tout
$\mathcal{P}\in Newt$ il existe une permutation $\s \in
\mathfrak{S}_n$ telle que 
$$
(\pi^{-a_1},\dots,\pi^{-a_n}). \left (
  (\pi^{-b_1},\dots,\pi^{-b_n}).\mathcal{P} \right ) = (\pi^{-(a_{\s
    (1)}+b_1)},\dots , \pi^{-(a_{\s (n)}+ b_n)}).\mathcal{P}
$$
\item Si $b_1 \leq b_2 \leq \dots \leq b_n$ on peut prendre $\s =Id$
\item En particulier il y a une action du monoïde 
$$
\{ (\pi^{-a_1},\dots,\pi^{-a_n})\; |\; 0 \leq a_1 \leq \dots \leq a_n \;\}
$$
sur $Newt$.
\end{itemize}
\end{coro}

\begin{rema}
Bien sûr $\forall b\in \N$ l'action de $(\pi^{-a_1},\dots,\pi^{-a_n})$
est la même que celle de $(\pi^{-(a_1+b)},\dots , \pi^{-(a_n+b)})$ ce
qui permet de définir une action du monoïde défini comme précédemment
en relâchant la condition $a_1 \geq 0$. 
\end{rema}

\subsection{Définition des simplexes maximaux en termes des pentes de Newton}

\begin{defi}
Soit $(e_1,\dots,e_n)$ la base canonique de $F^n$.
On note $\mathcal{A}$ l'appartement de l'immeuble de $\text{PGL}_n$
associé (cf. appendice \ref{kujitopo}) et $Q$ le quartier dont les sommets sont les
classes de réseaux 
$$
\{[<\pi^{a_1}e_1,\dots, \pi^{a_n} e_n>]\;|\; a_1\leq \dots \leq a_n \}
$$
On note $\text{pr}_Q : \mathcal{A} \ldrt Q$ la projection de
l'appartement sur son quartier (cf. la section \ref{Qdch} de
l'appendice \ref{kujitopo}). 
\end{defi}

Notons pour tout $i<j$ $\; \a_{ij} (<\pi^{a_1} e_1 , \dots ,
\pi^{a_n} e_n>) = a_j -a_i$. D'après la section \ref{rutoipm} de
l'appendice les simplexes maximaux dans $\mathcal{A}$ sont en
bijection avec les collections 
$$
(b_{ij})_{1\leq i < j \leq n}
$$
telles que 
$$
b_{ij} \in \Z \text{ et } \forall i<j<k \;\; b_{ik}\in \{ b_{ij}+b_{jk},
b_{ij} + b_{jk} +1 \}
$$
\`A $(b_{ij})_{i<j}$ est associé le simplexe dont les sommets sont
$$
\{ x\in Q\;|\; \forall i<j \; \a_{ij} (x)\in \{ b_{ij}, b_{ij}+1 \}\;\}
$$
Les simplexes maximaux dans $Q$ sont eux paramétrés par les
$(b_{ij})_{i<j}$ comme ci-dessus avec la condition supplémentaire
$\forall i<j\; b_{ij}\in \N$. 

\begin{defi}
Soit $S$ un simplexe maximal 
dans le quartier $Q$ associé à $(b_{ij})_{1\leq i <j \leq n}$.
Notons $Newt (S)$ le sous-ensemble de $Newt$ formé des polygones de
Newton dont les pentes satisfont 
$$
\forall i<j \;\;\; \l_i^{(b_{ij}+1)} \geq \l_j \geq \l_j^{(b_{ij}+2)}
$$
\end{defi}

\begin{lemm}
Les ensembles $Newt (S)$ sont convexes. De plus pour tout $i$ et $k$
l'application 
à valeurs dans $\R$ qui  à un polygone associé aux pentes $\l_1\geq
\dots \geq \l_n$ associe $\l_i^{(k)}$ est affine sur les $Newt (S)$.
\end{lemm}
\dem
Nous laissons au lecteur le soin de se convaincre de ce lemme en dessinant quelques polygones de Newton.
\qed

\begin{lemm}
Les convexes précédents recouvrent $Newt$ :
$$
Newt= \bigcup_S Newt (S)
$$
où $S$ parcourt les simplexes maximaux dans $Q$. 
\end{lemm}
\dem 
Soit un polygone de Newton associé aux pentes $\l_1\geq \dots \geq
\l_n$. Soient $(b_{ij})_{1\leq i < j \leq n}$, $b_{ij}\in \N$, tels
que 
$$
\forall i<j\;\; \l_i^{(b_{ij}+1)} \geq \l_j > \l_i^{(b_{ij}+2)}
$$
Soient $i<j<k$. On a les inégalités 
$$
\l_j^{(b_{jk}+1)} \geq \l_k > \l_j^{(b_{jk}+2)}
$$
et 
$$
\l_i^{(b_{ij}+1)} \geq \l_j > \l_i^{(b_{ij}+2)}
$$
Mais $$
\forall \a,\b\;\forall a,b,c\geq 0\;\; \l_{\a}^{(a)} \geq \l_\b^{(b)}
\limpl \l_\a^{(a+c)} \geq \l_\b^{(b+c)}
$$
(et de même avec $>$ au lieu de $\geq$) 
et donc
$$
\l_i^{(b_{ij}+b_{jk}+1)} \geq \l_k >\l_j^{(b_{ij}+b_{jk} + 3)}
$$
Ce qui implique que $b_{ik}\in \{ b_{ij} + b_{jk},  b_{ij} + b_{jk}+1
\}$. 
\qed

\begin{exem}\label{exwxssde}
Soit $S_0$ le ``simplexe fondamental'' de sommets 
$$
[<e_1,\dots, e_{i-1},\pi e_i ,\dots ,\pi e_n>] \;\;\; 1 \leq i \leq n
$$
Il est associé à $\forall i<j \; b_{ij} =0$. Mais en fait
$|S_0| = \{x\in Q\;|\; \a_{1 n } (x)\leq 1 \}$. Donc  
$$
Newt (S_0) = \{ \l_1^{(2)} \leq \l_n \} = \{ \frac{\l_1}{q^n} \leq
\l_n \}
$$

\end{exem}

\subsection{Lien avec l'espace de Lubin-Tate en niveau infini}

\begin{prop}
Soit $H$ un $\O$-module formel de dimension $1$ sur $\O_K$. Supposons
le muni d'une structure de niveau ``infinie'' 
$\eta : \O_F^n\iso T_p (H)$. Soit $S$ le simplexe de l'immeuble de $\text{PGL}
(V_p(H))$ construit dans la section \ref{kkjjuuii}, $S^\vee$ le
simplexe dual dans l'immeuble de $\text{PGL}
(V_p(H)^*)$ et $\eta^*S^\vee$ celui de l'immeuble de $\text{PGL}_n$
via l'isomorphisme $(F^n)^*\iso F^n$ donné par la base duale de la
base canonique. Alors, si $\mathcal{P}$ est le polygone de Newton
associé à $H$ et $S'$ un simplexe maximal dans le quartier $Q$
$$
\mathcal{P} \in Newt (S') \lssi \text{pr}_Q (\eta^* S^\vee)\subset S'
$$
où $\text{pr}_Q$ désigne la projection de l'immeuble sur le quartier $Q$.
\end{prop}
\dem
C'est une conséquence des résultats de la section \ref{kkjjuuii}.
\qed

\subsection{Action simpliciale des opérateurs
  $(\pi^{-a_1},\dots,\pi^{-a_n})$}

\begin{defi}\label{deffkkllcd}
On note $x\longmapsto (\pi^{a_1},\dots,\pi^{a_n}).x$ l'action par
translations sur $\mathcal{A}$ qui à  $[<\pi^{b_1} e_1,\dots, \pi^{b_n}
e_n>]$ associe $[<\pi^{a_1 +b_1} e_1,\dots , \pi^{a_n+b_n} e_n>]$. 
\end{defi}

\begin{theo}\label{koklmd}
Soit $S$ un simplexe maximal dans $Q$ et $(a_1,\dots,a_n)\in
\N^n$. Alors, l'opérateur $(\pi^{-a_1},\dots,\pi^{-a_n})$ induit une
bijection affine entre $Newt (S)$ et $Newt \left ( \text{pr}_Q
  \left ( (\pi^{a_1},\dots,\pi^{a_n}).S \right ) \right )$.
\end{theo}
\dem
Soit $\mathcal{P}\in S$ un polygone de Newton de pentes $\l_1\geq
\dots \geq \l_n$. Supposons $S$ associé aux $(b_{ij})_{i<j}$. Soit $H$
un $\O$-module formel de dimension $1$ ayant pour polygone de Newton
associé $\mathcal{P}$. Soit $(e_1,\dots,e_n)$ une base adaptée de
$T_p(H)$. Soit $\ph: H\ldrt H/C=H'$ l'isogénie associée à $(a_1,\dots
,a_n)$ comme dans le
corollaire \ref{kujokipm}. Soit $\s\in\mathfrak{S}_n$ une permutation
telle que 
$$
\l_{\s (1)}^{(a_{\s (1)})} \geq \dots \geq \l_{\s(n)}^{(a_{\s (n)})}
$$
Posons pour tout $i<j$
$$
b_{ij}'= a_{\s (j)} - a_{\s (i)} + b_{\s(i) \s (j)}
$$
où on a posé $b_{\s(i) \s (j)}= - b_{\s (j) \s (i)}$ si $\s(i)>\s
(j)$. Alors $\forall i<j\;\; b'_{ij}\in \N$. En effet, 
\begin{eqnarray*}
\text{si } \s (i) <\s (j) \;\; && \l_{\s  (i)}^{(b_{\s (i) \s (j)}+1)}
\geq \l_{\s (j)} \geq \l_{\s (i)}^{(b_{\s (i) \s (j)}+2)} \\
&\impl & \l_{\s(i)}^{(a_{\s(j)}+b_{\s(i)\s(j)})} \geq \l_{\s
  (j)}^{(a_{\s(j)})} \geq \l_{\s (i)}^{(a_{\s(j)}+ b_{\s(i)\s(j)}+1)} \\
 &\impl & a_{\s(i)} \leq \a_{\s (j)}+ b_{\s(i)\s(j)}
\end{eqnarray*}
\begin{eqnarray*}
\text{si } \s(i)>\s(j) \;\; && \l_{\s(j)}^{(b_{\s(j)\s(i)}+1)} \geq
\l_{\s (i)} \geq \l_{\s(j)}^{(b_{\s(j)\s(i)}+2)} \\
&& \l_{\s(j)}^{(a_{\s(i)}+b_{\s(j)\s(i)})} \geq
\l_{\s(i)}^{(a_{\s(i)})} \geq
\l_{\s(j)}^{(a_{\s(i)}+b_{\s(j)\s(i)}+1)} \\
&\impl & a_{\s(j)} \leq a_{\s(i)} +b_{\s(j) \s(i)} 
\end{eqnarray*}
Donc, $\s$  est l'unique élément du groupe de Weyl tel
que si $S'= (\pi^{a_1},\dots, \pi^{a_n}).S$
 est le simplexe dans $\mathcal{A}$ associé à la donnée
 $(b_{ij}+a_j-a_i)_{i<j}$ alors 
$$
\s.S'= \text{pr}_Q (S')
$$
Comme dans la démonstration du théorème \ref{kakuyata} on vérifie que
si
$\l'_1\geq \dots\geq \l'_n$ sont les pentes associées à $H'$ alors 
$$
\forall k\geq 1\;\; {\l'_i}^{(k)} = \sum_{0\leq k_j \leq a_{\s(j)} \atop
  1\leq j\leq n} \prod_{1\leq j \leq n} A(k_j)
\,\inf\{\l_{\s(i)}^{(a_{\s(i)}+k ) } \} \cup \{ \l_{\s(j)}^{(k_{j})}
  \;|\; j\neq i\}
$$
Cette formule permet de démontrer facilement que 
$$
\forall i<j\;\; {\l'_i}^{(b'_{ij}+1)} \geq \l'_j \geq   {\l'_i}^{(b'_{ij}+2)}
$$
De même elle permet de démontrer que sur $S$ l'application qui à un
polygone de pentes $\l_1\geq \dots \geq \l_n$ associe celui de pentes
$\l'_1\geq \dots \geq \l'_n$ est affine. En effet, sur $S$ les
applications qui à $(\l_1\geq \dots\geq \l_n)$ associent les quantités
  $\inf
\{\l_{\s(i)}^{(a_{\s(i)}+1)} \} \cup \{ \l_{\s(j)}^{(k_j)}\;|\; j\neq
i \}$ sont affines puisque l'ordre des $(\l_i^{(k)})_{1\leq i\leq n,
  k\geq 1}$ y reste inchangé.

Reste à voir que l'opérateur $(\pi^{-a_1},\dots,\pi^{-a_n})$ est un
isomorphisme. Mais pour $N\geq \underset{i}{\sup}\{ a_i\}$ 
sur $S$ la composée des opérateurs
$(\pi^{-a_1},\dots, \pi^{-a_n})$ et
$(\pi^{-(N-a_{\s(1)})},\dots,\pi^{-(N-a_{\s(n)})})$ correspond au
niveau des groupes de Lubin-Tate à $H\mapsto H/H[\pi^N]$ est est donc
l'identité. De même la composée sur $S'$ des opérateurs $(\pi^{-(N-a_{\s(1)})},\dots,\pi^{-(N-a_{\s(n)})})$
puis $(\pi^{-a_1},\dots, \pi^{-a_n})$ est l'identité.
\qed

\begin{rema}
Le lecteur effrayé par le fait que l'opérateur $(\pi^{-a_1},\dots
,\pi^{-a_n})$ sur $\text{Newt}$ agisse comme
$(\pi^{a_1},\dots,\pi^{a_n})$ sur l'appartement devrait penser au fait
que dans l'isomorphisme entre les tours de Lubin-Tate et de Drinfeld
l'action de $g\in \GL_n (F)$ d'un coté est transformée en celle de
$\,^t g$ de l'autre.
\end{rema}

\subsection{Définition des sommets}

\begin{defi}
Soit $x=[<\pi^{a_1},\dots,\pi^{a_n}>]$ un sommet dans
$\mathcal{A}$. On note
$\mathcal{P}(x)=(\pi^{-a_1},\dots,\pi^{-a_n}).\mathcal{P}_0$ où
$\mathcal{P}_0$ est le polygone de Newton ayant une seule pente 
$\l_1=\dots=\l_n=\frac{1}{q^n-1}$.
\end{defi}

\begin{prop}
Si $x=[<\pi^{a_1},\dots,\pi^{a_n}>]$ alors  $\mathcal{P}(x)$ est
l'unique polygone de pentes $\l_1\geq \dots\geq \l_n$ vérifiant
$$
\forall i<j\;\; \l_i^{(a_j-a_i+1)} =\l_j
$$
De plus 
$$
\forall i\;\; \l_i= \frac{q^{a_i}}{q^n-1} \sum_{0\leq k_j \leq a_j
  \atop j\neq i} \prod_{j\neq i} A(k_j) q^{-n (\sup( \{a_i+1 \}\cup
  \{k_j \;|\; j\neq i\})-1)}
$$
où $A(k)=q^k-q^{k-1}$ si $k>0$ et $A(0)=1$.
\end{prop}
\dem 
Le fait que les pentes de $\mathcal{P}(x)$ vérifient les égalités
annoncées résulte des formules données dans la section
\ref{abbhui} ainsi que dans la démonstration du théorème
\ref{kakuyata}. 

La démonstration concernant l'unicité est laissée au lecteur. Il
s'agit de vérifier qu'un certain système linéaire en les pentes
$(\l_i)_{1\leq i\leq n}$
possède une unique solution ce qui ne pose pas de problème particulier.
\qed

\begin{exem}\label{ootuiml}
Reprenons l'exemple \ref{exwxssde}. Le sommet $<e_1,\dots,e_n>$ correspond 
au polygone plat 
$$
\l_1= \dots= \l_n = \frac{1}{q^n-1}
$$
et pour $1\leq i\leq n-1$ le sommet $<e_1,\dots,e_i,\pi e_{i+1}, \dots,\pi e_n>$ correspond 
au polygone
$$
\l_1=\dots = \l_i = \frac{q^{n-i}}{q^n-1},\;\;\;\; \l_{i+1}=\dots = \l_n = \frac{1}{q^i (q^n-1)}
$$
\end{exem}

\subsection{Définition des chambres comme enveloppe convexe de
  sommets}

\begin{prop}\label{kkortuu}
Soit $S$ un simplexe maximal dans $\mathcal{A}$ vu comme la collection
de ses sommets. Alors $Newt (S)$ est l'enveloppe convexe des
$\mathcal{P}(x), \, x\in S$.
\end{prop}
\dem
Soit $S_0$ le simplexe de sommets 
$$
\{[<e_1,\dots,e_i,\pi e_{i+1},\dots,\pi e_n>]\;|\; 1\leq i\leq n\}
$$
Il existe $(a_1,\dots,a_n)\in \N^n$ tel que
$S=\text{proj}_Q ((\pi^{-a_1},\dots, \pi^{-a_n}).S_0)$. D'après le
théorème \ref{koklmd} il suffit de démontrer le résultat pour $S_0$. 

On a $Newt (S_0) = \{ \frac{\l_1}{q^n} \leq \l_n \}$. Donc
\begin{eqnarray*}
Newt(S_0)& =&\{ \l_1\geq \dots \geq \l_n>0\;|\; \sum_{i=1}^n (q^i
-q^{i-1})\l_i =1 \text{ et } \frac{\l_1}{q^n} \leq \l_n \} \\
&\simeq & \{ (\mu_1,\dots,\mu_n)\in\R_+^n\;|\; \sum_{i=1}^n \mu_i=1\}
\end{eqnarray*}
où l'on a posé $\mu_i =q^i(\l_i-\l_{i+1})$ pour $1\leq i\leq n-1$ et
$\mu_n=q^n\l_n-\l_1$ qui définit un changement de coordonnées
affines. Dans ces nouvelles coordonnées le simplexe est le simplexe standard de $\R_+^n$.  
 On vérifie alors 
que les sommets du simplexe standard dans ces
coordonnées correspondent aux pentes des $\mathcal{P}(x),x\in S_0$ telles que calculées dans l'exemple \ref{ootuiml}.
\qed
\\

\subsection{Définition des murs et demi-appartements}\label{muape}

\begin{defi}
Soit $b\in \N$ et $M$ le mur du quartier $Q$ égal  à 
$$
\{x\in Q\;|\; \a_{ij} (x) = b \}
$$
On définit alors un ``mur'' dans $Newt$ égal à 
$$
\{\mathcal{P}\in Newt \;|\; \l_i^{(b+1)} = \l_j \}
$$
Soit les demi-appartements dans $Q$ égaux à
$$
 \{ x\in Q\;|\; \a_{ij} (x) \geq b \}
$$
resp.
$$
 \{ x\in Q\;|\; \a_{ij} (x) \leq b \}
$$
On définit des ``demi-appartements'' associés dans $Newt$
$$
 \{ \mathcal{P}\in Newt \;|\; \l_i^{(b+1)} \geq  \l_j \}
$$
resp.
$$
 \{ \mathcal{P}\in Newt \;|\; \l_i^{(b+1)} \leq  \l_j \}
$$
\end{defi}

Ainsi avec les notations de la définition précédente si $S$ est un
simplexe maximal $Newt (S) $ est égal à l'intersection des demi-appartements le contenant.
De même si $x$ est un sommet de $Q$ alors $\{\mathcal{P}(x) \}$ est égal à
l'intersection des murs de $Newt$ contenant
$x$. 
\\

Dans la figure ci-dessous on a dessiné l'image réciproque de la
décomposition simpliciale de l'espace de Newton au squelette de
l'espace de Lubin-Tate dans le cas de $\GL_3$. Ce squelette est
l'ensemble des $(v(x_1),v(x_2))\in ]0,+\infty]^2$. La partie en haut à
droite correspond à la zone où le polygone de Newton est plat. Les
bandes horizontales et verticales correspondent à des zones où le
polygone de Newton ne détermine par de façon unique $v(x_1)$ et
$v(x_2)$. On  voit dans ce dessin le quartier dont l'origine, le
carré en
haut à droite,  correspond au polygone de Newton plat.

\begin{figure}
   \begin{center}
\begin{picture}(420,320)
\put(0,0){\vector(1,0){430}}
\put(412,-10){$v(x_1)$}
\put(0,0){\vector(0,1){330}}
\put(-25,295){$v(x_2)$}
\path(0,0)(325,190)
\path(325,190)(325,320)
\path(325,190)(420,190)
\path(325,190)(246,0)
\path(285,94)(420,94)
\path(265,46)(420,46)
\path(256,23)(420,23)
\path(251,11)(420,11)
\path(248,5)(420,5)
\path(247,2)(420,2)
\path(247,1)(420,1)
\path(246,0)(420,0)
\path(246,0)(420,0)
\path(246,0)(420,0)
\path(246,0)(420,0)
\path(229,134)(229,320)
\path(161,94)(161,320)
\path(113,66)(113,320)
\path(80,46)(80,320)
\path(56,33)(56,320)
\path(39,23)(39,320)
\path(27,16)(27,320)
\path(19,11)(19,320)
\path(13,8)(13,320)
\path(9,5)(9,320)
\path(6,4)(6,320)
\path(4,2)(4,320)
\path(3,2)(3,320)
\path(2,1)(2,320)
\path(285,94)(229,134)
\path(265,46)(201,66)
\path(201,66)(161,94)
\path(256,23)(187,33)
\path(187,33)(141,46)
\path(141,46)(113,66)
\path(251,11)(180,16)
\path(180,16)(131,23)
\path(131,23)(99,33)
\path(99,33)(80,46)
\path(248,5)(176,8)
\path(176,8)(126,11)
\path(126,11)(92,16)
\path(92,16)(70,23)
\path(70,23)(56,33)
\path(247,2)(175,4)
\path(175,4)(124,5)
\path(124,5)(89,8)
\path(89,8)(65,11)
\path(65,11)(49,16)
\path(49,16)(39,23)
\path(247,1)(174,2)
\path(174,2)(123,2)
\path(123,2)(87,4)
\path(87,4)(62,5)
\path(62,5)(46,8)
\path(46,8)(34,11)
\path(34,11)(27,16)
\path(246,0)(173,0)
\path(173,0)(122,1)
\path(122,1)(86,2)
\path(86,2)(61,2)
\path(61,2)(44,4)
\path(44,4)(32,5)
\path(32,5)(24,8)
\path(24,8)(19,11)
\path(246,0)(173,0)
\path(173,0)(122,0)
\path(122,0)(86,0)
\path(86,0)(61,1)
\path(61,1)(43,2)
\path(43,2)(31,2)
\path(31,2)(22,4)
\path(22,4)(17,5)
\path(17,5)(13,8)
\path(246,0)(173,0)
\path(173,0)(122,0)
\path(122,0)(86,0)
\path(86,0)(60,0)
\path(60,0)(43,0)
\path(43,0)(30,1)
\path(30,1)(21,2)
\path(21,2)(16,2)
\path(16,2)(12,4)
\path(12,4)(9,5)
\path(246,0)(173,0)
\path(173,0)(122,0)
\path(122,0)(86,0)
\path(86,0)(60,0)
\path(60,0)(42,0)
\path(42,0)(30,0)
\path(30,0)(21,0)
\path(21,0)(15,1)
\path(15,1)(11,2)
\path(11,2)(8,2)
\path(8,2)(6,4)
\path(229,134)(201,66)
\path(201,66)(187,33)
\path(161,94)(141,46)
\path(187,33)(180,16)
\path(141,46)(131,23)
\path(113,66)(99,33)
\path(180,16)(176,8)
\path(131,23)(126,11)
\path(99,33)(92,16)
\path(80,46)(70,23)
\path(176,8)(175,4)
\path(126,11)(124,5)
\path(92,16)(89,8)
\path(70,23)(65,11)
\path(56,33)(49,16)
\path(175,4)(174,2)
\path(124,5)(123,2)
\path(89,8)(87,4)
\path(65,11)(62,5)
\path(49,16)(46,8)
\path(39,23)(34,11)
\path(174,2)(173,0)
\path(123,2)(122,1)
\path(87,4)(86,2)
\path(62,5)(61,2)
\path(46,8)(44,4)
\path(34,11)(32,5)
\path(27,16)(24,8)
\path(173,0)(173,0)
\path(122,1)(122,0)
\path(86,2)(86,0)
\path(61,2)(61,1)
\path(44,4)(43,2)
\path(32,5)(31,2)
\path(24,8)(22,4)
\path(19,11)(17,5)
\path(173,0)(173,0)
\path(122,0)(122,0)
\path(86,0)(86,0)
\path(61,1)(60,0)
\path(43,2)(43,0)
\path(31,2)(30,1)
\path(22,4)(21,2)
\path(17,5)(16,2)
\path(13,8)(12,4)
\path(173,0)(173,0)
\path(122,0)(122,0)
\path(86,0)(86,0)
\path(60,0)(60,0)
\path(43,0)(42,0)
\path(30,1)(30,0)
\path(21,2)(21,0)
\path(16,2)(15,1)
\path(12,4)(11,2)
\path(9,5)(8,2)
\path(285,94)(201,66)
\path(265,46)(187,33)
\path(201,66)(141,46)
\path(256,23)(180,16)
\path(187,33)(131,23)
\path(141,46)(99,33)
\path(251,11)(176,8)
\path(180,16)(126,11)
\path(131,23)(92,16)
\path(99,33)(70,23)
\path(248,5)(175,4)
\path(176,8)(124,5)
\path(126,11)(89,8)
\path(92,16)(65,11)
\path(70,23)(49,16)
\path(247,2)(174,2)
\path(175,4)(123,2)
\path(124,5)(87,4)
\path(89,8)(62,5)
\path(65,11)(46,8)
\path(49,16)(34,11)
\path(247,1)(173,0)
\path(174,2)(122,1)
\path(123,2)(86,2)
\path(87,4)(61,2)
\path(62,5)(44,4)
\path(46,8)(32,5)
\path(34,11)(24,8)
\path(246,0)(173,0)
\path(173,0)(122,0)
\path(122,1)(86,0)
\path(86,2)(61,1)
\path(61,2)(43,2)
\path(44,4)(31,2)
\path(32,5)(22,4)
\path(24,8)(17,5)
\path(246,0)(173,0)
\path(173,0)(122,0)
\path(122,0)(86,0)
\path(86,0)(60,0)
\path(61,1)(43,0)
\path(43,2)(30,1)
\path(31,2)(21,2)
\path(22,4)(16,2)
\path(17,5)(12,4)
\path(246,0)(173,0)
\path(173,0)(122,0)
\path(122,0)(86,0)
\path(86,0)(60,0)
\path(60,0)(42,0)
\path(43,0)(30,0)
\path(30,1)(21,0)
\path(21,2)(15,1)
\path(16,2)(11,2)
\path(12,4)(8,2)
\end{picture}
\end{center}
\caption{\footnotesize La décomposition simpliciale de la boule
  $p$-adique ouverte de dimension $2$}
\end{figure}
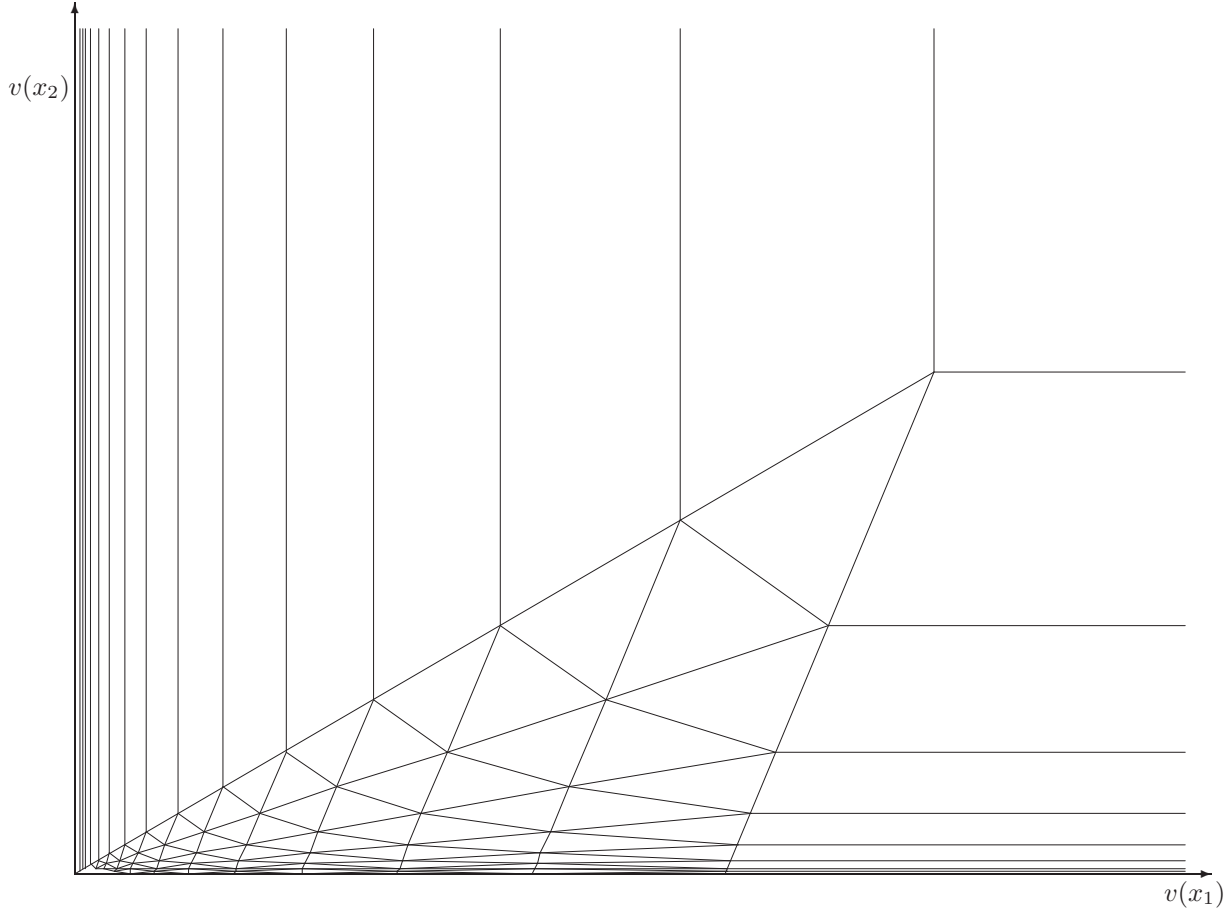
\newpage

\section{L'isomorphisme entre les tours de Lubin-Tate et de Drinfeld
  au niveau des squelettes après quotient par des sous-groupes
  compacts maximaux} \label{secotopp}

 Soient $|\M_\infty^{\LT,[0]}|$, resp. $|\M^{\D,[0]}_\infty|$, 
 l'espace de Berkovich associé à
 l'espace de Lubin-Tate, resp. l'espace de Drinfeld, en niveau infini
 (cf. \cite{Points} pour les notations). Ces espaces topologiques sont
 munis d'une action de $\GL_n(F)^1\times \O_D^\times$.
 D'après \cite{Points} il y a
 un homéomorphisme
$$
|\M_\infty^{\LT,[0]}|\iso |\M_\infty^{\D,[0]}|
$$
qui transforme l'action de $(g,d)\in \GL_n(F)^1\times \O_D^\times$ en
$(\,^t g ,d)$.  D'où en particulier 
$$
\xymatrix{
\GL_n(\O_F)\times \O_D^\times \bc \left |\M_\infty^{\LT,[0]} \right | \ar[r]^\sim
\ar@{=}[d] &
\GL_n(\O_F)\times \O_D^\times \bc \left  |\M_\infty^{\D,[0]} \right |  \ar@{=}[d] \\
\O_D^\times \bc | \M^{\LT,[0]} | \ar[r]^\sim & \GL_n (\O_F) \bc |
\Omega |
}
$$
où $|\M^{\LT,[0]}| $ est l'espace topologique noté $|\X^{an}|$ au
début dans la section \ref{ujikoln} et $\Omega$ l'espace symétrique
 de Drinfeld. Nous voulons décrire cette application au niveau des
 squelettes via la rétraction de ces espaces sur leur squelette : 
$$
\xymatrix{
  | \X^{an} | \ar[r] \ar[d]^r & \GL_n (\O_F) \bc
 |\Omega | \ar[d]^r \\
 S (\X^{an}) \ar@{..>}[r]^(.4){?} & \GL_n (\O_F) \bc
S(\Omega)
}
$$
où $S(\X^{an})$ a été décrit dans la section \ref{ujikoln},
$S(\Omega) = |\mathcal{I}(F^n)|$ est la réalisation géométrique de l'immeuble de $\text{PGL}_n
(F)$ et la rétraction $r: | \Omega| \ldrt S(\Omega) =|\mathcal{I}
(F^n)|$ est l'application qui à $[F^n\hookrightarrow K] \in \Omega (K)$
associe la norme sur $F^n$ déduite de la valuation sur $K$. 

Rappelons (cf. section \ref{ujikoln}) qu'il y a une application
$\O_D^\times$-invariante 
$S(\X^{an})\ldrt Newt$.

\begin{theo}
Il y a une application
$$
 S (\X^{an}) \xrig{\; \ph \;} \GL_n (\O_F) \bc
S(\Omega) 
$$
faisant commuter le diagramme précédent. Soit $Q$ le quartier 
de l'immeuble de $\text{PGL}_n (F)$ dont les
sommets sont les classes de réseaux $[<\pi^{a_1}e_1,\dots,\pi^{a_n}>]$
où $a_1\leq \dots \leq a_n$ et $(e_i)_i$ est la base canonique de
$F^n$. Alors l'application précédente se factorise en 
$$
\xymatrix{
 S(\X^{an}) \ar[rrr]^\ph \ar@{->>}[rd] &&& \GL_n (\O_F)
\bc S (\Omega) \\
& Newt \ar[r]^{\;\psi}_{\;\sim} & Q \ar[ru]^\simeq  
}
$$
où $\psi$ est une bijection simpliciale affine par morceaux lorsque
$Newt$ est muni de la structure simpliciale définie précédemment. Plus
précisément, pour tout sommet $x\in Q$ $\; \psi (\mathcal{P}(x))=x$,
pour toute chambre $S$ dans $Q$ $\; \psi (Newt(S))=S$ et 
$$
\psi_{|Newt(S)} :Newt(S) \iso S
$$
est affine.

Dans cette bijection les murs et demi-appartements
de $Newt$ tels que définis dans la section \ref{muape} correspondent à ceux du
quartier $Q$. 

De plus dans cette bijection l'action de l'opérateur
$(\pi^{-a_1},\dots,\pi^{-a_n}): \text{Newt}\ldrt \text{Newt}$
correspond à l'opérateur composé 
$$
Q\hookrightarrow |\mathcal{A}|\xrig{(\pi^{a_1},\dots,\pi^{a_n})}
|\mathcal{A}| \xrig{\; pr_Q\;} Q
$$
où $|\mathcal{A}|\xrig{(\pi^{a_1},\dots,\pi^{a_n})}
|\mathcal{A}|$ est une translation dans la partie vectorielle du groupe de Weyl affine sur
l'appartement et $pr_Q$ est la projection de l'appartement sur le
quartier, obtenue par application d'éléments du groupe de Weyl
vectoriel $\mathfrak{S}_n$, $Q=\mathfrak{S}_n\bc |\mathcal{A}|$. 
\end{theo}
\dem
L'application $|\X^{an}| \ldrt \GL_n (\O_F) \bc
|\Omega|$ se décrit de la façon suivante : soit $K|\breve{F}$ et $H$
un $\O$-module formel de dimension $1$ sur $\O_K$. Soit l'application
de Hodge-Tate de $H^\vee$ tordue, $\a_{H^\vee}^\O (-1): T_p (H)^* \ldrt \omega_H\otimes
K(-1)$. Fixons un isomorphisme $\O_F^n\iso T_p (H)^*$. Alors le point
de $\GL_n (\O_F)\bc \Omega (K)$ associé est la $\GL_n (\O_F)$-orbite
de $[\O_F^n \iso T_p (H)^*
\hookrightarrow \omega_H\otimes K(-1)]\in \Omega (K)$. Le point dans
$\GL_n (\O_F)\bc |\mathcal{I} (F^n)|$ associé est obtenu en prenant la
norme associée sur $F^n$. 

Les résultats des sections \ref{Filtraminf} et \ref{pohtmlsu} montrent 
que cette application se factorise par $|\X^{an}|\ldrt Newt$ (le
polygone de Newton détermine la filtration de ramification inférieure
dans un quartier de l'immeuble et celle-ci détermine le point de
Hodge-Tate associé dans $|\mathcal{I} (V_p (H)^* )|$). Soit $\psi :
Newt \ldrt Q$ cette application factorisée. 

Rappelons qu'on a définit une action par translations sur
l'appartement associé à la base canonique de $F^n$ (cf. définition
\ref{deffkkllcd}). 
Montrons que $\forall (a_1,\dots,a_n)\in \N\;\forall \mathcal{P}\in
Newt \;\; \psi ((\pi^{-a_1},\dots,\pi^{-a_n}).\mathcal{P}) =
\text{pr}_Q (
(\pi^{a_1},\dots,\pi^{a_n}).\psi (\mathcal{P}))$ où $\text{pr}_Q$ est
la projection de l'appartement sur le quartier.  

Le morphisme $|\M_\infty^{\LT}|\ldrt |\M_\infty^\D|$ transforme
l'action de $g\in \GL_n (F)$ en $\,^t g^{-1}$. Il existe donc
$\tau\in\mathfrak{S}_n$ telle que $\psi
((\pi^{-a_1},\dots,\pi^{-a_n}).\mathcal{P}) = \text{pr}_Q ((
\pi^{a_{\tau (1)}},\dots, \pi^{a_{\tau (n)}}).\psi (\mathcal{P}))$. 

Mais si $H$ sur $\O_K$ a pour polygone $\mathcal{P}$ et
$(\e_1,\dots,\e_n)$ est une base adaptée de $T_p (H)$ de base duale
$(\e_1^*,\dots, \e_n^*)$ alors d'après la proposition \ref{shibouzouk} 
le point de Hodge-Tate associé est dans un simplexe $S$ dans 
le quartier $\{[<\pi^{b_1} \e_1^*,\dots, \pi^{b_n} \e_n^*>] \;|\;
b_1\leq \dots \leq b_n \}$. Si $f:H\twoheadrightarrow H'$ est l'isogénie associée à
$(a_1,\dots, a_n)$ et $(\e_1,\dots,\e_n)$ comme dans le corollaire
\ref{kujokipm} alors il existe $\s\in \mathfrak{S}_n$ tel que si 
$f_* : T_p (H) \hookrightarrow T_p (H')$ $\;(f_* \pi^{-a_{\s(1)}-1}
e_1,\dots, f_* \pi^{-a_{\s(n)}-1})$ soit une base adaptée de $T_p
(H')$.
De plus, d'après le théorème \ref{koklmd}, $\s$ est un élément du
groupe de Weyl vérifiant $\s.(
(\pi^{a_1},\dots,\pi^{a_n}). \text{proj}_Q (S) ) =
\text{pr}_Q ((\pi^{a_1},\dots,\pi^{a_n}). \text{proj}_Q (S))$ où
$\text{proj}_Q: |\mathcal{I} (F^n)|\twoheadrightarrow Q$ et
$\text{pr}_Q$ est la restriction de $\text{proj}_Q$ à l'appartement.
D'où le résultat. 
\\

D'après le théorème \ref{koklmd} il suffit donc maintenant de montrer
que si $S_0$ est la chambre de sommets 
$$
\{[<e_1,\dots,e_i,\pi e_{i+1},\dots,\pi e_n>]\;|\; 1\leq i\leq n\}
$$
alors $\psi$ induit un isomorphisme affine entre $Newt (S_0)$ et
$S_0$. Il suffit pour cela de montrer que $\psi$ restreinte à $S_0$ est une application
affine. En effet, supposons le vérifié. Alors 
d'après l'exemple \ref{jjdfii} si $0$ désigne
le sommet de $Q$ alors $\psi (\mathcal{P} (0))=0$. Donc, d'après la
compatibilité de $\psi$ aux opérateurs $(\pi^{-a_1},\dots,\pi^{-a_n})$
pour tout sommet $x$ de $S_0$, $\psi (\mathcal{P} (x))=x$ ce qui
implique 
d'après la proposition \ref{kkortuu} que $\psi_{|Newt (S_0)}$ est un
isomorphisme entre $Newt (S_0)$ et $S_0$. 

Montrons donc que $\psi_{|Newt (S_0)}$ est affine. Pour un $H$ sur $\O_K$ 
notons $(\e_{1,H},\dots,\e_{n,H})$  une base adaptée de $T_p
(H)$ et $\|.\|_H$ désigne la norme additive associée sur $T_p (H)^*$.
Il suffit alors de montrer que l'application qui à un polygone de Newton $\mathcal{P}$ associe 
$\|\e_{i,H}^*\|_H$  est affine en restriction à $Newt (S_0)$. Mais on
vérifie directement avec le théorème \ref{tklopui} que cette
application est 
$$
\mathcal{P} \longmapsto \frac{q}{q-1} \left (1- (\mathcal{P} (q^i) -
  \mathcal{P} (q^{i-1})) \right )
$$
\qed

\section{Quelques applications}

\subsection{Sous-groupes canoniques généralisés}

Soit $K|F$ un corps valué complet pour une valuation $v:K\ldrt \R\cup
\{+\infty \}$ prolongeant celle de $F$. Soit $H$ un $\O$-module formel
de dimension $1$ et hauteur $n$ sur $\O_K$.

\begin{defi}
Un sous-groupe canonique de rang $r$ et niveau $k$ est un
sous $\O_F/\pi^k \O_F$ module  
$$
M\subset H[\pi^k] [\O_{\overline{K}}]
$$
libre de rang $r$
vérifiant
$$
\forall x\in M\;\forall y\in H[\pi^k] [\O_{\overline{K}}]\setminus
M\;\; v(x)>v(y) 
$$
\end{defi}

S'il existe, un tel sous-groupe est unique; il s'agit en effet d'un
cran de la filtration de ramification inférieure sur  $H[\pi^k]
[\O_{\overline{K}}]$ égal à l'ensemble des $q^{kr}$-éléments de plus
grande valuation dans $ H[\pi^k] [\O_{\overline{K}}]$. Son unicité
implique qu'il est stable sous $\Gal (\overline{K}|K)$ et définit donc
un sous-groupe étale de $H[\pi^k]\otimes_{\O_K} K$. 
\\

\begin{lemm}
Soit $\mathcal{P}\in Newt$ le polygone de Newton associé à $H$. Le
groupe $H$ possède un sous-groupe canonique de rang $r$ et niveau $k$
ssi $\mathcal{P}$ appartient au demi-appartement de Newt associé au
demi-appartement ouvert du quartier $Q$ défini par 
$$
\{x\in Q\;|\; \a_{r,r+1} (x)>k-1 \}
$$
\end{lemm}
\dem
D'après les résultats de la section \ref{Filtraminf} l'existence d'un tel sous-groupe canonique est équivalente à l'inégalité $\l_r^{(k)}>\l_{r+1}$.
\qed
\\

\begin{rema}
En particulier l'existence d'un sous-groupe canonique de rang et niveau donné peut se lire
sur l'application de Hodge-Tate de $H^\vee$.
\end{rema}

\begin{figure}[htbp]
   \begin{center}
      \includegraphics{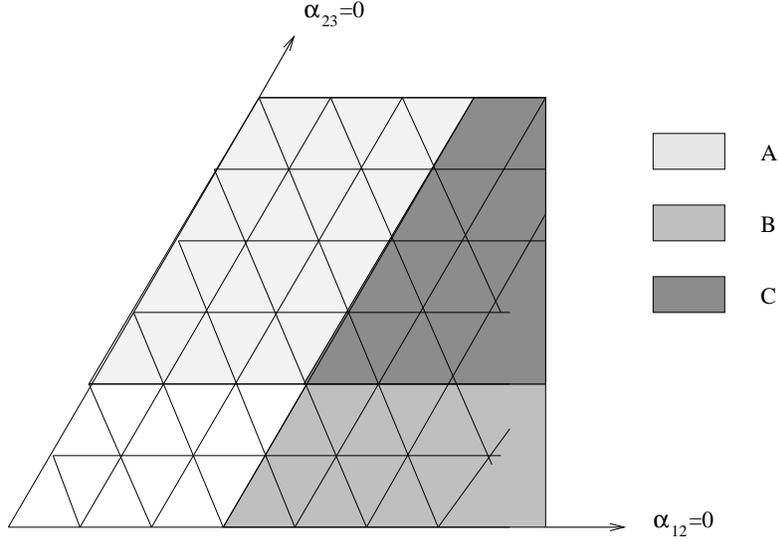}
   \end{center}
   \caption{\footnotesize Demi-appartement  $A=\exists$ sous-groupe canonique
     de rang $2$ et niveau $3$, $B=\exists$  sous-groupe canonique
     de rang $1$ et niveau $4$, $C=A\cap B$
     }
   \label{canonique_quartier}
\end{figure}

\begin{rema}
L'image réciproque dans l'espace de Lubin-Tate $\X^{rig}$ d'un demi-appartement ouvert définit un ouvert admissible de la boule ouverte.
\end{rema}

\begin{coro}
Pour tout $k\in \N^*\;\; \exists \e,0<\e<1,$  et des ouverts admissibles 
$(U_i^{(k)})_{1\leq i\leq n-1}$ dans l'espace de Lubin-Tate $\X^{rig}$ tels que 
$$
\left ( \bigcup_{i=1}^{n-1} U_i^{(k)} \right )^c\subset \mathbb{B} (0,\e)
$$
et $\forall x\in U_i^{(k)}$ le groupe $p$-divisible universel spécialisé en $x$  possède un sous-groupe canonique de rang $i$ et niveau $k$.
\end{coro}
\dem
Il suffit de vérifier que dans le quartier $Q$ pour $k\in \N^*$
$$
Q\setminus \bigcup_{1\leq i\leq n-1} \{ x\in Q\;|\; \a_{i,i+1} (x)>k-1 \}
$$

est un compact dans $Q$ et l'application $|\X^{an} |\twoheadrightarrow Newt$ est propre.
\qed
\\

On remarquera également le lemme suivant.

\begin{lemm}
Les assertions suivantes sont équivalentes : 
\begin{itemize}
\item $H$ possède un sous-groupe canonique de rang $r$ et niveau $k$
\item  $H$ possède un sous-groupe canonique de rang $r$ et niveau $1$ et si $H[\pi]^\mu \subset H[\pi]$ est le cran de la filtration de ramification  supérieure associé, $H[\pi^k]^\mu (\O_{\overline{K}})$ est libre de rang $r$ sur $\O/\pi^k \O$
\end{itemize}
Ainsi si c'est le cas $\forall 1\leq k'\leq k \;\; H[\pi^{k'}]^\mu (\O_{\overline{K}})$ est un sous-groupe canonique de rang $r$ et niveau $k'$ et il y a une suite exacte
$$
0\ldrt H[\pi^{k'}]^\mu (\O_{\overline{K}}) \ldrt H[\pi^k]^\mu (\O_{\overline{K}}) \xrig{\;\times \pi^{k'}\;} H[\pi^{k-k'}]^\mu (\O_{\overline{K}})\ldrt 0
$$
\end{lemm}
\dem
Elle ne pose pas de problème.
\qed

\subsection{Orbites de Hecke}

\begin{lemm}
Soit $H_1\ldrt H_2$ une isogénie de $\O$-modules formels de dimension
$1$. Soient $\mathcal{P}_1\in Newt$, resp. $\mathcal{P}_2\in Newt $, les
polygones de Newton associés à $H_1$, resp. $H_2$. Il existe alors 
$(a_1,\dots,a_n)\in \N^n$ tel que $\mathcal{P}_2 =
(\pi^{-a_1},\dots,\pi^{-a_n}).\mathcal{P}$. 
\end{lemm}
\dem
Toute isogénie peut s'écrire comme une composée d'isogénies cycliques
c'est à dire de noyau engendré par un point de $\pi$-torsion non-nul.
Pour une isogénie cyclique le résultat découle de la proposition
\ref{komopo}.
Le cas général est donc une conséquence du corollaire
\ref{compoheklp}.
\qed

\begin{defi}
Soit $\X$ l'espace de Lubin-Tate et $x\in |\X^{an}|$.  
Soit $H$ le groupe $p$-divisible universel sur $\X$,
  $K$ algébriquement clos tel que
$x$ provienne d'un élément $x'\in \X (\O_K)$ et  $H_{x'}$ la
spécialisation sur $\O_K$ de $H$. On appelle orbite de Hecke de
$x$ le sous-ensemble de $|\X^{an}|$ associé aux $y'\in \X(\O_K)$ tels
qu'il existe une isogénie $H_{x'}\ldrt H_{y'}$. 
\end{defi}

\begin{rema}
L'orbite de Hecke d'un point est un sous-ensemble $\O_D^\times$-invariant.
\end{rema}

\begin{prop}
Soit $x\in |\X^{an}|$. L'image dans $Newt$ de l'orbite de Hecke de $x$
ne dépend que de l'image de $x$ dans $Newt$. Via l'isomorphisme
$Newt\iso Q$ (où $Q$ est le quartier de l'immeuble défini dans les
section précédentes) cette image est l'orbite du point associé à $x$
dans $Q$ sous le groupe de Weyl affine (c'est à dire les points de l'orbite dans
l'appartement qui sont dans $Q$).  
\end{prop}
\dem
Soit $\mathcal{A}$ l'appartement et $\text{pr}_Q: \mathcal{A}\ldrt Q$
la projection. Soit $y\in Q$ le point associé à $x$. D'après le lemme
précédent l'image de l'orbite de Hecke de $x$ dans $Q$ est égale à
l'ensemble des 
$$
\text{pr}_Q (t.Q)
$$
où $t$ parcourt les translations dans le groupe de Weyl affine. 
Donc cette orbite est $W_{\text{aff}}.y\cap Q$.
\qed

\subsection{Domaines fondamentaux pour l'action des correspondances de
  Hecke}\label{kjyystsf286Mp}

Dans cette section nous montrons comment construire des domaines fondamentaux généralisant le domaine fondamental de Gross-Hopkins utilisé dans \cite{Cellulaire}. 
\\

L'orbite dans l'appartement $\mathcal{A}$ d'un point $y\in Q$ sous le groupe de Weyl affine admet la
description suivante. 
Soit $\e: \mathcal{A} \ldrt \Z/n\Z$ un étiquetage des sommets de
$\mathcal{A}$ donnée par $[\La]\longmapsto [\La_0:\La]\text{ mod }n$
où $\La_0$ est un réseau fixé. Si $S$ est une chambre de $\mathcal{A}$
(un simplexe maximal) on note $W(S)\subset \Aut (S)$ le groupe des
 rotations de $S$ de centre le barycentre de $|S|$, 
isomorphe à $\Z/n\Z$, défini par
$$
W(S) = \{ f\in \Aut (S)\;|\; \exists a \in \Z/n\Z\; \forall s\in S\;
\e (f(s))=\e (s) +a \;\}
$$
Alors,
$$
W(S) = \text{Stab}_{W_{\text{aff}}} (S)
$$
Soit donc maintenant $x\in |\mathcal{A}|$ et $S_0$ une chambre telle
 que $x\in |S_0|$. Soit $S$ une autre chambre
dans $\mathcal{A}$. Fixons un isomorphisme 
$$
\a: S_0\iso S
$$
tel que 
$$
\exists b\in \Z/n\Z\;\; \forall s\in S_0 \;\; \e (\a (s)) = \e (s) +b 
$$
Alors, si $|\a|: |S_0|\iso |S|$, on a 
$$
W_{\text{aff}}.x\cap |S| = W(S) . |\a| (x) = |\a| \left (
W(S_0).x \right )
$$

De là on déduit la proposition suivante 

\begin{prop}
Soit $S$ une chambre du quartier $Q$ et $D\subset |S|$ un
sous-ensemble tel que $W(S).D =|S|$. Supposons que $D$ est un polyèdre
rationnel dans le simplexe $|S|$, c'est à dire défini par des
inégalités linéaires à coefficients rationnels en coordonnées barycentriques dans
$|S|$. Alors, l'image réciproque à $|\X^{an}|$ de $D$ définit un
domaine analytique fermé $\mathcal{D}$ dans $\X^{an}$
associé à un ouvert admissible de $\X^{rig}$ tel que pour tout $K$
 et $x\in \X^{an} (K)$ il existe $K'|K$ finie, $y\in \mathcal{D}
(K')$ et une isogénie $H_x\ldrt H_y$. En d'autres termes les itérés de $\mathcal{D}$
par les opérateurs de Hecke sphériques  recouvrent $\X^{an}$.
\end{prop}
                 
Le problème est maintenant de comprendre comment se recollent les itérés d'un tel domaine $\mathcal{D}$ sous l'action des correspondances de Hecke. Cela peut se faire dans le simplexe fondamental du quartier comme dans la section 3 de \cite{Cellulaire}. 

Soit $S_0 = \text{convexe} ( 0,\omega_1,\dots,\omega_{n-1})$ le simplexe de $Q$ de sommets 
$$
0=<e_1,\dots,e_n>,\; \; \omega_i=<e_1,\dots,e_i,\pi e_{i+1},\dots, \pi e_n>\;\; 1\leq i \leq n-1
$$
Soit 
$$
\widetilde{S}_0 = S_0 \setminus \text{convexe} ( \omega_1,\dots, \omega_{n-1}) = S_0 \setminus \{ \a_{1n} =1 \}
$$
le simplexe privé du mur $\a_{1n}=1$. 
L'image réciproque à $Newt$ de $\widetilde{S}_0$ est
$$
Newt (\widetilde{S}_0 ) =\{\;\frac{\l_1}{q^n}< \l_n \;\} 
$$
Voici une généralisation de la proposition 3.5 de \cite{Cellulaire}. 

\begin{prop}
Soient $x,y\in \X^{an} (K)$ tels que les polygones de Newton associés appartiennent à $Newt(\widetilde{S}_0)$ et $\ph : H_x \ldrt H_y$ une isogénie ne se factorisant pas par $\pi$.
Alors $\ker \ph \subset H_x [\pi]$ et si $\dim_{\Fq} \ker \ph =r$ $\;\; \ker \ph$ est un sous-groupe canonique de rang $r$ et niveau $1$.
\end{prop}
\dem
Soit $M=\ker \ph$. Soient $\mu_1>\dots >\mu_r$ les pentes du polygone de Newton associé à $H_x$. On a des inclusions
$$
\xymatrix{
M[\pi^k]/M[\pi^{k-1}] \ar@{^(->}[r]^(.7){\times \pi^{k-1}} & M[\pi]  \ar@{^(->}[r] & H[\pi]
}
$$
De plus 
$$
\forall a\in M[\pi^k]\setminus M[\pi^{k-1}]\;\forall b\in M[\pi^{k-1}]\;\; v(a)<v(b)
$$
Posons
$$
\forall i, \; 1 \leq i\leq r, \forall k\geq 1\;\; d_i^{(k)} = \dim_{\Fq} \left ( \pi^{k-1} M[\pi^k] \cap \Fil_{\mu_i} M[\pi] \right )
$$
où l'on pose $\Fil_{\l} = \{ a\;|\; v(a)\geq \l \}$. 
On a donc $d_1^{(k)} \geq \dots \geq d_r^{(k)}$ et $\forall i\; d_i^{(k)}\geq d_i^{(k+1)}$.
\\
Les valuations des éléments non-nuls de $M$ sont
$$
\{ \frac{\mu_i}{q^{nk}}\;|\; k\geq 1,\; 1\leq i\leq r,\; d_i^{(k)} \neq 0 \} 
$$
et
$$
\forall i,k\;\; \text{long} ( \Fil_{\frac{\mu_i}{q^{nk}}} M ) = d_i^{(k)} |M[\pi^{k-1}]|
$$
Soit $(e_1,\dots, e_n)$ une base adaptée de $T_p (H)$. Il est alors aisé de vérifier qu'il existe $(a_1,\dots,a_n) \in \N^n$ tels que si $N=<\pi^{-a_1}\bar{e}_1 ,\dots,\pi^{-a_n} \bar{e}_n> \subset H_x[\pi^\infty]$ alors $M\simeq N$ comme groupes abstraits et 
$$
\forall \mu\; \; \text{long}\; \Fil_\mu N  = \text{long } \Fil_\mu M
$$
Il suffit pour cela de choisir les $(a_i)_i$ tels que
$$
\forall i,k \;\; |\{j\;|\; v(\pi^{-1} e_j)\geq \mu_i\text{ et } a_j \geq k\}|= d_i^{(k)}
$$
Soient donc $\mathcal{P}_x, \mathcal{P}_y \in Newt$ les polygones de Newton de $H_x$, resp. $H_y$. On  a donc d'après la proposition \ref{komopo}
$$
\mathcal{P}_y = (\pi^{-a_1},\dots,\pi^{-a_n}).\mathcal{P}_x
$$
Mais il est aisé de vérifier que 
$$
\forall (a_1,\dots,a_n)\in \N^n\;\; \inf \{ a_i \} =0 \;\; \text{pr}_Q ((\pi^{a_1},\dots, \pi^{a_n}).\widetilde{S}_0) \cap \widetilde{S}_0 \limpl \exists i\; (a_1,\dots,a_n) = (
\underbrace{1,\dots,1}_{i},0,\dots 0)
$$
De plus si $u\in \widetilde{S}_0$ vérifie 
$$
\text{pr}_Q ( \underbrace{( \pi,\dots,\pi)}_{i},1,\dots 1).u) \in \widetilde{S}_0
$$
alors $\a_{i,i+1} (u)\neq 0$ et on conclut facilement.
 \qed

\begin{coro}[Généralisation du domaine fondamental de Gross-Hopkins]
Soit $S_0$ le simplexe maximal du quartier $Q$ contenant
l'origine. Soit $D$ un polyèdre rationnel dans $|S_0|$ tel que
$D\subset |\widetilde{S}_0|$ c'est à dire $D$ ne rencontre pas le mur
$\{ \a_{1n} =1 \}$. Supposons que $$|S_0|= W(S_0).D$$
où $W(S_0) = \text{Stab}_{W_{\text{aff}}} (S_0)$. 
 Notons $\forall \s \in W (S_0) \simeq \Z/n\Z\;\; \partial_\s D= D\cap \s .D$. Soit 
$\mathcal{D}$ l'image réciproque à $|\X^{an}|$ de $D$ et $\forall s \in W(S_0)\;\; \partial_\s \mathcal{D}$ celle de $\partial_\s D$. Ce sont des domaines analytiques fermés associés à des ouverts admissibles de $\X^{rig}$. De plus les itérés de $\mathcal{D}$ sous les correspondances de Hecke sphériques forment un recouvrement admissible de $\X^{rig}$. Les itérés de $\mathcal{D}$ ayant une intersection non-vides avec $\mathcal{D}$ sont associés aux correspondances de Hecke $(\underbrace{\pi^{-1},\dots,\pi^{-1}}_{i},1,\dots,1)$ avec $1\leq i\leq n-1$. Ces intersections sont les $\partial_\s \mathcal{D}$ et les relations de faces sont obtenues par quotient par un sous-groupe canonique généralisé dans les points de $\pi$-torsion.
\end{coro}

\begin{coro}
Pour tout domaine fondamental $\mathcal{D}$ comme dans le corollaire précédent on peut refaire la construction de \cite{Cellulaire} : on peut reconstruire la tour de Lubin-Tate en niveau infini à partir d'un tel domaine fondamental. 
\end{coro}

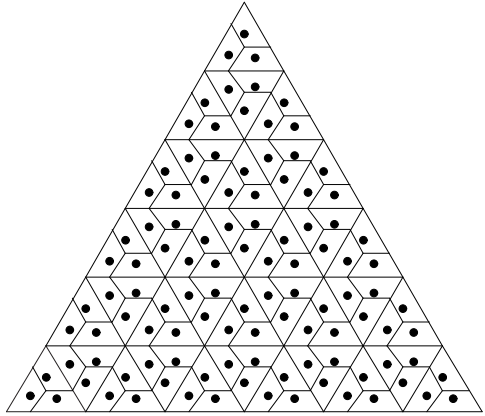
\begin{figure}[h]
   \begin{center}
\begin{picture}(300,200)
\path(0,0)(30,0)
\path(30,0)(15,25)
\path(0,0)(15,25)
\put(9,6){\circle*{3}}
\put(19,5){\circle*{3}}
\put(15,13){\circle*{3}}
\path(15,8)(9,0)
\path(15,8)(10,17)
\path(15,8)(25,8)
\path(30,0)(60,0)
\path(60,0)(45,25)
\path(30,0)(45,25)
\put(39,6){\circle*{3}}
\put(49,5){\circle*{3}}
\put(45,13){\circle*{3}}
\path(45,8)(39,0)
\path(45,8)(40,17)
\path(45,8)(55,8)
\path(15,25)(45,25)
\path(45,25)(30,51)
\path(15,25)(30,51)
\put(24,31){\circle*{3}}
\put(34,30){\circle*{3}}
\put(30,39){\circle*{3}}
\path(30,34)(24,25)
\path(30,34)(25,43)
\path(30,34)(40,34)
\path(60,0)(90,0)
\path(90,0)(75,25)
\path(60,0)(75,25)
\put(69,6){\circle*{3}}
\put(79,5){\circle*{3}}
\put(75,13){\circle*{3}}
\path(75,8)(69,0)
\path(75,8)(70,17)
\path(75,8)(85,8)
\path(45,25)(75,25)
\path(75,25)(60,51)
\path(45,25)(60,51)
\put(54,31){\circle*{3}}
\put(64,30){\circle*{3}}
\put(60,39){\circle*{3}}
\path(60,34)(54,25)
\path(60,34)(55,43)
\path(60,34)(70,34)
\path(30,51)(60,51)
\path(60,51)(45,77)
\path(30,51)(45,77)
\put(39,57){\circle*{3}}
\put(49,56){\circle*{3}}
\put(45,65){\circle*{3}}
\path(45,60)(39,51)
\path(45,60)(40,69)
\path(45,60)(55,60)
\path(90,0)(120,0)
\path(120,0)(105,25)
\path(90,0)(105,25)
\put(99,6){\circle*{3}}
\put(109,5){\circle*{3}}
\put(105,13){\circle*{3}}
\path(105,8)(99,0)
\path(105,8)(100,17)
\path(105,8)(115,8)
\path(75,25)(105,25)
\path(105,25)(90,51)
\path(75,25)(90,51)
\put(84,31){\circle*{3}}
\put(94,30){\circle*{3}}
\put(90,39){\circle*{3}}
\path(90,34)(84,25)
\path(90,34)(85,43)
\path(90,34)(100,34)
\path(60,51)(90,51)
\path(90,51)(75,77)
\path(60,51)(75,77)
\put(69,57){\circle*{3}}
\put(79,56){\circle*{3}}
\put(75,65){\circle*{3}}
\path(75,60)(69,51)
\path(75,60)(70,69)
\path(75,60)(85,60)
\path(45,77)(75,77)
\path(75,77)(60,103)
\path(45,77)(60,103)
\put(54,83){\circle*{3}}
\put(64,82){\circle*{3}}
\put(60,91){\circle*{3}}
\path(60,86)(54,77)
\path(60,86)(55,95)
\path(60,86)(70,86)
\path(120,0)(150,0)
\path(150,0)(135,25)
\path(120,0)(135,25)
\put(129,6){\circle*{3}}
\put(139,5){\circle*{3}}
\put(135,13){\circle*{3}}
\path(135,8)(129,0)
\path(135,8)(130,17)
\path(135,8)(145,8)
\path(105,25)(135,25)
\path(135,25)(120,51)
\path(105,25)(120,51)
\put(114,31){\circle*{3}}
\put(124,30){\circle*{3}}
\put(120,39){\circle*{3}}
\path(120,34)(114,25)
\path(120,34)(115,43)
\path(120,34)(130,34)
\path(90,51)(120,51)
\path(120,51)(105,77)
\path(90,51)(105,77)
\put(99,57){\circle*{3}}
\put(109,56){\circle*{3}}
\put(105,65){\circle*{3}}
\path(105,60)(99,51)
\path(105,60)(100,69)
\path(105,60)(115,60)
\path(75,77)(105,77)
\path(105,77)(90,103)
\path(75,77)(90,103)
\put(84,83){\circle*{3}}
\put(94,82){\circle*{3}}
\put(90,91){\circle*{3}}
\path(90,86)(84,77)
\path(90,86)(85,95)
\path(90,86)(100,86)
\path(60,103)(90,103)
\path(90,103)(75,129)
\path(60,103)(75,129)
\put(69,109){\circle*{3}}
\put(79,108){\circle*{3}}
\put(75,117){\circle*{3}}
\path(75,112)(69,103)
\path(75,112)(70,121)
\path(75,112)(85,112)
\path(150,0)(180,0)
\path(180,0)(165,25)
\path(150,0)(165,25)
\put(159,6){\circle*{3}}
\put(169,5){\circle*{3}}
\put(165,13){\circle*{3}}
\path(165,8)(159,0)
\path(165,8)(160,17)
\path(165,8)(175,8)
\path(135,25)(165,25)
\path(165,25)(150,51)
\path(135,25)(150,51)
\put(144,31){\circle*{3}}
\put(154,30){\circle*{3}}
\put(150,39){\circle*{3}}
\path(150,34)(144,25)
\path(150,34)(145,43)
\path(150,34)(160,34)
\path(120,51)(150,51)
\path(150,51)(135,77)
\path(120,51)(135,77)
\put(129,57){\circle*{3}}
\put(139,56){\circle*{3}}
\put(135,65){\circle*{3}}
\path(135,60)(129,51)
\path(135,60)(130,69)
\path(135,60)(145,60)
\path(105,77)(135,77)
\path(135,77)(120,103)
\path(105,77)(120,103)
\put(114,83){\circle*{3}}
\put(124,82){\circle*{3}}
\put(120,91){\circle*{3}}
\path(120,86)(114,77)
\path(120,86)(115,95)
\path(120,86)(130,86)
\path(90,103)(120,103)
\path(120,103)(105,129)
\path(90,103)(105,129)
\put(99,109){\circle*{3}}
\put(109,108){\circle*{3}}
\put(105,117){\circle*{3}}
\path(105,112)(99,103)
\path(105,112)(100,121)
\path(105,112)(115,112)
\path(75,129)(105,129)
\path(105,129)(90,155)
\path(75,129)(90,155)
\put(84,135){\circle*{3}}
\put(94,134){\circle*{3}}
\put(90,143){\circle*{3}}
\path(90,138)(84,129)
\path(90,138)(85,147)
\path(90,138)(100,138)
\path(30,17)(25,8)
\path(30,17)(24,25)
\path(30,17)(40,17)
\put(30,11){\circle*{3}}
\put(24,18){\circle*{3}}
\put(34,19){\circle*{3}}
\path(60,17)(55,8)
\path(60,17)(54,25)
\path(60,17)(70,17)
\put(60,11){\circle*{3}}
\put(54,18){\circle*{3}}
\put(64,19){\circle*{3}}
\path(45,43)(40,34)
\path(45,43)(39,51)
\path(45,43)(55,43)
\put(45,36){\circle*{3}}
\put(39,44){\circle*{3}}
\put(49,45){\circle*{3}}
\path(90,17)(85,8)
\path(90,17)(84,25)
\path(90,17)(100,17)
\put(90,11){\circle*{3}}
\put(84,18){\circle*{3}}
\put(94,19){\circle*{3}}
\path(75,43)(70,34)
\path(75,43)(69,51)
\path(75,43)(85,43)
\put(75,36){\circle*{3}}
\put(69,44){\circle*{3}}
\put(79,45){\circle*{3}}
\path(60,69)(55,60)
\path(60,69)(54,77)
\path(60,69)(70,69)
\put(60,62){\circle*{3}}
\put(54,70){\circle*{3}}
\put(64,71){\circle*{3}}
\path(120,17)(115,8)
\path(120,17)(114,25)
\path(120,17)(130,17)
\put(120,11){\circle*{3}}
\put(114,18){\circle*{3}}
\put(124,19){\circle*{3}}
\path(105,43)(100,34)
\path(105,43)(99,51)
\path(105,43)(115,43)
\put(105,36){\circle*{3}}
\put(99,44){\circle*{3}}
\put(109,45){\circle*{3}}
\path(90,69)(85,60)
\path(90,69)(84,77)
\path(90,69)(100,69)
\put(90,62){\circle*{3}}
\put(84,70){\circle*{3}}
\put(94,71){\circle*{3}}
\path(75,95)(70,86)
\path(75,95)(69,103)
\path(75,95)(85,95)
\put(75,88){\circle*{3}}
\put(69,96){\circle*{3}}
\put(79,97){\circle*{3}}
\path(150,17)(145,8)
\path(150,17)(144,25)
\path(150,17)(160,17)
\put(150,11){\circle*{3}}
\put(144,18){\circle*{3}}
\put(154,19){\circle*{3}}
\path(135,43)(130,34)
\path(135,43)(129,51)
\path(135,43)(145,43)
\put(135,36){\circle*{3}}
\put(129,44){\circle*{3}}
\put(139,45){\circle*{3}}
\path(120,69)(115,60)
\path(120,69)(114,77)
\path(120,69)(130,69)
\put(120,62){\circle*{3}}
\put(114,70){\circle*{3}}
\put(124,71){\circle*{3}}
\path(105,95)(100,86)
\path(105,95)(99,103)
\path(105,95)(115,95)
\put(105,88){\circle*{3}}
\put(99,96){\circle*{3}}
\put(109,97){\circle*{3}}
\path(90,121)(85,112)
\path(90,121)(84,129)
\path(90,121)(100,121)
\put(90,114){\circle*{3}}
\put(84,122){\circle*{3}}
\put(94,123){\circle*{3}}
\end{picture}
\caption{\footnotesize Une décomposition simpliciale du quartier relativement à un domaine fondamental pour l'action des correspondances de Hecke et une orbite de Hecke (points noirs)}
\end{center}
\end{figure}

\subsection{L'image du domaine fondamental de Gross-Hopkins dans 
  l'espace de Drinfeld}

Soit 
$$
\mathcal{D} = \{ (x_1,\dots,x_{n-1})\in \X^{an}\;|\; \forall i\;
v(x_i)\geq 1-\frac{i}{n} \} 
$$
le domaine fondamental de Gross-Hopkins (\cite{HopkinsGross},
\cite{Cellulaire}). 
Si pour $x\in \X^{an}$ on note $\mathcal{P}_x\in Newt$ le polygone de
Newton correspondant et si $\mathcal{P}_{GH}$ désigne le polygone
défini par $\forall i\; v(x_i)=1-\frac{i}{n}$ alors
$$
\mathcal{D} = \{ x\in \X^{an}\;|\; \mathcal{P}_x \geq \mathcal{P}_{GH}
\} 
$$
et $\mathcal{D}$ est donc l'image réciproque dans $\X^{an}$ de
l'ensemble 
$$
D=\{\mathcal{P}\; |\; \mathcal{P}\geq \mathcal{P}_{GH} \}
$$ 

On se propose de calculer l'image de ce domaine fondamental de
Gross-Hopkins dans l'espace de Drinfeld via l'isomorphisme entre les
deux tours.

\begin{prop}
L'ensemble convexe $D\subset Newt$ est un polytope possédant $2^{n-1}$
points extrémaux en bijection avec les parties $A\subset \{1,\dots,
n-1 \}$. \`A $A = \emptyset$ est associé le polygone plat
$\mathcal{P}_\emptyset$ défini par
$$
\l_1=\dots =\l_n = \frac{1}{q^n-1}
$$
Pour $A\neq \emptyset$ si $A=\{ i_1<\dots <i_r \}$ le polygone associé
$\mathcal{P}_A$ est défini par  
$$
\l_1=\dots = \l_{i_1}, \; \l_{i_1+1} = \dots = \l_{i_2},\; \dots,
\;\l_{i_{r-1}+1} = \dots = \l_{i_r}
$$
$$
\text{et }\; \forall k\in \{1,\dots,r\}\;\;\; v(x_{i_k}) = 1 - \frac{i_k}{n}
$$
Ainsi $\mathcal{P}_{\{1,\dots, n\}} = \mathcal{P}_{GH}$. On a donc 
$$
D=\text{Convexe} ( \mathcal{P}_A )_{A\subset \{1,\dots, n-1 \}}
$$
\end{prop}  

\begin{prop}
Soit pour tout $i, \;1\leq i\leq n-1$, le sommet $\omega_i
=<e_1,\dots,e_i,\pi e_{i+1},\dots, \pi e_n>$, $\mathcal{P}
(\omega_i)\in Newt$ le polygone associé (cf. exemple \ref{ootuiml}) et
$\mathcal{P}_0$ le polygone ``plat'' associé au sommet
$<e_1,\dots,e_n>$. Alors $\forall A=\{ i_1<\dots< i_r \}
\subset \{1,\dots, n-1 \}$, avec
les notations de la proposition précédente, 
$$
\mathcal{P}_A\in \text{Convexe} ( \mathcal{P}_0,\mathcal{P}
(\omega_{i_1}),
\dots, \mathcal{P} (\omega_{i_r}))
$$
plus précisément
$$
\mathcal{P}_A= a_0 \mathcal{P}_0 + \sum_{k=1}^r a_k \mathcal{P} (\omega_{i_k})
$$
où
$$
\forall k\neq 0\;\; a_k=\frac{q^{i_k}}{n} \left ( 
\frac{i_k-i_{k-1}}{q^{i_k}-q^{i_{k-1}}} -
\frac{i_{k+1}-i_k}{q^{i_{k+1}} -q^{i_k}}
\right )
$$
où l'on a posé $i_0=0,\; i_n=n$, et
$$
a_0= 1-\sum_{k\neq 0} a_k 
$$
En particulier
$$
\mathcal{P}_{\{1,\dots,n\}} = \frac{1}{n} \left ( \mathcal{P}_0 +
  \sum_{i=1}^{n-1} \mathcal{P} ( \omega_i) \right )
$$
\end{prop}

\begin{coro}
L'image du domaine fondamental de Gross-Hopkins dans le simplexe fondamental est l'enveloppe convexe de l'origine $[<e_1,\dots,e_n>]$ et des 
$$
a_0 [<e_1,\dots,e_n>] + \sum_{k=1}^r a_k \omega_{i_k}
$$
où $\{i_1<\dots<i_r\} \subset \{ 1,\dots, n\}$ 
et les $(a_i)_i$ sont comme dans la proposition précédente. Le sommet correspondant à 
$\{i_1<\dots<i_r\} = \{ 1,\dots, n\}$ est le barycentre du simplexe.

La bijection de \cite{Points} induit une bijection entre les points de la tour de Lubin-Tate en niveau infini au dessus du domaine fondamental de Gross-Hopkins et ceux de la tour de Drinfeld en niveau infini au dessus d'un point de $\Omega$ s'envoyant dans le polytope précédent via l'application quotient $|\Omega|\twoheadrightarrow |\mathcal{I} (\text{PGL}_n)| \twoheadrightarrow \GL_n (\O_F)  \bc |\mathcal{I} (\text{PGL}_n)|=
Q$. 
\end{coro}

\begin{figure}[htbp]
   \begin{center}
      \includegraphics{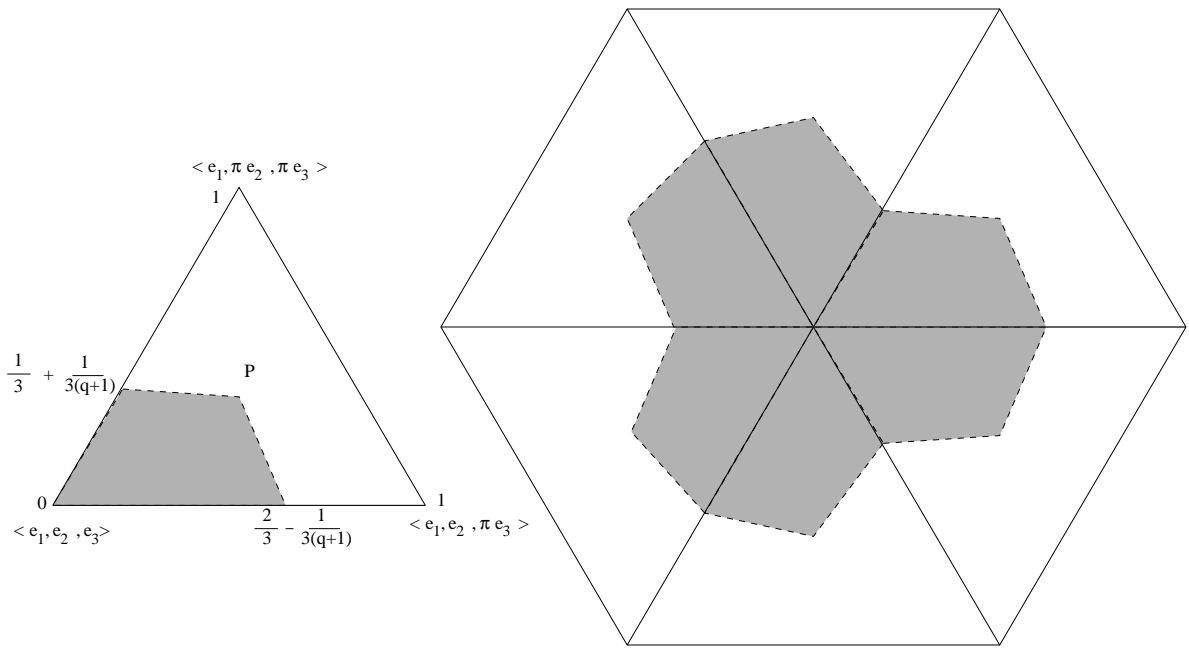}
   \end{center}
   \caption{\footnotesize L'image du domaine fondamental de
     Gross-Hopkins dans le ``simplexe fondamental''  du quartier et son image dans l'appartement
     }
   \label{image_grosshopkins}
\end{figure}

\subsection{Application au morphisme des périodes}

L'application des périodes de Gross-Hopkins a été étudiée en détails
dans \cite{Yu}. Dans cette section on
explique comment obtenir de nouveau les résultats de \cite{Yu} à partir des
calculs de la section 2 de \cite{Points}. 

Commençons par remarquer aves les notations de la section
\ref{kjyystsf286Mp}
que si $S_0$ est le simplexe fondamental dand le quartier $Q$ et
$\widetilde{S}_0$ désigne $S_0$ privé du mur $\{\a_{1n}=0\}$ alors
$$
\forall x\in |\widetilde{S}_0| \;\forall (a_1,\dots,a_n)\in
\N^n\setminus \N.(1,\dots,1) \;
n|\sum_{i=1}^n a_i\limpl (\pi^{a_1},\dots,\pi^{a_n}).x\notin |\widetilde{S}_0|
$$ 
et que donc si $U\subset \X^{an}$ désigne l'image réciproque dans
l'espace de Lubin-Tate de $|\widetilde{S}_0|$ 
$$
U= \{ \frac{\l_1}{q^n}<\l_n \}
$$
pour toute correspondance de Hecke $T$ non-triviale de degré un multiple de $n$ on
a $T.U\cap U=\emptyset$. On vient donc de donner une démonstration
géométrique sur l'immeuble de la proposition 4.2 de
\cite{Cellulaire}. On en déduit comme dans la section 4 de
\cite{Cellulaire} que l'application des périodes de Gross-Hopkins 
$\breve{\pi} : \X^{an}\ldrt \mathbb{P}^{n-1}$ induit un isomorphisme
entre $U$ et sont image. 

Intéressons-nous maintenant aux propriétés ``métriques'' de
l'application des périodes. On sait d'après le théorème 2.3 de
\cite{Cellulaire} que si $\breve{\pi} = [f_0:\dots : f_{n-1}]$ alors
$$
\forall \underline{x}\in U\;\; v\left (\frac{f_i (\underline{x})}{f_0
  (\underline{x})}\right ) = v (x_i)
$$
Mais maintenant si $y\in Q$ $\exists \; (a_1,\dots \a_n)\;\;
(\pi^{a_1},\dots , \pi^{a_n}).y \in |\widetilde{S}_0|$. De cela on
déduit que l'on a des estimations sur la position $\breve{\pi} (x)$
pour $x\in \X^{an}$ en fonction de l'image de $x$ dans le quartier
$Q$. En effet, si $y\in Q$ correspond à $x$ alors si
$z=(\pi^{a_1},\dots,\pi^{a_n}).y \in |\widetilde{S}_0|$ on a une
estimation sur la position de $ \Pi^{\sum_i a_i}. \breve{\pi} ( x)$
à un élément de $\O_D^\times$ près où $\Pi$ est une uniformisante de
l'alèbre à division $D$. 

On peut par exemple déduire de tout cela que si $S$ est un simplexe de
$Q$ et $U(|S|)\subset \X^{an}$ le domaine analytique fermé de l'espace
de Lubin-Tate associé alors $\breve{\pi} (U(|S|) =
\mathbb{P}^{n-1}$. De même on vérifie que l'application $\breve{\pi}$
restreinte à l'ouvert $U(\mathring{|S|})$ est finie, et on peut calculer son
degré. 

Nous laissons au lecteur la liberté de trouver d'autres applications.

\section{Le cas des espaces de Lubin-Tate en niveau Iwahori}

Nous avons précédemment étudié l'isomorphisme au niveau des squelettes induit par l'isomorphisme entre les tours de Lubin-Tate et Drinfeld après quotient par $\GL_n (\O_F)\times \O_D^\times$.
Nous allons faire de même après quotient par $I\times \O_D^\times$ où $I$ désigne un sous-groupe d'Iwahori. 
\\

Rappelons que si $(e_1;\dots,e_n)$ désigne la base canonique de $F^n$ alors $\mathcal{A}$ est l'appartement de sommets $$[<\pi^{a_1} e_1,\dots, \pi^{a_n} e_n>],\; (a_1,\dots,a_n)\in\Z^n$$ $Q$ le quartier de sommets 
$$
[<\pi^{a_1} e_1,\dots, \pi^{a_n} e_n>]\;\;\; a_1\leq \dots \leq a_n
$$
et $S_0$ le ``simplexe fondamental'' de sommets 
$$
[<e_1,\dots,e_n>],\;\;\; [<e_1,\dots, e_i, \pi e_{i+1},\dots ,\pi e_n>]\;\; 1\leq i\leq n-1
$$
Soit $I$ le sous-groupe d'Iwahori associé au simplexe $S_0$
$$
I=\{ g\in \GL_n (\O_F)\;|\; g \text{ mod } \pi = \left (
\begin{array}{ccc}
\times &0 & 0 \\
\times & \ddots & 0 \\
\times & \times & \times
\end{array}
\right ) \}
$$
Soit $I^{opp}$ le sous-groupe d'Iwahori 
$$
I^{opp} = \{ g\in \GL_n (\O_F)\;|\; g \text{ mod } \pi = \left (
\begin{array}{ccc}
\times & \times & \times \\
0 & \ddots & \times \\
0 & 0 & \times
\end{array}
\right ) \}
$$
Soit $\mathcal{Y}$ le schéma formel sur $\spf (\breve{\O})$ classifiant les classes d'isomorphismes d'objets 
$$
(H_1\xrig{\;\ph_1\;} H_2 \xrig{\;\ph_2\;}\dots \ldrt H_{n-1}\xrig{\ph_{n-1}} H_n \xrig{\;\ph_n\; } H_1,\rho)
$$
où $(H_1,\rho)$ est dans l'espace de Lubin-Tate sans niveau ($\rho$ étant la déformation), $\forall i\;\; \ph_i:H_i\ldrt H_{i+1}$ est une isogénie de degré $q$ entre $\O$-modules formels et $\ph_{n}\circ \dots \circ \ph_1 =\pi$. Le schéma formel $\mathcal{Y}$ est un modèle entier de l'espace de Lubin-Tate en niveau $I^{opp}$. Le morphisme d'oubli de la structure de niveau est 
$$
\mathcal{Y} \ni (H_1\drt \dots\drt H_n \drt H_1,\rho)\longmapsto (H_1,\rho)\in \X
$$
Il y a de plus un isomorphisme
$$
\mathcal{Y} \simeq \spf \left ( \breve{\O}[[y_1,\dots, y_n]]/(y_1\dots y_n -\pi)\right )
$$
où si $(H_1\drt \dots\drt H_n \drt H_1,\rho)$ est l'objets universel sur $\mathcal{Y}$ il y a un diagramme commutatif
$$
\xymatrix{
\Lie \,H_1 \ar[r]^{\Lie \,\ph_1} & \dots \ar[r] & \Lie\, H_i \ar[r]^{\Lie \, \ph_i} & \Lie\, H_{i+1} \ar[r] & \dots \ar[r] & \Lie\, H_n\ar[r]^{\Lie\, \ph_n} & \Lie\, H_1 \\ 
\O_\mathcal{Y} \ar[r]^{\times y_1} \ar[u]^\simeq & \dots \ar[r] & \O_\mathcal{Y} \ar[r]^{\times y_i} \ar[u]^\simeq & \O_\mathcal{Y} \ar[r]\ar[u]^\simeq  & \dots \ar[r] & \O_\mathcal{Y} \ar[u]^\simeq \ar[r]^{\times y_n} & \O_\mathcal{Y} \ar[u]^\simeq 
}
$$
Cela découle de la théorie de la déformation de Grothendieck-Messing.

\begin{defi}
Posons $\Delta = S(\mathcal{Y}^{an})$ le squelette de la fibre générique $\mathcal{Y}^{an}$ égal à 
$$
\Delta = \{ (v_1,\dots,v_n)\in ]0 \, 1[^n\;|\; \sum_{i=1}^n v_i =1 \}
$$
\end{defi}

Il y a une rétraction 
\begin{eqnarray*}
r : |\mathcal{Y}^{an} | & \twoheadrightarrow & \Delta \\
(y_1,\dots, y_n) & \mapsto & (v(y_1),\dots, v(y_n))
\end{eqnarray*}
qui est $\O_D^\times$-invariante, c'est à dire ne dépend que de la
chaîne d'isogénies 
$H_1\drt \dots\drt H_n\drt H_1$ et pas de la déformation $\rho$.

Comme au début de la section \ref{secotopp} l'isomorphisme entre les deux tours devrait induire un diagramme 
$$
\xymatrix{
|\mathcal{Y}^{an} | \ar[r] \ar@{->>}[d] & I\bc |\Omega | \ar@{->>}[d] \\
S (\mathcal{Y}^{an}) \ar@{-->}[r]^{?} & I\bc S (\Omega) 
}
$$
où $ I\bc S (\Omega) $ s'identifie à $|\mathcal{A}|$ et donc un morphisme
$$
\xymatrix{
\Delta \ar@{-->}[r] & |\mathcal{A}|
}
$$

\subsection{Action du groupe de Weyl affine sur $\Delta$}

Soit $W_{\text{aff}} = \Z^n/\Z \rtimes \mathfrak{S}_n$ où $\Z\hookrightarrow \Z^n$ est le plongement diagonal. 

De la même façon que l'on a défini des sections des correspondances de Hecke au niveau de l'espace $Newt$ on peut définir de telles sections sur $\Delta$.

\begin{prop}
Soient $(a_1,\dots,a_n)\in \N^n$ et $\s\in \mathfrak{S}_n$. Soit $K$ algébriquement clos, 
$$
H_1\xrig{\; \ph_1\;} \dots \ldrt H_n \xrig{\;\ph_n\;} H_1,\;\;\;\ph_n \circ \dots \circ \ph_1 =\pi
$$
une chaîne d'isogénies 
sur $\spf (\O_K)$ et $y_1,\dots,y_n\in\O_K$ tels que $\prod_{i} y_i =\pi$ un uplet associé. Soit $(e_1,\dots,e_n)$ une base de $T_p (H_1)$ telle que 
\begin{itemize}
\item $\exists \tau\in \mathfrak{S}_n\;\; (e_{\tau (1)},\dots,e_{\tau (n)})$ soit une base adaptée de $T_p (H_1)$
\item la filtration de $H_1 [\pi] (\O_K)$
 $$(<\pi^{-1} \bar{e}_1,\dots,\pi^{-1} \bar{e}_i, \bar{e}_{i+1},\dots,\bar{e}_n>)_{1\leq i\leq n}$$ soit égale à la filtration $$\ker \ph_1\subset \ker (\ph_2\circ \ph_1) \subset \dots \subset \ker (\ph_{n-1}\circ \dots \ph_1) \subset H_1 [\pi](\O_K)$$ 
\end{itemize}
Soit $C$ l'adhérence schématique du sous-groupe fini
$$
<\pi^{-a_1}\bar{e}_{\s (1)},\dots,\pi^{-{a_n}}\bar{e}_{\s(n)}>\subset T_p (H)\otimes F/\O_F
$$
Soit $H'_1=H_1/C$, $(\e_1,\dots,\e_n)$ la base de $T_p (H'_1)$ image de $(\pi^{-a_1} e_{\s(1)},\dots,\pi^{-a_n} e_{\s (n)})$ via $T_p (H_1)\hookrightarrow T_p (H'_1)$. Soit la chaîne d'isogénies 
$$
H'_1\xrig{\; \ph'_1\;} \dots \ldrt H'_n \xrig{\;\ph'_n\;} H'_1
$$
telle que $\forall i\;\; \ker ( \ph'_i\circ \dots\circ \ph'_1)$ soit l'adhérence schématique du groupe $$
<\pi^{-1}\bar{\e}_1,\dots,\pi^{-1}\bar{\e}_i,\bar{\e}_{i+1},\dots, \bar{\e}_n>
\subset T_p ({H'}_1)\otimes F/\O_F
$$
Soient ${y'}_1,\dots,{y'}_n\in \O_K$, $\prod_i y'_i=\pi$, des nombres associés à $H'_1\drt \dots \drt H'_n\drt H'_1$. 

Alors, $(v(y'_1),\dots,v(y'_n))\in \Delta$ ne dépend que de $(a_1,\dots,a_n)\rtimes \s\in W_{\text{aff}}$ et de $(v(y_1),\dots,v(y_n))\in \Delta$.
Cela définit une action de $W_{\text{aff}}$ sur $\Delta$. L'opérateur associé à $(a_1,\dots,a_n)\rtimes \s$ est noté $$(\pi^{-a_1},\dots,\pi^{-a_n})\rtimes \s$$
\end{prop}

\subsection{L'application $\Delta \ldrt Newt$}

\begin{defi}
Soit $\l \in ]0,\frac{1}{q-1}[$. On note $\eta_\l:[0,+\infty [\ldrt
[0,+\infty [$ la fonction continue linéaire par morceaux définie par 
$$
\eta_\l (x) = \left \{ \begin{array}{c}
 q x\text{ si } x\in [0,\l ] \\
(x-\l)+q\l \text{ si } x\in [\l,+\infty [
\end{array} \right.
$$
\end{defi}

\begin{lemm}
Pour $(v_1,\dots,v_n)\in \Delta$ soit 
$$
\eta (v_1,\dots,v_n) = \eta_{\frac{v_1}{q-1}}\circ \dots \circ
\eta_{\frac{v_n}{q-1}} :[0,+\infty [ \ldrt [0,+\infty [
$$
Soit 
$$
\mathcal{P} (v_1,\dots,v_n) = \left (\eta (v_1,\dots,v_n)^*\right
)_{|[1, q^n]}
$$
où on pose $f^* (x) = \sup \{ -xt + f(t)\;|\; t\in [0,+\infty [
\}$. Alors $\mathcal{P} (v_1,\dots,v_n) \in Newt$ et cela définit une
application
$$
\mathcal{P}: \Delta \ldrt Newt
$$
\end{lemm}

\begin{lemm}
Soit $H_1\xrig{\; \ph_1\;} \dots \ldrt H_n \xrig{\;\ph_n \;} H_1$ une
chaîne d'isogénies sur $\O_K$ et $(v_1,\dots,v_n)\in \Delta$ le uplet
associé. Alors $\mathcal{P} (v_1,\dots,v_n)$ est le polygone de Newton
de $H_1$.
\end{lemm}

\begin{coro}
Pour tous $(a_1,\dots,a_n)\rtimes \s \in W_{\text{aff}}$ et 
$\underline{v}\in \Delta\;\; \exists \tau\in \mathfrak{S}_n$
$$
\mathcal{P} ((\pi^{-a_1},\dots,\pi^{-a_n})\rtimes
\s. \underline{v})
= (\pi^{-a_{\tau (1)}},\dots,\pi^{-a_{\tau (n)}}).\mathcal{P} (\underline{v})
$$
\end{coro}

\subsection{Quartiers}

\begin{defi}
On note 
$$
Q(\Delta) =\{ (v_1,\dots,v_n)\in \Delta \; |\; v_1\geq \frac{v_2}{q}
\geq \dots \geq \frac{v_i}{q^i}\geq \dots \geq \frac{v_n}{q^n}\}
$$
\end{defi}

\begin{prop}
L'application $\mathcal{P}_{|Q(\Delta)}$ induit une bijection entre $Q
(\Delta)$ et $Newt$. \`A $(v_1,\dots,v_n)\in Q (\Delta)$ elle associe
le polygone de pentes $\l_1\geq \dots \geq \l_n$ où $\forall i\;
\l_i=\frac{v_i}{q^i -q^{i-1}}$. 

Soit $H_1\xrig{\; \ph_1\;} \dots \ldrt H_n \xrig{\;\ph_n \;} H_1$ une
chaîne d'isogénies sur $\O_K$. Le point associé dans $\Delta$ est dans
$Q(\Delta)$ ssi la filtration 
$$
(0)\subsetneq \ker \ph_1 \subsetneq \dots \subsetneq \ker (\ph_n \circ
\dots \circ \ph_1) =H[\pi]
$$
raffine la filtration de ramification inférieure (donnée par la
valuation des points de $\pi$-torsion) sur $H[\pi]$. 
\end{prop}

Rappelons qu'étant donnés deux drapeaux complets sur un $\Fq$-espaces vectoriel de dimension $n$ on peut définir leur invariant dans $\mathfrak{S}_n$ mesurant la position relative des sous-groupes de Borel du groupe linéaire associé.

\begin{prop}
On a la décomposition
$$
\Delta = \bigcup_{\s\in\mathfrak{S}_n} \s.Q(\Delta)
$$
Une chaine d'isogénies $H_1\xrig{\; \ph_1\;} \dots \ldrt H_n \xrig{\;\ph_n \;} H_1$ sur $\O_K$ donne lieu à un élément de $\s.Q(\Delta)$ ssi l'invariant des drapeaux 
$$\ker \ph_1\subset \ker (\ph_2\circ \ph_1) \subset \dots \subset \ker (\ph_{n-1}\circ \dots \ph_1) \subset H_1 [\pi](\O_K)$$ 
et un drapeaux complet raffinant la filtration de ramification sur $H_1[\pi]$ est $\s$.
L'application $\mathcal{P}_{|\s.Q(\Delta)}$ est l'isomorphisme composé 
$$
\s.Q(\Delta) \xrig{\;\s^{-1}\;} Q(\Delta) \xrig{\;\mathcal{P}\;} Newt
$$
Les applications $\forall \tau\in \mathfrak{S}_n\;\; \tau : \s.Q(\Delta) \iso \tau\s .Q(\Delta)$ sont affines. 
\end{prop}

\begin{rema}
Pour $\s\in\mathfrak{S}_n$ l'application affine $\s:Q(\Delta)\iso \s.Q(\Delta)$ est donnée par les formules suivantes. Soit $(v_1,\dots,v_n)\in Q(\Delta)$ et $\s.(v_1,\dots,v_n)=(v'_1,\dots,v'_n)$. Alors
$$
v'_1=\frac{v_{\s(1)}}{q^{\s(1)-1}}
$$
et $\forall i>1$ si $\tau_1\in \mathfrak{S}_i,\; \tau_2\in \mathfrak{S}_{i+1}$ sont tels que 
$$
\begin{array}{c}
\s \tau_1 (1)<\dots <\s\tau_1 (i-1) \\
\s \tau_2 (1)<\dots <\s\tau_2 (i)
\end{array}
$$
alors
$$
v'_i = \sum_{k=1}^i \frac{v_{\s\tau_2 (k)}}{q^{\s\tau_2 (k)-1}} -  \sum_{k=1}^{i-1} \frac{v_{\s\tau_1 (k)}}{q^{\s\tau_1 (k)-1}}
$$
\end{rema}

\subsection{La structure simpliciale sur $\Delta$}

\begin{defi}
On munit $\Delta$ de l'unique structure simpliciale $W_{\text{aff}}$-invariante pour laquelle 
l'application $\mathcal{P}_{|Q(\Delta)}:Q(\Delta)\iso Newt$ soit un isomorphisme simplicial. 
\end{defi}

\begin{prop}
Les simplexes de la structure simpliciale précédente sont des simplexes standards pour la structure affine naturelle de $\Delta$. L'action de $W_{\text{aff}}$ est affine sur chaque simplexe. 
\end{prop}

\subsection{\'Enoncé du théorème principal}

\begin{theo}
L'isomorphisme entre les tours de Lubin-Tate et de Drinfeld induit un isomorphisme simpliciale $W_{\text{aff}}$-équivariant 
$$
\Delta \iso |\mathcal{A}|
$$
affine sur chaque simplexe
tel que l'application d'oubli du niveau $\mathcal{Y}^{an}\ldrt \X^{an}$ s'identifie au quotient par $\mathfrak{S}_n$ 
$$
\xymatrix{
|\mathcal{Y}^{an}| \ar@{->>}[r]\ar[d] & \Delta \ar@{->>}[d]^{\mathcal{P}} \ar[r]^\sim & |\mathcal{A}| \ar@{->>}[d]^{\text{pr}_Q} \\
|\X^{an}|\ar@{->>}[r] &  Newt  \ar[r]^{\sim} & Q
}
$$
Dans cette bijection $Q(\Delta)$ correspond au quartier $Q\subset |\mathcal{A}|$ et le simplexe fondamental dans $Q$ correspond à
$$
\{ (v_1,\dots,v_n)\in \Delta \;|\; v_1\geq \frac{v_2}{q}\geq \dots\geq \frac{v_n}{q^{n-1}} \geq \frac{v_1}{q^n} \}
$$
Soit $G$ le barycentre du simplexe fondamental dans $|\mathcal{A}|$. 
Pour $1\leq i\leq n$ la
correspondance de Hecke $(H_1\drt\dots \drt H_n \drt H_1)
\longmapsto H_i$ est donnée par une rotation autour de $G$ dans
l'appartement $\mathcal{A}$ composée avec l'application quotient
$|\mathcal{A}|\twoheadrightarrow Q$. 

Soit $Q'\subset Q$ le quartier privé des murs $\a_{i,i+1} =0$ pour
$1\leq i\leq n-1$. Sur l'image réciproque de $Q'$ dans l'espace de
Lubin-Tate sans niveau le groupe $p$-divisible universel possède un
drapeau complet de sous-groupes canoniques. Cela fournit une section
sur cet ouvert admissible du revêtement en niveau Iwahori. Cette
section est donnée par le diagramme 
$$
\xymatrix{
|\mathcal{A} | \ar@{->>}[d]^{\text{pr}_Q} \\
Q & \ar@{_(->}[l]  Q'\ar@{_(->}[lu]
}
$$
\end{theo}

\begin{figure}[htbp]
   \begin{center}
      \input{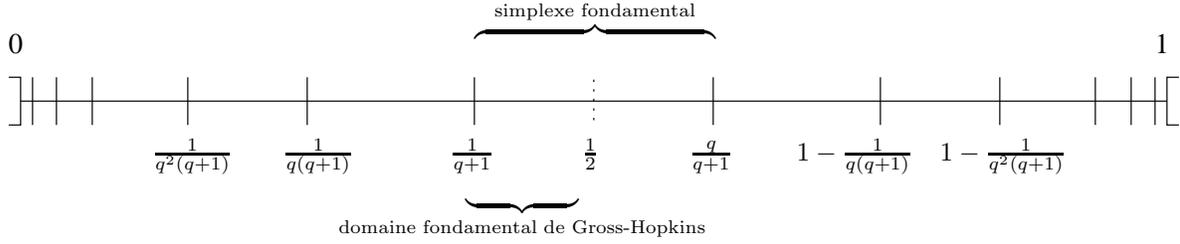}
   \end{center}
\caption{\footnotesize La structure simpliciale de $\Delta=]0 1[$ dans le
  cas de $\GL_2$
}
\end{figure}

\begin{figure}[htbp]
   \begin{center}
      \input{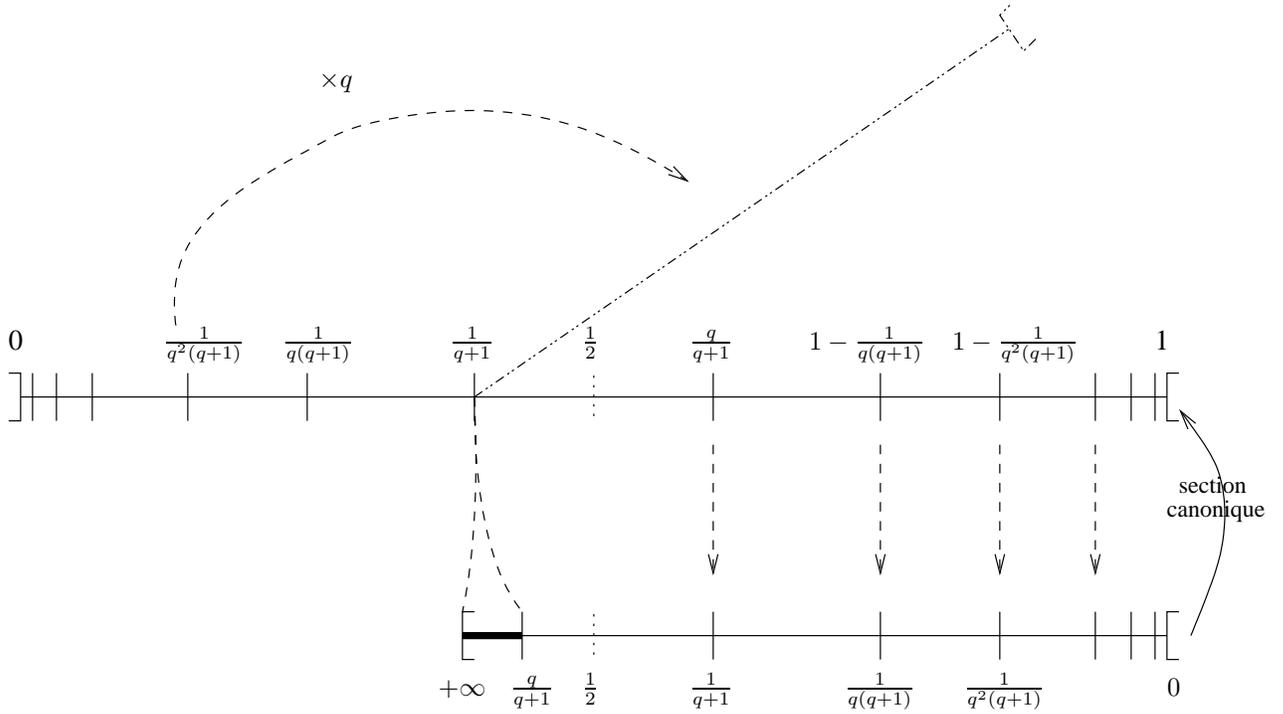}
   \end{center}
\caption{\footnotesize L'application entre les squelettes de l'espace de Lubin-Tate en
  niveau Iwahori et celui sans niveau pour $\GL_2$ (cette application est
  ``multivaluée'' en $\frac{1}{q+1}$). La section canonique est
  définie sur le segment $]0,\frac{q}{q+1}[$ d'image $]\frac{1}{q+1},
1[$. L'opérateur de Hecke donné par  l'élément non-trivial dans le
normalisateur de l'Iwahori modulo l'Iwahori qui agit comme automorphisme de l'espace en niveau
Iwahori  correspond à la symétrie de centre $\frac{1}{2}$, $x\mapsto 1-x$.}
\end{figure}

\newpage

\begin{figure}[htbp]
   \begin{center}
\begin{picture}(400,400)
\path(0,0)(420,0)
\path(0,0)(210,363)
\path(210,363)(420,0)
\path(238,81)(161,116)
\path(238,81)(286,116)
\path(161,116)(167,57)
\path(229,165)(176,189)
\path(161,116)(229,165)
\path(229,165)(238,81)
\path(242,40)(167,57)
\path(242,40)(291,57)
\path(124,133)(114,81)
\path(263,189)(224,224)
\path(124,133)(176,189)
\path(263,189)(286,116)
\path(167,57)(114,81)
\path(291,57)(325,81)
\path(114,81)(118,40)
\path(224,224)(186,241)
\path(176,189)(224,224)
\path(286,116)(291,57)
\path(244,20)(170,28)
\path(244,20)(294,28)
\path(105,141)(87,94)
\path(280,201)(247,241)
\path(105,141)(149,201)
\path(280,201)(309,133)
\path(170,28)(118,40)
\path(294,28)(329,40)
\path(87,94)(80,57)
\path(247,241)(219,265)
\path(149,201)(186,241)
\path(309,133)(325,81)
\path(118,40)(80,57)
\path(329,40)(353,57)
\path(80,57)(83,28)
\path(219,265)(193,277)
\path(186,241)(219,265)
\path(325,81)(329,40)
\path(245,10)(172,14)
\path(245,10)(296,14)
\path(95,146)(74,99)
\path(288,207)(259,249)
\path(95,146)(136,207)
\path(288,207)(321,141)
\path(172,14)(120,20)
\path(296,14)(331,20)
\path(74,99)(61,66)
\path(259,249)(236,277)
\path(136,207)(167,249)
\path(321,141)(342,94)
\path(120,20)(83,28)
\path(331,20)(356,28)
\path(61,66)(56,40)
\path(236,277)(216,294)
\path(167,249)(193,277)
\path(342,94)(353,57)
\path(83,28)(56,40)
\path(356,28)(373,40)
\path(56,40)(58,20)
\path(216,294)(198,302)
\path(193,277)(216,294)
\path(353,57)(356,28)
\path(245,4)(172,7)
\path(245,4)(297,7)
\path(91,148)(67,102)
\path(292,210)(265,253)
\path(91,148)(129,210)
\path(292,210)(327,146)
\path(172,7)(121,10)
\path(297,7)(332,10)
\path(67,102)(52,70)
\path(265,253)(244,283)
\path(129,210)(158,253)
\path(327,146)(350,99)
\path(121,10)(84,14)
\path(332,10)(357,14)
\path(52,70)(43,46)
\path(244,283)(228,302)
\path(158,253)(179,283)
\path(350,99)(365,66)
\path(84,14)(58,20)
\path(357,14)(375,20)
\path(43,46)(39,28)
\path(228,302)(214,314)
\path(179,283)(198,302)
\path(365,66)(373,40)
\path(58,20)(39,28)
\path(375,20)(387,28)
\path(39,28)(41,14)
\path(214,314)(201,320)
\path(198,302)(214,314)
\path(373,40)(375,20)
\path(246,2)(173,3)
\path(246,2)(297,3)
\path(89,149)(64,104)
\path(294,211)(268,255)
\path(89,149)(126,211)
\path(294,211)(330,148)
\path(173,3)(121,4)
\path(297,3)(333,4)
\path(64,104)(47,72)
\path(268,255)(248,286)
\path(126,211)(153,255)
\path(330,148)(354,102)
\path(121,4)(85,7)
\path(333,4)(358,7)
\path(47,72)(36,49)
\path(248,286)(234,307)
\path(153,255)(173,286)
\path(354,102)(371,70)
\path(85,7)(59,10)
\path(358,7)(376,10)
\path(36,49)(30,32)
\path(234,307)(223,320)
\path(173,286)(188,307)
\path(371,70)(381,46)
\path(59,10)(41,14)
\path(376,10)(388,14)
\path(30,32)(28,20)
\path(223,320)(213,329)
\path(188,307)(201,320)
\path(381,46)(387,28)
\path(41,14)(28,20)
\path(388,14)(396,20)
\path(28,20)(29,10)
\path(213,329)(204,333)
\path(201,320)(213,329)
\path(387,28)(388,14)
\path(246,1)(173,1)
\path(246,1)(297,1)
\path(87,149)(62,105)
\path(295,212)(269,256)
\path(87,149)(124,212)
\path(295,212)(331,149)
\path(173,1)(121,2)
\path(297,1)(333,2)
\path(62,105)(45,73)
\path(269,256)(251,287)
\path(124,212)(151,256)
\path(331,149)(356,104)
\path(121,2)(85,3)
\path(333,2)(359,3)
\path(45,73)(33,51)
\path(251,287)(237,309)
\path(151,256)(170,287)
\path(356,104)(374,72)
\path(85,3)(60,4)
\path(359,3)(376,4)
\path(33,51)(25,34)
\path(237,309)(227,323)
\path(170,287)(184,309)
\path(374,72)(385,49)
\path(60,4)(41,7)
\path(376,4)(389,7)
\path(25,34)(21,23)
\path(227,323)(219,333)
\path(184,309)(195,323)
\path(385,49)(392,32)
\path(41,7)(29,10)
\path(389,7)(397,10)
\path(21,23)(19,14)
\path(219,333)(212,339)
\path(195,323)(204,333)
\path(392,32)(396,20)
\path(29,10)(19,14)
\path(397,10)(403,14)
\path(19,14)(20,7)
\path(212,339)(205,342)
\path(204,333)(212,339)
\path(396,20)(397,10)
\path(246,0)(173,0)
\path(246,0)(297,0)
\path(87,150)(61,105)
\path(296,213)(270,257)
\path(87,150)(124,213)
\path(296,213)(332,149)
\path(173,0)(122,1)
\path(297,0)(333,1)
\path(61,105)(44,74)
\path(270,257)(252,288)
\path(124,213)(150,257)
\path(332,149)(357,105)
\path(122,1)(85,1)
\path(333,1)(359,1)
\path(44,74)(31,51)
\path(252,288)(238,310)
\path(150,257)(168,288)
\path(357,105)(375,73)
\path(85,1)(60,2)
\path(359,1)(377,2)
\path(31,51)(23,35)
\path(238,310)(229,325)
\path(168,288)(181,310)
\path(375,73)(387,51)
\path(60,2)(42,3)
\path(377,2)(389,3)
\path(23,35)(18,24)
\path(229,325)(222,335)
\path(181,310)(191,325)
\path(387,51)(395,34)
\path(42,3)(29,4)
\path(389,3)(398,4)
\path(18,24)(15,16)
\path(222,335)(216,342)
\path(191,325)(199,335)
\path(395,34)(400,23)
\path(29,4)(20,7)
\path(398,4)(404,7)
\path(15,16)(13,10)
\path(216,342)(211,346)
\path(199,335)(205,342)
\path(400,23)(403,14)
\path(20,7)(13,10)
\path(404,7)(408,10)
\path(13,10)(14,4)
\path(211,346)(207,348)
\path(205,342)(211,346)
\path(403,14)(404,7)
\path(246,0)(173,0)
\path(246,0)(297,0)
\path(87,150)(61,105)
\path(296,213)(270,257)
\path(87,150)(123,213)
\path(296,213)(332,150)
\path(173,0)(122,0)
\path(297,0)(333,0)
\path(61,105)(43,74)
\path(270,257)(252,288)
\path(123,213)(149,257)
\path(332,150)(358,105)
\path(122,0)(85,0)
\path(333,0)(359,0)
\path(43,74)(31,52)
\path(252,288)(239,310)
\path(149,257)(167,288)
\path(358,105)(376,74)
\path(85,0)(60,1)
\path(359,0)(377,1)
\path(31,52)(22,36)
\path(239,310)(230,326)
\path(167,288)(180,310)
\path(376,74)(388,51)
\path(60,1)(42,1)
\path(377,1)(389,1)
\path(22,36)(16,25)
\path(230,326)(223,336)
\path(180,310)(190,326)
\path(388,51)(397,35)
\path(42,1)(29,2)
\path(389,1)(398,2)
\path(16,25)(12,17)
\path(223,336)(218,343)
\path(190,326)(197,336)
\path(397,35)(402,24)
\path(29,2)(20,3)
\path(398,2)(404,3)
\path(12,17)(10,11)
\path(218,343)(214,348)
\path(197,336)(202,343)
\path(402,24)(406,16)
\path(20,3)(14,4)
\path(404,3)(409,4)
\path(10,11)(9,7)
\path(214,348)(211,351)
\path(202,343)(207,348)
\path(406,16)(408,10)
\path(14,4)(9,7)
\path(409,4)(411,7)
\path(9,7)(10,3)
\path(211,351)(207,353)
\path(207,348)(211,351)
\path(408,10)(409,4)
\path(246,0)(173,0)
\path(246,0)(297,0)
\path(86,150)(61,105)
\path(296,213)(270,257)
\path(86,150)(123,213)
\path(296,213)(333,150)
\path(173,0)(122,0)
\path(297,0)(333,0)
\path(61,105)(43,74)
\path(270,257)(252,289)
\path(123,213)(149,257)
\path(333,150)(358,105)
\path(122,0)(86,0)
\path(333,0)(359,0)
\path(43,74)(30,52)
\path(252,289)(239,310)
\path(149,257)(167,289)
\path(358,105)(376,74)
\path(86,0)(60,0)
\path(359,0)(377,0)
\path(30,52)(21,36)
\path(239,310)(230,326)
\path(167,289)(180,310)
\path(376,74)(389,52)
\path(60,0)(42,0)
\path(377,0)(389,0)
\path(21,36)(15,25)
\path(230,326)(224,337)
\path(180,310)(189,326)
\path(389,52)(397,36)
\path(42,0)(29,1)
\path(389,0)(398,1)
\path(15,25)(11,17)
\path(224,337)(219,344)
\path(189,326)(196,337)
\path(397,36)(403,25)
\path(29,1)(20,1)
\path(398,1)(404,1)
\path(11,17)(9,12)
\path(219,344)(216,349)
\path(196,337)(201,344)
\path(403,25)(408,17)
\path(20,1)(14,2)
\path(404,1)(409,2)
\path(9,12)(7,8)
\path(216,349)(213,353)
\path(201,344)(204,349)
\path(408,17)(410,11)
\path(14,2)(10,3)
\path(409,2)(412,3)
\path(7,8)(6,4)
\path(213,353)(210,355)
\path(204,349)(207,353)
\path(410,11)(411,7)
\path(10,3)(6,4)
\path(412,3)(414,4)
\path(6,4)(7,2)
\path(210,355)(208,356)
\path(207,353)(210,355)
\path(411,7)(412,3)
\path(246,0)(173,0)
\path(246,0)(297,0)
\path(86,150)(61,105)
\path(296,213)(271,257)
\path(86,150)(123,213)
\path(296,213)(333,150)
\path(173,0)(122,0)
\path(297,0)(333,0)
\path(61,105)(43,74)
\path(271,257)(252,289)
\path(123,213)(149,257)
\path(333,150)(358,105)
\path(122,0)(86,0)
\path(333,0)(359,0)
\path(43,74)(30,52)
\path(252,289)(240,311)
\path(149,257)(167,289)
\path(358,105)(376,74)
\path(86,0)(60,0)
\path(359,0)(377,0)
\path(30,52)(21,36)
\path(240,311)(231,326)
\path(167,289)(179,311)
\path(376,74)(389,52)
\path(60,0)(42,0)
\path(377,0)(389,0)
\path(21,36)(15,25)
\path(231,326)(224,337)
\path(179,311)(189,326)
\path(389,52)(398,36)
\path(42,0)(29,0)
\path(389,0)(398,0)
\path(15,25)(11,18)
\path(224,337)(220,345)
\path(189,326)(195,337)
\path(398,36)(404,25)
\path(29,0)(21,0)
\path(398,0)(405,0)
\path(11,18)(8,12)
\path(220,345)(216,350)
\path(195,337)(200,345)
\path(404,25)(408,17)
\path(21,0)(14,1)
\path(405,0)(409,1)
\path(8,12)(6,8)
\path(216,350)(214,353)
\path(200,345)(203,350)
\path(408,17)(411,12)
\path(14,1)(10,1)
\path(409,1)(412,1)
\path(6,8)(5,5)
\path(214,353)(212,356)
\path(203,350)(206,353)
\path(411,12)(413,8)
\path(10,1)(7,2)
\path(412,1)(414,2)
\path(5,5)(4,3)
\path(212,356)(210,357)
\path(206,353)(208,356)
\path(413,8)(414,4)
\path(7,2)(4,3)
\path(414,2)(415,3)
\path(4,3)(5,1)
\path(210,357)(208,358)
\path(208,356)(210,357)
\path(414,4)(414,2)
\path(229,165)(238,81)
\path(229,165)(238,81)
\path(238,81)(161,116)
\path(161,116)(229,165)
\path(238,81)(161,116)
\path(161,116)(229,165)
\path(238,81)(242,40)
\path(238,81)(242,40)
\path(161,116)(124,133)
\path(229,165)(263,189)
\path(161,116)(124,133)
\path(229,165)(263,189)
\path(161,116)(167,57)
\path(286,116)(291,57)
\path(167,57)(114,81)
\path(176,189)(224,224)
\path(229,165)(176,189)
\path(238,81)(286,116)
\path(242,40)(244,20)
\path(242,40)(244,20)
\path(124,133)(105,141)
\path(263,189)(280,201)
\path(124,133)(105,141)
\path(263,189)(280,201)
\path(167,57)(170,28)
\path(291,57)(294,28)
\path(114,81)(87,94)
\path(224,224)(247,241)
\path(176,189)(149,201)
\path(286,116)(309,133)
\path(114,81)(118,40)
\path(325,81)(329,40)
\path(118,40)(80,57)
\path(186,241)(219,265)
\path(224,224)(186,241)
\path(291,57)(325,81)
\path(244,20)(245,10)
\path(244,20)(245,10)
\path(105,141)(95,146)
\path(280,201)(288,207)
\path(105,141)(95,146)
\path(280,201)(288,207)
\path(170,28)(172,14)
\path(294,28)(296,14)
\path(87,94)(74,99)
\path(247,241)(259,249)
\path(149,201)(136,207)
\path(309,133)(321,141)
\path(118,40)(120,20)
\path(329,40)(331,20)
\path(80,57)(61,66)
\path(219,265)(236,277)
\path(186,241)(167,249)
\path(325,81)(342,94)
\path(80,57)(83,28)
\path(353,57)(356,28)
\path(83,28)(56,40)
\path(193,277)(216,294)
\path(219,265)(193,277)
\path(329,40)(353,57)
\path(245,10)(245,4)
\path(245,10)(245,4)
\path(95,146)(91,148)
\path(288,207)(292,210)
\path(95,146)(91,148)
\path(288,207)(292,210)
\path(172,14)(172,7)
\path(296,14)(297,7)
\path(74,99)(67,102)
\path(259,249)(265,253)
\path(136,207)(129,210)
\path(321,141)(327,146)
\path(120,20)(121,10)
\path(331,20)(332,10)
\path(61,66)(52,70)
\path(236,277)(244,283)
\path(167,249)(158,253)
\path(342,94)(350,99)
\path(83,28)(84,14)
\path(356,28)(357,14)
\path(56,40)(43,46)
\path(216,294)(228,302)
\path(193,277)(179,283)
\path(353,57)(365,66)
\path(56,40)(58,20)
\path(373,40)(375,20)
\path(58,20)(39,28)
\path(198,302)(214,314)
\path(216,294)(198,302)
\path(356,28)(373,40)
\path(245,4)(246,2)
\path(245,4)(246,2)
\path(91,148)(89,149)
\path(292,210)(294,211)
\path(91,148)(89,149)
\path(292,210)(294,211)
\path(172,7)(173,3)
\path(297,7)(297,3)
\path(67,102)(64,104)
\path(265,253)(268,255)
\path(129,210)(126,211)
\path(327,146)(330,148)
\path(121,10)(121,4)
\path(332,10)(333,4)
\path(52,70)(47,72)
\path(244,283)(248,286)
\path(158,253)(153,255)
\path(350,99)(354,102)
\path(84,14)(85,7)
\path(357,14)(358,7)
\path(43,46)(36,49)
\path(228,302)(234,307)
\path(179,283)(173,286)
\path(365,66)(371,70)
\path(58,20)(59,10)
\path(375,20)(376,10)
\path(39,28)(30,32)
\path(214,314)(223,320)
\path(198,302)(188,307)
\path(373,40)(381,46)
\path(39,28)(41,14)
\path(387,28)(388,14)
\path(41,14)(28,20)
\path(201,320)(213,329)
\path(214,314)(201,320)
\path(375,20)(387,28)
\path(246,2)(246,1)
\path(246,2)(246,1)
\path(89,149)(87,149)
\path(294,211)(295,212)
\path(89,149)(87,149)
\path(294,211)(295,212)
\path(173,3)(173,1)
\path(297,3)(297,1)
\path(64,104)(62,105)
\path(268,255)(269,256)
\path(126,211)(124,212)
\path(330,148)(331,149)
\path(121,4)(121,2)
\path(333,4)(333,2)
\path(47,72)(45,73)
\path(248,286)(251,287)
\path(153,255)(151,256)
\path(354,102)(356,104)
\path(85,7)(85,3)
\path(358,7)(359,3)
\path(36,49)(33,51)
\path(234,307)(237,309)
\path(173,286)(170,287)
\path(371,70)(374,72)
\path(59,10)(60,4)
\path(376,10)(376,4)
\path(30,32)(25,34)
\path(223,320)(227,323)
\path(188,307)(184,309)
\path(381,46)(385,49)
\path(41,14)(41,7)
\path(388,14)(389,7)
\path(28,20)(21,23)
\path(213,329)(219,333)
\path(201,320)(195,323)
\path(387,28)(392,32)
\path(28,20)(29,10)
\path(396,20)(397,10)
\path(29,10)(19,14)
\path(204,333)(212,339)
\path(213,329)(204,333)
\path(388,14)(396,20)
\path(246,1)(246,0)
\path(246,1)(246,0)
\path(87,149)(87,150)
\path(295,212)(296,213)
\path(87,149)(87,150)
\path(295,212)(296,213)
\path(173,1)(173,0)
\path(297,1)(297,0)
\path(62,105)(61,105)
\path(269,256)(270,257)
\path(124,212)(124,213)
\path(331,149)(332,149)
\path(121,2)(122,1)
\path(333,2)(333,1)
\path(45,73)(44,74)
\path(251,287)(252,288)
\path(151,256)(150,257)
\path(356,104)(357,105)
\path(85,3)(85,1)
\path(359,3)(359,1)
\path(33,51)(31,51)
\path(237,309)(238,310)
\path(170,287)(168,288)
\path(374,72)(375,73)
\path(60,4)(60,2)
\path(376,4)(377,2)
\path(25,34)(23,35)
\path(227,323)(229,325)
\path(184,309)(181,310)
\path(385,49)(387,51)
\path(41,7)(42,3)
\path(389,7)(389,3)
\path(21,23)(18,24)
\path(219,333)(222,335)
\path(195,323)(191,325)
\path(392,32)(395,34)
\path(29,10)(29,4)
\path(397,10)(398,4)
\path(19,14)(15,16)
\path(212,339)(216,342)
\path(204,333)(199,335)
\path(396,20)(400,23)
\path(19,14)(20,7)
\path(403,14)(404,7)
\path(20,7)(13,10)
\path(205,342)(211,346)
\path(212,339)(205,342)
\path(397,10)(403,14)
\path(246,0)(246,0)
\path(246,0)(246,0)
\path(87,150)(87,150)
\path(296,213)(296,213)
\path(87,150)(87,150)
\path(296,213)(296,213)
\path(173,0)(173,0)
\path(297,0)(297,0)
\path(61,105)(61,105)
\path(270,257)(270,257)
\path(124,213)(123,213)
\path(332,149)(332,150)
\path(122,1)(122,0)
\path(333,1)(333,0)
\path(44,74)(43,74)
\path(252,288)(252,288)
\path(150,257)(149,257)
\path(357,105)(358,105)
\path(85,1)(85,0)
\path(359,1)(359,0)
\path(31,51)(31,52)
\path(238,310)(239,310)
\path(168,288)(167,288)
\path(375,73)(376,74)
\path(60,2)(60,1)
\path(377,2)(377,1)
\path(23,35)(22,36)
\path(229,325)(230,326)
\path(181,310)(180,310)
\path(387,51)(388,51)
\path(42,3)(42,1)
\path(389,3)(389,1)
\path(18,24)(16,25)
\path(222,335)(223,336)
\path(191,325)(190,326)
\path(395,34)(397,35)
\path(29,4)(29,2)
\path(398,4)(398,2)
\path(15,16)(12,17)
\path(216,342)(218,343)
\path(199,335)(197,336)
\path(400,23)(402,24)
\path(20,7)(20,3)
\path(404,7)(404,3)
\path(13,10)(10,11)
\path(211,346)(214,348)
\path(205,342)(202,343)
\path(403,14)(406,16)
\path(13,10)(14,4)
\path(408,10)(409,4)
\path(14,4)(9,7)
\path(207,348)(211,351)
\path(211,346)(207,348)
\path(404,7)(408,10)
\path(246,0)(246,0)
\path(246,0)(246,0)
\path(87,150)(86,150)
\path(296,213)(296,213)
\path(87,150)(86,150)
\path(296,213)(296,213)
\path(173,0)(173,0)
\path(297,0)(297,0)
\path(61,105)(61,105)
\path(270,257)(270,257)
\path(123,213)(123,213)
\path(332,150)(333,150)
\path(122,0)(122,0)
\path(333,0)(333,0)
\path(43,74)(43,74)
\path(252,288)(252,289)
\path(149,257)(149,257)
\path(358,105)(358,105)
\path(85,0)(86,0)
\path(359,0)(359,0)
\path(31,52)(30,52)
\path(239,310)(239,310)
\path(167,288)(167,289)
\path(376,74)(376,74)
\path(60,1)(60,0)
\path(377,1)(377,0)
\path(22,36)(21,36)
\path(230,326)(230,326)
\path(180,310)(180,310)
\path(388,51)(389,52)
\path(42,1)(42,0)
\path(389,1)(389,0)
\path(16,25)(15,25)
\path(223,336)(224,337)
\path(190,326)(189,326)
\path(397,35)(397,36)
\path(29,2)(29,1)
\path(398,2)(398,1)
\path(12,17)(11,17)
\path(218,343)(219,344)
\path(197,336)(196,337)
\path(402,24)(403,25)
\path(20,3)(20,1)
\path(404,3)(404,1)
\path(10,11)(9,12)
\path(214,348)(216,349)
\path(202,343)(201,344)
\path(406,16)(408,17)
\path(14,4)(14,2)
\path(409,4)(409,2)
\path(9,7)(7,8)
\path(211,351)(213,353)
\path(207,348)(204,349)
\path(408,10)(410,11)
\path(9,7)(10,3)
\path(411,7)(412,3)
\path(10,3)(6,4)
\path(207,353)(210,355)
\path(211,351)(207,353)
\path(409,4)(411,7)
\path(246,0)(246,0)
\path(246,0)(246,0)
\path(86,150)(86,150)
\path(296,213)(296,213)
\path(86,150)(86,150)
\path(296,213)(296,213)
\path(173,0)(173,0)
\path(297,0)(297,0)
\path(61,105)(61,105)
\path(270,257)(271,257)
\path(123,213)(123,213)
\path(333,150)(333,150)
\path(122,0)(122,0)
\path(333,0)(333,0)
\path(43,74)(43,74)
\path(252,289)(252,289)
\path(149,257)(149,257)
\path(358,105)(358,105)
\path(86,0)(86,0)
\path(359,0)(359,0)
\path(30,52)(30,52)
\path(239,310)(240,311)
\path(167,289)(167,289)
\path(376,74)(376,74)
\path(60,0)(60,0)
\path(377,0)(377,0)
\path(21,36)(21,36)
\path(230,326)(231,326)
\path(180,310)(179,311)
\path(389,52)(389,52)
\path(42,0)(42,0)
\path(389,0)(389,0)
\path(15,25)(15,25)
\path(224,337)(224,337)
\path(189,326)(189,326)
\path(397,36)(398,36)
\path(29,1)(29,0)
\path(398,1)(398,0)
\path(11,17)(11,18)
\path(219,344)(220,345)
\path(196,337)(195,337)
\path(403,25)(404,25)
\path(20,1)(21,0)
\path(404,1)(405,0)
\path(9,12)(8,12)
\path(216,349)(216,350)
\path(201,344)(200,345)
\path(408,17)(408,17)
\path(14,2)(14,1)
\path(409,2)(409,1)
\path(7,8)(6,8)
\path(213,353)(214,353)
\path(204,349)(203,350)
\path(410,11)(411,12)
\path(10,3)(10,1)
\path(412,3)(412,1)
\path(6,4)(5,5)
\path(210,355)(212,356)
\path(207,353)(206,353)
\path(411,7)(413,8)
\path(6,4)(7,2)
\path(414,4)(414,2)
\path(7,2)(4,3)
\path(208,356)(210,357)
\path(210,355)(208,356)
\path(412,3)(414,4)
\path(229,165)(161,116)
\path(229,165)(286,116)
\path(238,81)(167,57)
\path(161,116)(176,189)
\path(238,81)(229,165)
\path(161,116)(238,81)
\path(238,81)(167,57)
\path(238,81)(291,57)
\path(161,116)(114,81)
\path(229,165)(224,224)
\path(161,116)(176,189)
\path(229,165)(286,116)
\path(161,116)(114,81)
\path(286,116)(325,81)
\path(167,57)(118,40)
\path(176,189)(186,241)
\path(229,165)(224,224)
\path(238,81)(291,57)
\path(242,40)(170,28)
\path(242,40)(294,28)
\path(124,133)(87,94)
\path(263,189)(247,241)
\path(124,133)(149,201)
\path(263,189)(309,133)
\path(167,57)(118,40)
\path(291,57)(329,40)
\path(114,81)(80,57)
\path(224,224)(219,265)
\path(176,189)(186,241)
\path(286,116)(325,81)
\path(114,81)(80,57)
\path(325,81)(353,57)
\path(118,40)(83,28)
\path(186,241)(193,277)
\path(224,224)(219,265)
\path(291,57)(329,40)
\path(244,20)(172,14)
\path(244,20)(296,14)
\path(105,141)(74,99)
\path(280,201)(259,249)
\path(105,141)(136,207)
\path(280,201)(321,141)
\path(170,28)(120,20)
\path(294,28)(331,20)
\path(87,94)(61,66)
\path(247,241)(236,277)
\path(149,201)(167,249)
\path(309,133)(342,94)
\path(118,40)(83,28)
\path(329,40)(356,28)
\path(80,57)(56,40)
\path(219,265)(216,294)
\path(186,241)(193,277)
\path(325,81)(353,57)
\path(80,57)(56,40)
\path(353,57)(373,40)
\path(83,28)(58,20)
\path(193,277)(198,302)
\path(219,265)(216,294)
\path(329,40)(356,28)
\path(245,10)(172,7)
\path(245,10)(297,7)
\path(95,146)(67,102)
\path(288,207)(265,253)
\path(95,146)(129,210)
\path(288,207)(327,146)
\path(172,14)(121,10)
\path(296,14)(332,10)
\path(74,99)(52,70)
\path(259,249)(244,283)
\path(136,207)(158,253)
\path(321,141)(350,99)
\path(120,20)(84,14)
\path(331,20)(357,14)
\path(61,66)(43,46)
\path(236,277)(228,302)
\path(167,249)(179,283)
\path(342,94)(365,66)
\path(83,28)(58,20)
\path(356,28)(375,20)
\path(56,40)(39,28)
\path(216,294)(214,314)
\path(193,277)(198,302)
\path(353,57)(373,40)
\path(56,40)(39,28)
\path(373,40)(387,28)
\path(58,20)(41,14)
\path(198,302)(201,320)
\path(216,294)(214,314)
\path(356,28)(375,20)
\path(245,4)(173,3)
\path(245,4)(297,3)
\path(91,148)(64,104)
\path(292,210)(268,255)
\path(91,148)(126,211)
\path(292,210)(330,148)
\path(172,7)(121,4)
\path(297,7)(333,4)
\path(67,102)(47,72)
\path(265,253)(248,286)
\path(129,210)(153,255)
\path(327,146)(354,102)
\path(121,10)(85,7)
\path(332,10)(358,7)
\path(52,70)(36,49)
\path(244,283)(234,307)
\path(158,253)(173,286)
\path(350,99)(371,70)
\path(84,14)(59,10)
\path(357,14)(376,10)
\path(43,46)(30,32)
\path(228,302)(223,320)
\path(179,283)(188,307)
\path(365,66)(381,46)
\path(58,20)(41,14)
\path(375,20)(388,14)
\path(39,28)(28,20)
\path(214,314)(213,329)
\path(198,302)(201,320)
\path(373,40)(387,28)
\path(39,28)(28,20)
\path(387,28)(396,20)
\path(41,14)(29,10)
\path(201,320)(204,333)
\path(214,314)(213,329)
\path(375,20)(388,14)
\path(246,2)(173,1)
\path(246,2)(297,1)
\path(89,149)(62,105)
\path(294,211)(269,256)
\path(89,149)(124,212)
\path(294,211)(331,149)
\path(173,3)(121,2)
\path(297,3)(333,2)
\path(64,104)(45,73)
\path(268,255)(251,287)
\path(126,211)(151,256)
\path(330,148)(356,104)
\path(121,4)(85,3)
\path(333,4)(359,3)
\path(47,72)(33,51)
\path(248,286)(237,309)
\path(153,255)(170,287)
\path(354,102)(374,72)
\path(85,7)(60,4)
\path(358,7)(376,4)
\path(36,49)(25,34)
\path(234,307)(227,323)
\path(173,286)(184,309)
\path(371,70)(385,49)
\path(59,10)(41,7)
\path(376,10)(389,7)
\path(30,32)(21,23)
\path(223,320)(219,333)
\path(188,307)(195,323)
\path(381,46)(392,32)
\path(41,14)(29,10)
\path(388,14)(397,10)
\path(28,20)(19,14)
\path(213,329)(212,339)
\path(201,320)(204,333)
\path(387,28)(396,20)
\path(28,20)(19,14)
\path(396,20)(403,14)
\path(29,10)(20,7)
\path(204,333)(205,342)
\path(213,329)(212,339)
\path(388,14)(397,10)
\path(246,1)(173,0)
\path(246,1)(297,0)
\path(87,149)(61,105)
\path(295,212)(270,257)
\path(87,149)(124,213)
\path(295,212)(332,149)
\path(173,1)(122,1)
\path(297,1)(333,1)
\path(62,105)(44,74)
\path(269,256)(252,288)
\path(124,212)(150,257)
\path(331,149)(357,105)
\path(121,2)(85,1)
\path(333,2)(359,1)
\path(45,73)(31,51)
\path(251,287)(238,310)
\path(151,256)(168,288)
\path(356,104)(375,73)
\path(85,3)(60,2)
\path(359,3)(377,2)
\path(33,51)(23,35)
\path(237,309)(229,325)
\path(170,287)(181,310)
\path(374,72)(387,51)
\path(60,4)(42,3)
\path(376,4)(389,3)
\path(25,34)(18,24)
\path(227,323)(222,335)
\path(184,309)(191,325)
\path(385,49)(395,34)
\path(41,7)(29,4)
\path(389,7)(398,4)
\path(21,23)(15,16)
\path(219,333)(216,342)
\path(195,323)(199,335)
\path(392,32)(400,23)
\path(29,10)(20,7)
\path(397,10)(404,7)
\path(19,14)(13,10)
\path(212,339)(211,346)
\path(204,333)(205,342)
\path(396,20)(403,14)
\path(19,14)(13,10)
\path(403,14)(408,10)
\path(20,7)(14,4)
\path(205,342)(207,348)
\path(212,339)(211,346)
\path(397,10)(404,7)
\path(246,0)(173,0)
\path(246,0)(297,0)
\path(87,150)(61,105)
\path(296,213)(270,257)
\path(87,150)(123,213)
\path(296,213)(332,150)
\path(173,0)(122,0)
\path(297,0)(333,0)
\path(61,105)(43,74)
\path(270,257)(252,288)
\path(124,213)(149,257)
\path(332,149)(358,105)
\path(122,1)(85,0)
\path(333,1)(359,0)
\path(44,74)(31,52)
\path(252,288)(239,310)
\path(150,257)(167,288)
\path(357,105)(376,74)
\path(85,1)(60,1)
\path(359,1)(377,1)
\path(31,51)(22,36)
\path(238,310)(230,326)
\path(168,288)(180,310)
\path(375,73)(388,51)
\path(60,2)(42,1)
\path(377,2)(389,1)
\path(23,35)(16,25)
\path(229,325)(223,336)
\path(181,310)(190,326)
\path(387,51)(397,35)
\path(42,3)(29,2)
\path(389,3)(398,2)
\path(18,24)(12,17)
\path(222,335)(218,343)
\path(191,325)(197,336)
\path(395,34)(402,24)
\path(29,4)(20,3)
\path(398,4)(404,3)
\path(15,16)(10,11)
\path(216,342)(214,348)
\path(199,335)(202,343)
\path(400,23)(406,16)
\path(20,7)(14,4)
\path(404,7)(409,4)
\path(13,10)(9,7)
\path(211,346)(211,351)
\path(205,342)(207,348)
\path(403,14)(408,10)
\path(13,10)(9,7)
\path(408,10)(411,7)
\path(14,4)(10,3)
\path(207,348)(207,353)
\path(211,346)(211,351)
\path(404,7)(409,4)
\path(246,0)(173,0)
\path(246,0)(297,0)
\path(87,150)(61,105)
\path(296,213)(270,257)
\path(87,150)(123,213)
\path(296,213)(333,150)
\path(173,0)(122,0)
\path(297,0)(333,0)
\path(61,105)(43,74)
\path(270,257)(252,289)
\path(123,213)(149,257)
\path(332,150)(358,105)
\path(122,0)(86,0)
\path(333,0)(359,0)
\path(43,74)(30,52)
\path(252,288)(239,310)
\path(149,257)(167,289)
\path(358,105)(376,74)
\path(85,0)(60,0)
\path(359,0)(377,0)
\path(31,52)(21,36)
\path(239,310)(230,326)
\path(167,288)(180,310)
\path(376,74)(389,52)
\path(60,1)(42,0)
\path(377,1)(389,0)
\path(22,36)(15,25)
\path(230,326)(224,337)
\path(180,310)(189,326)
\path(388,51)(397,36)
\path(42,1)(29,1)
\path(389,1)(398,1)
\path(16,25)(11,17)
\path(223,336)(219,344)
\path(190,326)(196,337)
\path(397,35)(403,25)
\path(29,2)(20,1)
\path(398,2)(404,1)
\path(12,17)(9,12)
\path(218,343)(216,349)
\path(197,336)(201,344)
\path(402,24)(408,17)
\path(20,3)(14,2)
\path(404,3)(409,2)
\path(10,11)(7,8)
\path(214,348)(213,353)
\path(202,343)(204,349)
\path(406,16)(410,11)
\path(14,4)(10,3)
\path(409,4)(412,3)
\path(9,7)(6,4)
\path(211,351)(210,355)
\path(207,348)(207,353)
\path(408,10)(411,7)
\path(9,7)(6,4)
\path(411,7)(414,4)
\path(10,3)(7,2)
\path(207,353)(208,356)
\path(211,351)(210,355)
\path(409,4)(412,3)
\path(246,0)(173,0)
\path(246,0)(297,0)
\path(86,150)(61,105)
\path(296,213)(271,257)
\path(86,150)(123,213)
\path(296,213)(333,150)
\path(173,0)(122,0)
\path(297,0)(333,0)
\path(61,105)(43,74)
\path(270,257)(252,289)
\path(123,213)(149,257)
\path(333,150)(358,105)
\path(122,0)(86,0)
\path(333,0)(359,0)
\path(43,74)(30,52)
\path(252,289)(240,311)
\path(149,257)(167,289)
\path(358,105)(376,74)
\path(86,0)(60,0)
\path(359,0)(377,0)
\path(30,52)(21,36)
\path(239,310)(231,326)
\path(167,289)(179,311)
\path(376,74)(389,52)
\path(60,0)(42,0)
\path(377,0)(389,0)
\path(21,36)(15,25)
\path(230,326)(224,337)
\path(180,310)(189,326)
\path(389,52)(398,36)
\path(42,0)(29,0)
\path(389,0)(398,0)
\path(15,25)(11,18)
\path(224,337)(220,345)
\path(189,326)(195,337)
\path(397,36)(404,25)
\path(29,1)(21,0)
\path(398,1)(405,0)
\path(11,17)(8,12)
\path(219,344)(216,350)
\path(196,337)(200,345)
\path(403,25)(408,17)
\path(20,1)(14,1)
\path(404,1)(409,1)
\path(9,12)(6,8)
\path(216,349)(214,353)
\path(201,344)(203,350)
\path(408,17)(411,12)
\path(14,2)(10,1)
\path(409,2)(412,1)
\path(7,8)(5,5)
\path(213,353)(212,356)
\path(204,349)(206,353)
\path(410,11)(413,8)
\path(10,3)(7,2)
\path(412,3)(414,2)
\path(6,4)(4,3)
\path(210,355)(210,357)
\path(207,353)(208,356)
\path(411,7)(414,4)
\path(6,4)(4,3)
\path(414,4)(415,3)
\path(7,2)(5,1)
\path(208,356)(208,358)
\path(210,355)(210,357)
\path(412,3)(414,2)
\end{picture}
\caption{\footnotesize}
\end{center}
\caption{\footnotesize
L'appartement de $\text{PGL}_3$ dans le squelette de l'espace de Lubin-Tate
en niveau Iwahori !
}
\end{figure}

\newpage
\appendix

\section{L'immeuble de $\text{PGL}_n$} \label{kujitopo}
\subsection{Définitions}

On fixe $V$ un $F$-espace vectoriel de dimension finie. On note $G =
\text{PGL} (V)$ le groupe algébrique sur $F$ associé.

\begin{defi}
L'immeuble de $G$ est le complexe simplicial $G(F)$-équivariant dont les
sommets sont les classes d'homothéties de réseaux  dans $V$ et dont les
simplexes sont les ensembles finis $(a_1,\dots,a_d)$ de sommets tels
que les $(a_i)_{1\leq i \leq d}$ possèdent des représentants $(\La_i)_{1\leq i\leq
  d}$ tels que 
$$
\pi\La_1 \subsetneq \La_d \subsetneq \dots \subsetneq \La_2 \subsetneq \La_1
$$
On note $\mathcal{I}$ cet ensemble simplicial. 
Pour un réseau $\La$ on notera $[\La]\in\mathcal{I}$ sa classe d'homothétie. 
\end{defi}

Nous commettrons parfois l'abus de notation qui consiste à noter
$\mathcal{I}$ l'ensemble des sommets du complexe simplicial
$\mathcal{I}$ et écrirons ainsi $x\in \mathcal{I}$ pour $x$ un sommet
de $\mathcal{I}$.
\\

Rappelons que le choix d'un sommet $[\La_0]\in \mathcal{I}$ définit un
étiquetage 
des sommets de $\mathcal{I}$ par des éléments de $\Z/n\Z$ :
\begin{eqnarray*}
\mathcal{I} & \ldrt & \Z/n\Z \\
 {[ \La ]} & \longmapsto & [\La_0 : \La] 
\end{eqnarray*}
Le choix d'un autre sommet translate cet étiquetage par une constante
dans $\Z/n\Z$. Si $g\in G(F)$ alors l'action de $g$ translate l'étiquetage 
par une constante égale à $v (\det (g))\in \Z/n\Z$.

\begin{defi}
Une norme sur $V$ est une application $\|.\|: E \ldrt \R\cup \{+\infty \}$
vérifiant
\begin{itemize}
\item $\forall x\in E\;\;\|x\|=+\infty \ssi x=0$
\item $\forall (x,y)\in E^2\;\; \|x+y \| \geq \inf \{ \|x\|,\|y\| \}$
\item $\forall a\in F\;\forall x\in E\;\; \| ax \| = v(a)+ \| x \|$ 
\end{itemize}
\end{defi}

La définition d'une norme que nous prenons n'est pas la définition
usuelle (celle donnée page 58 de \cite{DelHus}). 
Nous prenons la version ``additive'' de la définition usuelle
 (si $\|.\|$ est une norme en notre sens alors la norme usuelle associée est $q^{-\|.\|}$).

\begin{defi}
Deux normes $\|.\|_1, \|.\|_2$  sont équivalentes si il existe $A\in
\R$ tel que $\|.\|_1 = \|.\|_2 +A$. On notera $[\|.\|]$ la classe
d'équivalence de $\|.\|$. On définit une action de $G(F)$ sur les
normes et leurs classes d'équivalence en posant $g.\|.\| = \|g^{-1}.\|$.
\end{defi}

Rappelons également la proposition suivante.

\begin{prop}
La réalisation géométrique $|\mathcal{I}|$ 
de l'immeuble de $G$ 
s'identifie comme ensemble $G(F)$-équivariant 
à l'ensemble 
des classes d'équivalence de normes sur $V$. Si $[\La]$ est une
classe d'homothétie de réseaux dans $V$ le sommet associé dans la 
réalisation géométrique de l'immeuble est la classe de la norme
$$
\|x\|_\La = -\inf \{\;k\in\Z\;|\; \pi^k x\in\La \;\}
$$
\end{prop}

Rappelons également que si $\|.\|$ est une norme alors l'ensemble 
$$
\{\; \left [ \{\; x\in V\;|\; \|x\|\geq a \;\} \right ]\;|\; a\in \R \;\}
$$
forme un simplexe dans l'immeuble et que dans la réalisation
géométrique le point $[\|.\|]$ est dans la réalisation géométrique de ce
simplexe.

Réciproquement, si $\pi\La_1\subsetneq \La_d \subsetneq \dots
\subsetneq \La_1$ définit un simplexe $\s$, si 
$$
\Delta^d= \{\; (t_1,\dots,t_d)\in\R_+^d\;|\; \sum_{i=1}^d t_i =1 \; \}
$$
est le simplexe usuel il y a alors un homéomorphisme 
$$
\Delta^d \iso |\s | 
$$
qui à $(t_1,\dots,t_d)$ associe la norme 
$$
\|.\| = \inf \{ \; \|.\|_{\La_i} + t_1+\dots + t_{i-1}\;\}_{1\leq i
  \leq d}
$$
(où l'on a posé $t_1+\dots + t_{i-1}=0$ si $i=1$). Cet homéomorphisme
est une isométrie si l'immeuble est muni de sa métrique définie à
partir de celle sur ses appartements et $\Delta^d$ de la métrique 
$\frac{1}{\sqrt{2}}\sqrt{\sum_i t_i^2 }$.  

\subsection{Appartements}
\subsubsection{Définition}

\begin{defi}
Soit $T$ un tore maximal déployé dans $G$. L'appartement associé à $T$
est le sous-complexe simplicial dont les sommets sont les points fixes
de $T(F)^1= \dpt{\bigcap_{\chi\in X^* (T)} \ker |\chi|}$. On le note
$\mathcal{A}(T)$.
\end{defi}

La réalisation géométrique $|\mathcal{A}(T)|$ de $\mathcal{A}(T)$  est
le sous-espace fermé de $|\mathcal{I}|$ formé des points fixes de $T(F)^1$. 

L'ensemble des tores maximaux déployés dans $G$ est en bijection avec
l'ensemble des n-uplets de droites $(D_1,\dots,D_n)$ dans $V$ tels que 
$V=D_1\oplus\dots\oplus D_n$. A $(D_1,\dots,D_n)$ est associé le tore 
$T=\Gm^n/\Gm \hookrightarrow G$ agissant diagonalement relativement à
la décomposition en sommes de droites de $V$. 
 L'ensemble des appartements est donc en bijection avec l'ensemble
 des classes d'équivalences de bases $(e_1,\dots,e_n)$ de $V$ où 
$\forall (x_i)_i \in (F^\times)^n\;\;
(e_1,\dots,e_n)\sim (x_1 e_1,\dots,x_n e_n)$. A la base
$(e_1,\dots,e_n)$ est associé le complexe simplicial dont les sommets
sont 
$$
\{\; \left [ <\pi^{a_1}e_1,\dots, \pi^{a_n} e_n >\right ]\;|\;
(a_i)_i\in \Z^n\; \}
$$
qui est l'appartement associé.

\begin{exem}
 Pour $\text{PGL}_2$ les appartements sont les droites
     simpliciales infinies à droite et à gauche dans l'arbre.
\end{exem}

\begin{figure}[htbp]
   \begin{center}
      \includegraphics{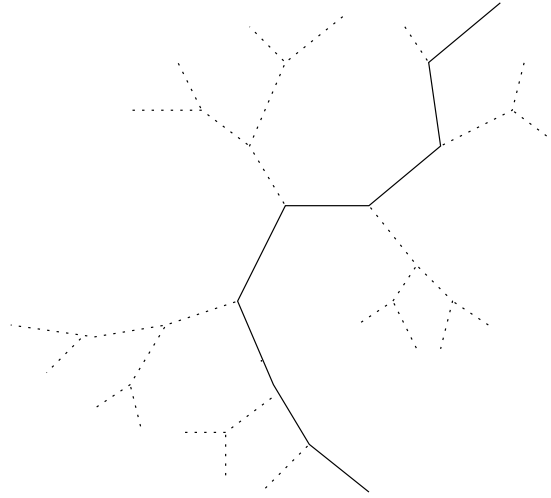}
   \end{center}
   \caption{\footnotesize Un appartement  dans l'immeuble de $\text{PGL}_2$
     sur $\Q_2$.
     }
   \label{figklsm}
\end{figure}

\subsubsection{Structure affine}

Le groupe $T(F)$ stabilise l'appartement $\mathcal{A}(T)$ et fait des
sommets de $\mathcal{A}(T)$ un espace principal homogène sous
$T(F)/T(F)^1$. Il y a un isomorphisme 
\begin{eqnarray*}
X_* (T) & \iso & T(F)/T(F)^1 \\
\omega & \longmapsto & \omega (\pi) T(F)^1
\end{eqnarray*}
Si $[\La]\in \mathcal{A}(T)$ alors si $\underline{G}$ désigne le
modèle entier de $G$ associé, $\underline{G}= \text{PGL} (\La)$, si $\underline{T}$
désigne l'adhérence schématique de $T$ dans $\underline{G}$,
$\underline{T}$ est un sous-tore de $\underline{G}$. Alors, 
$T(F)^1= \underline{T} (\O_F)$ et le choix de $[\La]$ induit  au niveau des sommets 
le diagramme suivant 
$$
\xymatrix{
T(F)/\underline{T}(\O_F) \ar[r]^(.6)\sim \ar@{^(->}[d] &  \mathcal{A} (T) \ar@{^(->}[d]
\\
G(F)/\underline{G} (\O_F)  \ar[r]^(.6)\sim & \mathcal{I}
}
$$
qui identifie les sommets de $\mathcal{A} (T)$ avec l'immeuble de
$T$. 

L'action de $X_* (T)$ sur les sommets de $\mathcal{A}(T)$ s'étend en
une structure d'espace affine sur $|\mathcal{A} (T)|$
 d'espace vectoriel sous-jacent
$X_*(T)\otimes \R$. Si $T$ est associé à la base $(e_1,\dots,e_n)$ de
$V$, $X_* (T) \simeq \Z^n/\Z$, $\Z$ agissant diagonalement sur $\Z^n$ et
 où si $(\e_i)_i$ est la base canonique de
$\Z^n$, $\e_i$ correspond au co-caractère $t\longmapsto [e_j \mapsto
t^{\delta_{ij}} e_j]$. En choisissant comme origine
$[<e_1,\dots,e_n>]$ dans $|\mathcal{A} (T) |$, 
$$
|\mathcal{A}
(T)|\simeq X_* (T)\otimes \R \simeq \R^n/\R
$$
et à $(a_1,\dots,a_n)\in \R^n/\R$ est associé la classe de la norme 
$$
\|\sum_{i=1}^n \l_i e_i \|=\inf \{ v(\l_i)-a_i \;|\; 1\leq i\leq n\; \}
$$

Et en fait on a la réciproque  suivante :

\begin{lemm}\label{katea}
Une classe d'équivalence de normes $[\|.\|]$ est dans l'appartement
associé à la base $(e_1,\dots,e_n)$ ssi $\forall (x_1,\dots,x_n)\in
F^n$ $\;\; \| \sum_{i=1}^n x_i e_i \|$ ne dépend que des $v(x_i)$ pour
$1\leq i \leq n$.
\end{lemm}
\dem 
Soit $\|.\|$ une norme comme dans l'énoncé. Il suffit de montre que pour tout 
$a\in \R$ $\;\{\; x\in V \;|\; |\ x\| \geq a \;\}$ est un réseau de la
forme $<\pi^{a_1} e_1,\dots ,\pi^{a_n} e_n>$ puisque $[\|.\|]$
est dans l'enveloppe convexe de ces classes de réseaux lorsque $a$ varie dans
$\R$.
Cela résulte du lemme qui suit dont la démonstration élémentaire est
laissée au lecteur.
\qed

\begin{lemm}\label{kopghut}
Soit un sous-ensemble $\Delta \subset ( \Z \cup  \{+\infty \}
)^n$. L'ensemble 
$$
\{ \; (x_i)_{1\leq i \leq n} \in F^n \;|\; v(x_i)\in \Delta \;\}
$$
est un réseau de $F^n$ ssi il existe $(a_1,\dots,a_n)\in \Z^n$ tels
que $\Delta = \{\; (k_1,\dots,k_n)\;|\; \forall i \; k_i \geq a_i \;\}$
\end{lemm}

\subsubsection{Le système de racines affines}

On a vu qu'on peut retrouver les sommets de l'appartement associé au
tore $T$ à partir de celui-ci et
mettre une structure affine dessus. En fait, à partir du système de
racines associé on peut retrouver toute la structure simpliciale. 
\\

Soit $\Phi \subset X^* (T)$ le système de racines de $T$ dans $\Lie\,
G$. On appelle racine affine l'ensemble des fonctions 
$$
\Phi_{\text{aff}} = \{\; \a+k\;|\; \a\in \Phi\; k\in \Z\;\}
$$
Soit $x_0\in |\mathcal{A} (T)|$ un sommet (un point spécial dans la
terminologie de Bruhat-Tits, mais pour $\text{PGL}_n$ tous les sommets
sont spéciaux). Le choix de $x_0$ permet d'identifier $|\mathcal{A}
(T)|$ à $X_* (T)\otimes \R$ et fournit donc un ensemble de fonctions
affines encore appelées racines affines sur $|\mathcal{A} (T) |$ et
noté $\Phi_{\text{aff}}$. On vérifie aussitôt que cet ensemble de fonctions
affines ne dépend pas du choix de $x_0$.

\begin{defi}
Un mur est un hyperplan de la forme $\a^{-1} (\{0\})$ où $\a\in
\Phi_{\text{aff}}$.
Un demi-appartement est un ensemble de la forme $\a^{-1}([0,+\infty [)$.
\end{defi}

\begin{exem}
Dans l'immeuble de $\text{PGL}_2$ les demi-appartements sont les
demi-droites simpliciales.
\end{exem}

La structure simpliciale se retrouve maintenant de la façon suivante.

\begin{defi}
Une facette est une classe dans $|\mathcal{A} (T)|$ pour la relation
d'équivalence $x\sim y$ si $x$ et $y$ sont contenus dans les même
demi-appartements. On met la relation d'ordre suivante sur les
facettes : $F_1<F_2$ si $F_1 \subset \overline{F}_2$.  
\end{defi}

Appelons simplexe ouvert $\{(x_i)_i \in \R_+^d\; |\; \sum_i x_i=1
\text{ et } \forall i \; x_i\neq 0\;\}$. Ainsi pour un complexe
simplicial $\mathcal{C}$, et $\s$ un simplexe de $\mathcal{C}$ on peut
définir le simplexe ouvert associé à $\s$ dans la réalisation
géométrique $|\s |\subset |\mathcal{C}|$. Ces simplexes ouverts
forment une partition de $|\mathcal{C}|$. 

\begin{Fait}
Les facettes s'identifient aux simplexes ouverts dans la réalisation
géométrique du complexe simplicial $\mathcal{A}(T)$. Ce complexe
simplicial s'identifie au complexe simplicial associé à l'ensemble
ordonné des facettes.
\end{Fait}

Soit $\Delta$ un ensemble de racines simples dans $\Phi$
associé à un ordre sur les racines et
$\{\omega_1,\dots,\omega_{n-1}\}$ l'ensemble des copoids fondamentaux
associé à $\Delta$. Alors, $X_* (T) = <\omega_1,\dots,\omega_{n-1}>$
et un simplexe maximal associé est 
$$
\text{Conv} (0,\omega_1,\dots,\omega_{n-1}) = \{\; \sum_{i=1}^{n-1} a_i
\omega_i\;|\; a_i\in \R_+ \text{ et } \sum_i a_i \leq 1 \;\} \subset X_* (T)\otimes \R
$$
dont les sommets sont $(0,\omega_1,\dots,\omega_{n-1})$. Les autres
simplexes (les adhérences des facettes) s'en déduisent en prenant 
les itérés sous le groupe de Weyl affine des faces de ce simplexe.

\begin{figure}[htbp]
   \begin{center}
      \includegraphics{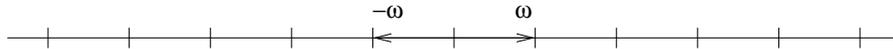}
   \end{center}
   \caption{\footnotesize Un appartement  de l'immeuble de $\text{PGL}_2$. $\omega =$
     copoids fondamental
     }
   \label{fiportj}
\end{figure}

\begin{figure}[htbp]
   \begin{center}
      \includegraphics{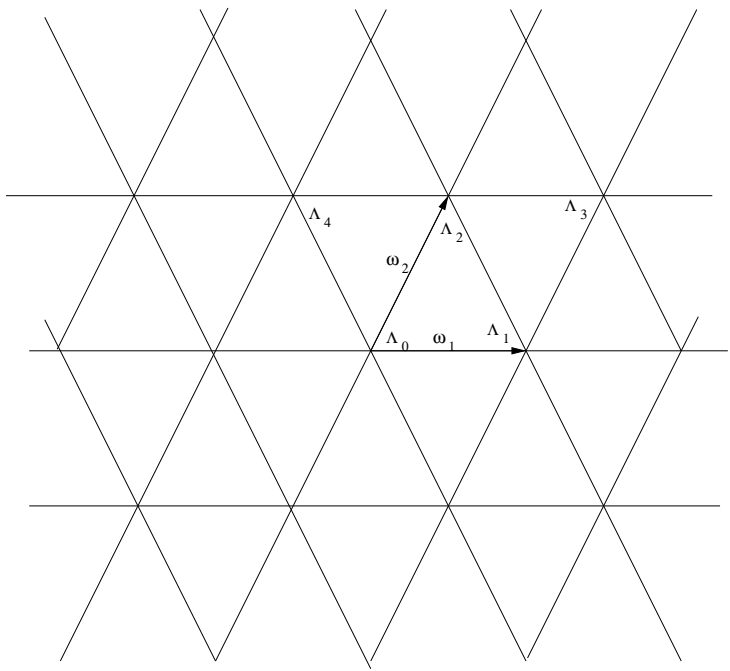}
   \end{center}
   \caption{\footnotesize Un appartement  dans l'immeuble de $\text{PGL}_3$
     }
   \label{fiopmds}
\end{figure}

Dans la figure ci-dessous $\omega_1 (t)= (t,1,1)$ et $\omega_2
(t)=(t,t,1)$ sont les copoids fondamentaux associés au système de
racines simples $\{\a_1,\a_2\}$ pour le tore diagonal de
$\text{PGL}_3$ où $\a_1 (t_1,t_2,t_3)= t_1 t_2^{-1}$ et $\a_2
(t_1,t_2,t_3)= t_2 t_3^{-1}$. Les réseaux indiqués sont 
$\La_0=<e_1,e_2,e_2>$, $\La_1=<\pi e_1,e_2,e_3>$, $\La_2=<\pi e_1,\pi e_2,
e_3>$, $\La_3=<\pi^2 e_1, \pi e_2, e_3>$, $\La_4=< e_1,\pi e_2 ,
e_3>$. Le choix fait de la structure euclidienne pour représenter
l'appartement de $\text{PGL}_3$ n'est pas quelconque. On a pris, à un
scalaire près, 
l'unique structure euclidienne invariante sous l'action du groupe de Weyl
affine (cf. la section qui suit). Pour cette métrique les simplexes
sont des triangles équilatéraux.

\subsubsection{Groupe de Weyl affine}

Soit $T$ un tore déployé maximal dans $G$.

\begin{defi}
On note $W_{\text{aff}} =N(T) (F)/T(F)^1$.
\end{defi}

Le groupe $W_{\text{aff}}$ agit sur $|\mathcal{A}(T)|$ par des
transformations affines. Il y a de plus une suite exacte scindable 
$$
1 \ldrt X_* (T) \ldrt W_{\text{aff}} \ldrt W \ldrt 1
$$
où $W=N(T)(F)/T(F)$ est le groupe de Weil vectoriel du système de
racines $\Phi$. L'application $W_{\text{aff}}\ldrt W$ s'identifie à
l'application qui à une application affine associe sa partie
vectorielle. Le choix d'une origine $x_0\in\mathcal{A}(T)$ induit un
scindage $W_{\text{aff}} \simeq X_* (T)\rtimes W$ via $W=
\text{Stab}_{W_{\text{aff}}} (x_0)$. L'action de $W_{\text{aff}}$ sur
$|\mathcal{A} (T)|$ permet également de retrouver la structure
simpliciale puisque les murs sont exactement les hyperplans fixes des
symétries non-triviales dans $W_{\text{aff}}$. Ce complexe simplicial
ne dépend que du système de Coxeter associé à $\Phi_{\text{aff}}$. 

On muni $|\mathcal{A} (T)|$ d'une métrique euclidienne invariante sous 
l'action du groupe de Weyl affine (une telle métrique est unique à un
scalaire près). 
 Bruhat et Tits montrent qu'une
telle métrique peut s'étendre à tout l'immeuble en une métrique 
$G(F)$-invariante de manière à redonner 
la métrique que l'on a fixée en restriction aux appartements. 
On renvoie à la section 3 du chapitre 3 de  \cite{DelHus} pour une
formule donnant la distance entre deux normes dans l'immeuble. 

\subsubsection{Une paramétrisation des chambres dans un
  appartement}\label{rutoipm}

\begin{figure}[htbp]
   \begin{center}
      \includegraphics{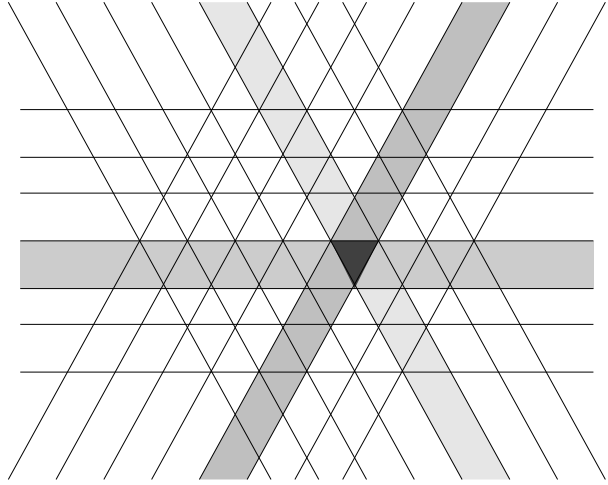}
   \end{center}
   \caption{\footnotesize Un simplexe paramétré par une intersection de demi-appartements
     }\label{rutypo}
\end{figure}

Soit $\mathcal{A} (T)$ un appartement dans $\mathcal{I}$. Soit $\Phi
\subset X^* (T)$ l'ensemble des racines associées et $\Phi^+$ un
ensemble de racines positives associé au choix d'un sous-groupe de
Borel contenant $T$. Nous allons donner une paramétrisation des
simplexes maximaux dans $\mathcal{A}(T)$ comme intersection de
demi-appartements (cf. figure \ref{rutypo}).

\begin{prop}
Fixons une origine dans $\mathcal{A} (T)$ et identifions donc les sommets
de $\mathcal{A} (T)$ à $X_* (T)$. Soit 
$$
A= \{ (b_\a)_{\a\in \Phi_+}\in \Z^{\Phi_+}\;|\; \forall
\a,\b\in\Phi_+\; \a+\b \in \Phi_+ \impl b_{\a+\b} \in \{ b_\a + b_\b,
b_\a+ b_\b +1 \} \;\}
$$
et $B$ l'ensemble des simplexes maximaux (vus comme les collections de leurs
sommets) dans $\mathcal{A} (T)$. Il y a alors une bijection entre $A$
et $B$ donnée par les deux applications inverses 
\begin{eqnarray*}
A & \ldrt & B \\
(b_\a)_{\a\in \Phi_+} & \longmapsto & \{ x\in X_* (T)\;|\; \forall \a
\in \Phi_+\; \a (x) \in \{ b_\a, b_\a+1\}\;\}
\end{eqnarray*}
et
\begin{eqnarray*}
B & \ldrt & A \\
S &\longmapsto & \left ( \inf \{\a (x)\;|\; x\in S \} \right )_{\a \in \Phi_+}
\end{eqnarray*}
La réalisation géométrique du simplexe associé à $(b_\a)_\a$ est 
$$
\{ x\in X_* (T)\otimes \R \;|\; \forall \a \in \Phi_+ \; \a (x) \in
[b_\a, b_\a +1]\;\}
$$
\end{prop}
\dem
On peut supposer que $T$ est le tore diagonal de $\text{PGL}_n$
$$
T=\{\text{diag} (t_1,\dots,t_n)\;|\; t_i\in \Gm \}
$$
et que $\Phi_+ = \{ t_i t_j^{-1}\; |\; i<j \}$. On identifie alors
$X_* (T)$ à $\Z^n/\Z$ où $\Z\hookrightarrow \Z^n$ est le plongement
diagonal. Soit donc une collection d'entiers relatifs $(b_{ij})_{1\leq
  i < j \leq n}$ telle que
$$
\forall i<j<k \;\; b_{ik}\in \{ b_{ij}+b_{jk},b_{ij}+b_{jk}+1\} 
$$
On veut montrer que 
$$
S=\{ (a_1,\dots,a_n)\in\Z^n/\Z\;|\; \forall i<j\; a_j-a_i \in\{ b_{ij}
, b_{ij}+1 \}
$$
est un simplexe maximal et que cela en donne une paramétrisation. 
Le point $P=(0,b_{12},\dots, b_{1i},\dots, b_{1n})$ appartient à
$S$. On vérifie facilement que quitte à translater par $-P$ et
transformer les $b_{ij}$ en $b_{ij}-b_{1j}+b_{1i}$ pour $1<i<j\leq n$
on doit montrer que les simplexes maximaux possédant $(0,\dots,0)$
comme sommet sont paramétrés par 
$$
\{ (b_{ij})_{1\leq i<j\leq n}\;|\; \forall i>1\; b_{1i}=0 \text{ et }
\forall 1\leq i<j<k\; b_{ik}\in \{ b_{ij}+b_{jk},b_{ij}+b_{jk}+1 \} \;\}
$$
Cet ensemble s'identifie à $$
X=\{ (b_{ij})_{1<i<j\leq n}\;|\; \forall i<j\; b_{ij}\in \{ -1,0
\}\text{ et } \forall i<j<k\; b_{ik}\in
\{b_{ij}+b_{jk},b_{ij}+b_{jk}+1\}\; \}
$$
Soit le simplexe maximal
\begin{eqnarray*}
S_0 & = &  \{(0,\dots, 0 ,\underbrace{1,\dots,1}_{i})\;|\; 0 \leq i <n \}
\\
&=& \{ (a_1,\dots,a_n)\;|\; \forall i<j \; a_j-a_i\in\{0,1 \}\;\}
\end{eqnarray*}
qui est associé à la donnée $\forall i<j\; b_{ij}=0$.
Les autres simplexes maximaux possédant $(0,\dots,0)$ comme sommet
sont les orbites sous le groupe de Weyl $\mathfrak{S}_n$ de $S_0$. De
plus
\begin{eqnarray*}
\forall \s \in \mathfrak{S}_n\;\; \s.S_0 &=& \{
(a_1,\dots,a_n)\in\Z^n/\Z\;|\; \forall i<j\; a_{\s^{-1} (j)}-
a_{\s^{-1} (i)} \in \{0,1 \}\;\} \\
&=& \{ (a_1,\dots,a_n)\in\Z^n/\Z\;|\; \forall i<j\; a_j-a_i\in
\{0,1\}\text{ si } \s^{-1}(i)<\s^{-1} (j) \\
&& \hspace{3.7cm} \text{ et } 
 a_j-a_i\in
\{-1,0\}\text{ si } \s^{-1}(i)>\s^{-1} (j)\;\}
\end{eqnarray*}
Soit donc l'application qui à $\s\in\mathfrak{S}_n$ associe
$(b_{ij})_{1\leq i<j\leq n}$ où
$$
b_{ij}= \left \{ 0 \text{ si } \s^{-1} (i)< \s^{-1} (j) \atop
-1 \text{ si } \s^{-1} (i)> \s^{-1} (j) 
\right.
$$
Il est asié de vérifier que cela définit une application
$\mathfrak{S}_n \ldrt X$ qui d'après le lemme qui suit est une
bijection.
\qed

\begin{lemm}
Soit l'application de $\mathfrak{S}_n$ dans l'ensemble des parties de 
$\{(i,j)\;|\; 1 \leq i<j\leq n\}$ qui à $\s$ associe 
$\{i<j\;|\; \s (i) < \s (j)\}$. Elle induit une bijection entre 
$\mathfrak{S}_n$ et
$$
\left \{ B\; | \; { (i,j)\in B\text{ et } (j,k)\in B \impl (i,k)\in B \atop
 (i,j)\notin B\text{ et } (j,k)\notin B \impl (i,k)\notin B}
\;\right \}
$$
\end{lemm}

\subsection{Immeuble dual}\label{immdual}

Notons $\mathcal{I} (V)$ l'immeuble de $\text{PGL} (V)$ et
$\mathcal{I} (V^*)$ celui de $\text{PGL} (V^*)$. 

\begin{defi}
Pour  $\La$ un réseau dans $V$ notons  $\La^\vee =\{ \ph\in V^*\;|\; \ph
(\La) \subset \O_F \}$ le réseaux dual.
\end{defi}
 
L'isomorphisme de complexes
simpliciaux 
\begin{eqnarray*} 
\mathcal{I} (V) & \ldrt & \mathcal{I} (V^*) \\
{[\La]} & \longmapsto &
[\La^\vee]
\end{eqnarray*}
  s'étend en
une isométrie $|\mathcal{I} (V) | \iso |\mathcal{I} (V^*)|$ en posant
pour toute norme $\|.\|$ sur $V$ 
$$
\forall \ph\in V^*\;\;
\|\ph \|^\vee = \inf \{ v(\ph (x)) + \| x\| \; |\; \|x\|\geq 0 \}
$$
qui définit une norme $\|.\|^\vee$ sur $V^*$, et définit l'isométrie
$[\|.\|] \longmapsto [\|.\|^\vee ]$. Cette isométrie transforme
l'action de $g\in \text{PGL} (V)$ en celle de $(\,^t g)^{-1}$ et envoie
les appartements sur des appartements via des isométries affines.

\subsection{Quartiers}

\begin{defi}
Un quartier dans l'appartement $|\mathcal{A} (T)|$ de $|\mathcal{I}|$
est un sous-ensemble fermé dans $|\mathcal{A}(T)|$ de la forme $x+C$
où $C$ est une chambre vectorielle dans $X_*(T)\otimes \R$, l'espace vectoriel directeur
de l'espace affine $|\mathcal{A}(T)|$. L'élément $x$ est appelé le
sommet du quartier $x+C$.  
\end{defi}

\begin{figure}[htbp]
   \begin{center}
      \includegraphics{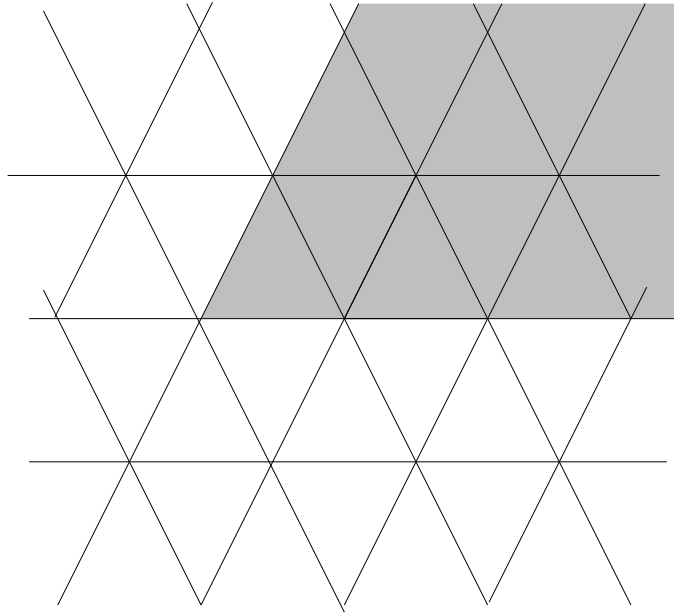}
   \end{center}
   \caption{\footnotesize Un quartier  dans l'immeuble de $\text{PGL}_3$
     }
   \label{quartierfigure}
\end{figure}

On renvoie à la figure \ref{quartierfigure}. 

Soit $[\La_0]\in\mathcal{I}$. Soit $Q$ un quartier de sommet
$[\La_0]$ dans $|\mathcal{I}|$. Il existe alors une base
$(e_1,\dots,e_n)$ de $\La_0$ telle que les sommets de $\mathcal{I}$
dans $Q$ soient
$$
\{\; \left [ <\pi^{a_1} e_1,\dots,\pi^{a_n} e_n >\right ] \;|\;
  a_1\geq \dots \geq a_n \; \}
$$
Si $Q$ est un quartier de sommet $x$ dans l'appartement $|\mathcal{A}(T)|$
il existe alors un ensemble de
co-poids fondamentaux $(\omega_1,\dots,\omega_{n-1})$ dans $X_* (T)$ tel
que 
$$
Q= x+   \R_+ \omega_1 + \dots + \R_+ \omega_{n-1}
$$
Les quartiers de sommet $x$ dans l'appartement $\mathcal{A} (T)$ sont en bijection avec les sous-groupes de
Borel de $G$ contenant $T$. A un sous-groupe de Borel est associé un
ensemble de racines positives, donc de racines simples et  de co-poids fondamentaux.

\subsection{Enclos}\label{enclosA}

\begin{defi}[\cite{BruhatTits1}]
Soit $M$ un sous-ensemble fini de sommets  contenu dans un
appartement $\mathcal{A}(T)$.
On appelle enclos de $M$ l'intersection des demi-appartements dans
$|\mathcal{A}(T)|$ contenant $M$. 
\end{defi}

A priori cette définition est ambiguë puisque l'enclos de $M$ semble
dépendre du choix d'un appartement dans lesquels $M$ est
contenu. Cependant ce n'est pas le cas grâce à la proposition
suivante.

\begin{prop}
Soit $M$ comme dans la définition précédente. Alors l'enclos de $M$
est égal à l'intersection des appartements dans lesquels $M$ est
contenu. Il est aussi égal à la réalisation géométrique de l'enveloppe
convexe simpliciale de $M$. 
\end{prop}
\dem
Soient $C_1$ l'enclos de $M$, $C_2$ l'intersection des appartements
contenant $M$ et $C_3$ la réalisation géométrique de l'enveloppe
convexe simpliciale. Il est facile de vérifier ``à la main'' que les
demi-appartements sont des intersections de deux appartements. De plus
l'intersection de deux appartements ainsi que les demi-appartements
sont simplicialement convexes. On en déduit les inclusions 
$$
C_3  \subset C_2 \subset C_1
$$
Il suffit donc de montrer que si $C$ est la réalisation géométrique
d'un sous-complexe simplicial convexe fini d'un appartement et
$x\notin C$ un
sommet dans cet appartement  il existe alors $\a\in \Phi_{\text{aff}}$
telle que $\a (x)\leq 0$ et $\forall y\in C \; \a (y)\geq 0$. Bien que
pénible à vérifier cela ne pose pas de problèmes.
\qed

\begin{figure}[htbp]
   \begin{center}
      \includegraphics{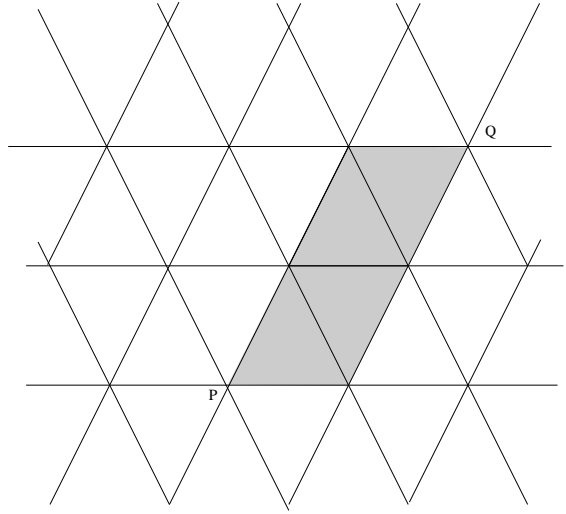}
   \end{center}
   \caption{\footnotesize L'enclos délimité par les points $P$ et $Q$ dans l'immeuble de $\text{PGL}_3$
     }
   \label{Enclosfigure}
\end{figure}

\begin{rema}
Pour l'immeuble de $\text{PGL}_2$ l'enclos coïncide avec l'enveloppe
convexe pour la structure affine de l'appartement, mais en général
cette enveloppe convexe affine est contenue strictement dans
l'enveloppe convexe simpliciale. 
\end{rema}

\begin{prop}\ref{Fokkomp}
Soit $S$ un ensemble fini de sommets contenus dans un même appartement. Une classe de normes $[\|.\|] \in |\mathcal{I}|$ est dans l'enclos délimité par les éléments de $S$ ssi pour un
ensemble 
de représentants $\widetilde{S}$ des classes d'homothéties des
réseaux dans  $S$ la fonction $x\mapsto \| x\|$ ne dépende que des $
\|x\|_\La$, $\La\in \widetilde{S}$. 
\end{prop} 
\dem
Si les éléments de $S$ sont contenus dans l'appartement associé
à la base $(e_1,\dots,e_n)$ de $V$ alors les $\|.\|_{\La}$ ne
dépendent que des valuations des coordonnées dans une telle base. La
proposition résulte donc du lemme \ref{katea} couplé à la proposition précédente.
\qed

\begin{defi}\label{ifgopmu}
Soient $\|.\|_1$ et $\|. \|_2$ deux normes sur $V$. Posons $$\|.\|_1
\wedge \|.\|_2  =
\inf \{ \|.\|_1, \|.\|_2 \}$$ et 
$$
\|.\|_1\vee \|.\|_2 = \left ( \|.\|_1^\vee \wedge \|.\|_2^\vee \right )^\vee
$$
qui définit deux opérations binaires sur les normes sur $V$.
\end{defi}

Ces deux opérations étendent les opérations d'intersection et de sommes de deux réseaux : 
$$
\|.\|_{\La_1} \wedge \|.\|_{\La_2} = \|.\|_{\La_1\cap \La_2}, \;\;\; \|.\|_{\La_1} \vee \|.\|_{\La_2} = \|.\|_{\La_1 + \La_2}
$$

On peut alors montrer la proposition suivante :

\begin{prop}
 L'enclos d'un ensemble fini de sommets est l'ensemble des classes de
 normes $[\|.\|]$ telle que $\|.\|$ puisse s'obtenir à partir des
 $\|.\|_\La$, $[\La]\in S$ et l'itération d'un nombre fini d'opérations qui sont 
$$
\forall A\in \R\;\;
\|.\| \mapsto \|.\| + A \text{ , } \vee \text{ et } \wedge
$$
\end{prop}

\subsection{Quotients de l'immeuble}

\subsubsection{Quotient par un sous-groupe compact maximal de $\GL_n$} \label{Qdch}

Soit $[\La_0] \in \mathcal{I}$ et $\GL (\La_0)$ ls sous-groupe compact
de $\GL(V)$ associé. Soit $Q$ un quartier de sommet $[\La_0]$ dans un
appartement contenant $[\La_0]$. Alors, 
$$
Q \hookrightarrow |\mathcal{I} | \twoheadrightarrow \GL (\La_0)
\bc |\mathcal{I} |
$$
induit un isomorphisme 
$$
Q \iso \GL (\La_0) \bc |\mathcal{I} |
$$
qui au niveau des sommets n'est rien d'autre que la décomposition de
Cartan de $\GL_n$
$$
\GL_n (F)= \coprod_{a_1\geq \dots \geq a_n} \GL_n (\O_F)
\text{diag}(\pi^{a_1},\dots,\pi^{a_n}) \GL_n (\O_F)
$$
Cela définit une projection 
$$
\text{pr}_Q : |\mathcal{I}| \ldrt Q
$$
Si $Q'$ est un autre quartier de sommet $[\La_0]$ (pas forcément dans
le même appartement) 
alors la restriction de $\text{pr}_Q$ à $Q'$ induit une isométrie
affine (au sens où elle conserve les barycentres) entre $Q'$ et $Q$. 
De plus, si on fixe un étiquetage des sommets $\e: \mathcal{I} \ldrt
\Z/n\Z$ alors la projection $\text{pr}_Q$ conserve cet
étiquetage. 
Ainsi si $Q$, resp. $Q'$, est contenu dans l'appartement
$|\mathcal{A}(T)|$, resp . $|\mathcal{A} (T')|$, si 
$$
Q = [\La_0] + \R_+ \omega_1 +\dots + \R_+ \omega_{n-1}, \text{ resp. } 
Q' = [\La_0] + \R_+ \omega_1' +\dots + \R_+ \omega_{n-1}'
$$
pour des co-poids fondamentaux $(\omega_i)_i$ dans $X_* (T)$ , resp.
$(\omega_i')_i$ dans $X_* (T')$, il existe alors une permutation $\s
\in \mathfrak{S}_{n-1}$ telle que 
$$
\forall i\;\; \e ([\La_0]+ \omega_i') = \e ([\La_0] + \omega_{\s (i)})
$$
Cette permutation se calcule en utilisant le fait que $\text{pr}_Q$ conserve un étiquetage des sommets. 
Alors, via les identifications précédentes 
\begin{eqnarray*}
\text{pr}_{Q|Q'} : Q' & \ldrt & Q \\
\sum_{i=1}^{n-1} x_i \omega_i' & \longmapsto & \sum_{i=1}^{n-1} x_i
\omega_{\s (i)}
\end{eqnarray*}
Parmi les autres propriétés de cette projection on vérifie que si $S$
est un ensemble de sommets contenus dans un même quartier d'un
appartement de $\mathcal{I}$ alors l'image par $\text{pr}_Q$ de
l'enclos de $S$ est l'enclos de $\text{pr}_Q (S)$. 

\subsection{Quotient par un sous-groupe d'Iwahori}

Soit $S$ un simplexe maximal dans $\mathcal{I}$ et 
$$
I= \{ g\in G(F)\;|\; \forall x\in S\; g.x=x \}
$$
le sous-groupe d'Iwahori associé.

Soit $\mathcal{A}$ un appartement contenant $S$. L'application
composée 
$$
|\mathcal{A} | \hookrightarrow |\mathcal{I}| \twoheadrightarrow I\bc |\mathcal{I}|
$$
induit un isomorphisme
$$
|\mathcal{A}| \iso I\bc |\mathcal{I}|
$$

\section{Ramification supérieure/inférieure dans le cas monogène}\label{Ramitofm}

Dans cet appendice nous explicitions les définitions d'Abbes et Saito
dans le cas des algèbres localement intersection compléte monogènes,
ce qui est le cas des sous-schémas en groupes finis localement libres des groupes
formels de dimension $1$.

Notons $f$ un polynôme unitaire séparable à coefficients dans
un corps valué complet non-archimédien $K$ de valuation $v: K\ldrt \R\cup \{+\infty \}$.
 On supposera de plus que $f\in \O_K [T]$.  
Soit $\overline{K}$ une clôture algébrique de $K$.
On note
$$
f(T) = \prod_{i\in I} (T-\a_i)
$$
où $\forall i\; \a_i\in \O_{\overline{K}}$.

La fonction $f$ définit un morphisme étale fini d'espaces rigides 
$$
f: \mathbb{B}^1\xrig{\text{ étale fini }} \mathbb{B}^1
$$
Lorsque $f(0)=0$ on va s'intéresser à la façon dont varient les
composantes connexes géométriques des images réciproques des boules de
rayon $\e$ lorsque $\e$ varie dans $[01]$ 
$$
f^{-1} ( \mathbb{B} (0,\e)) \ldrt \mathbb{B} (0,\e)
$$
qui forme une famille de revêtements étales finis lorsque $\e$
varie. On va voir que les composantes connexes géométriques sont des
boules et que si $f^{-1} ( \mathbb{B} (0,\e))^0$ désigne la
composantes connexe neutre de $0$ alors
$$
f^{-1} ( \mathbb{B} (0,\e))^0 = \mathbb{B} (0,\psi (\e)) \xrig{\; f\;} \mathbb{B}(0,\e)
$$
où $\psi$, une fonction de Herbrand, est une fonction convexes affine par morceaux qui se calcule
en termes du polygone de Newton de $f$. 

Cela définit une filtration sur les racines de $f$  via 
$$
\e \longmapsto f^{-1} (\{0 \}) \cap f^{-1} ( \mathbb{B} (0,\e))^0
$$
qui lorsque $G=\spec (\O_K[T]/(f(T)))$ est un schéma en groupes 
fini localement libre 
de
section unité donnée par $T=0$ forme une filtration par des sous-groupes
étales de $G_\eta$.
Un théorème de fibration de Abbes et Saito montre plus généralement
que cete filtration peut être définie indépendamment d'une
présentation (i.e. dans le cas monogène le choix d'un polynome $f$)
pour des algèbres syntomiques finies.  Nous explicitions donc cete
filtration dans le cas particulier très simple des algèbres monogènes.

\subsection{Composantes connexes géométriques de $\{x\;|\; v(f(x))\geq \e\}$ }

\begin{defi}
Notons pour tout $\e\in v(\overline{K}^\times)_{\geq 0}\;\; Y^\e= \{ x\in \mathbb{B}^1
(\overline{K}) \;|\; 
v(f(x))\geq \e \}$ comme espace rigide sur $K$ (un ouvert affinoïde dans la boule unité). 
\end{defi}

Si $\mathbb{B} (0,\e) =\{x|v(x)\geq \e\}$
on a donc en termes d'espaces rigides $Y^\e = f^{-1} (\mathbb{B}
(0,\e))$ où $f:\mathbb{B}^1\ldrt \mathbb{B}^1$. 

\begin{defi}
Soit $\mathcal{R}_{0,\e}$ la relation d'équivalence  sur $I$ définie par
$\forall i,j\in I \; i\underset{\mathcal{R}_{0,\e}}{\sim} j$. 
Définissons par récurrence sur $k$ la relation d'équivalence
$\mathcal{R}_{k,\e}$ sur $I$ de la façon suivante :
$\mathcal{R}_{k,\e}$ est plus fine que $\mathcal{R}_{k-1,\e}$ définie
par 
$$
\forall i,j\in I\;\; i\underset{\mathcal{R}_{k,\e}}{\sim} j \text{ si } 
  i\underset{\mathcal{R}_{k-1,\e}}{\sim} j \text{ et } v(\a_i-\a_j) \geq
  \frac{1}{\text{Card}[i]_{ \mathcal{R}_{k-1,\e}}}
\left ( \e -\sum_{j' \underset{\mathcal{R}_{k-1,\e}}{\nsim} j}
  v(\a_{j'}-\a_j) \right )
$$
où $[i]_{ \mathcal{R}_{k-1,\e}}$ désigne la classe d'équivalence de $i$ pour la
relation $ \mathcal{R}_{k-1,\e}$. 
On note $\mathcal{R}_{\infty,\e}=\mathcal{R}_{k,\e}$ pour $k>>0$ la relation d'équivalence limite sur
$I$. 
\end{defi}

\begin{lemm}
Notons pour tout $r\in v(\overline{K}^\times)$ et $\a\in \overline{K}$ $\;\mathbb{B}(\a,r)= \{ x\;|\; v(x-\a)\geq r \}$
comme espace rigide sur $\widehat{\overline{K}}$. Les composantes connexes géométrique de $Y^\e$ sont des boules  
$$
Y^\e\otimes_K \widehat{\overline{K}} = \coprod_{[i]\in I/\mathcal{R}_{\infty,\e}} \mathbb{B}\left (\a_i,
\frac{1}{[i]_{\mathcal{R}_{\infty,\e}}} (
\e-\sum_{j \underset{\mathcal{R}_{\infty,\e}}{\nsim} i} v(\a_i - \a_j))\right)
$$
\end{lemm}

\begin{rema}
Bien sûr dans le lemme précédent il est inutile d'aller jusqu'à $\widehat{\overline{K}^\times}$. Il suffit d'étendre les scalaires à n'importe quelle extension de $K$ contenant toutes les racines de $f$.
\end{rema}

La démonstration du lemme précédent repose elle même sur le lemme suivant appliqué de façon récurrente.

\begin{lemm}
Soit $J$ un ensemble fini, $(\beta_j)_{j\in J}\in \overline{K}^J$ et
$\eta\in v(\overline{K})_{\geq 0}$. Définissons la relation d'équivalence suivante sur $J$ : $j_1\sim j_2$ si $v(\b_{j_1}-\b_{j_2})\geq \frac{\eta}{|J|}$. Alors
$$
\{ x\in\overline{K}\;|\; \sum_{j\in J} v(x-\beta_j)\geq \eta \}
= \coprod_{[j]\in J/\sim } \{ x\in\overline{K}\;|\; \sum_{j'\sim j} v(x-\beta_{j'})\geq \eta - \sum_{j'\nsim j} v(\b_{j}-\b_{j'}) \}
$$
Si $\forall j_1,j_2\in J\; \; j_1\sim j_2$ alors 
$$
\forall j_0\in J\; \;\; \{x\in \overline{K}\;|\; \sum_{j\in J} v(x-\b_j)\geq \eta \} =
\{x\in\overline{K}\;|\; v(x-\beta_{j_0}) \geq \frac{\eta}{|J|}\}
$$
\end{lemm}

\subsection{Dualité convexe et polygone de Newton}

\begin{defi}
Soit $\psi: \R \ldrt \R \cup \{\infty \}$ une fonction convexe
 semi-continue inférieurement. 
 Posons 
$$
\psi^* (s)= \text{sup} \{t\;|\; \psi(\bullet) \geq -s \bullet +t \}
$$
\end{defi}

La fonction $\psi$ se déduit de $\psi^*$ par 
$$
\psi (t)= \text{sup} \{- xt + \psi^* (x) \;|\; x\in \R \}
$$
(le graphe de $\psi$ est une ``enveloppe'' de droites). On écrit cela sous la forme $\psi=\psi^{**}$. 

\begin{defi}
Notons $\text{Newt}:[0 , \deg f]\ldrt [0,+\infty]$ le polygone de Newton de $f$. Il s'agit d'une fonction convexe linéaire par morceaux. Si $f=\sum_{k} a_k T^k$ son graphe est l'enveloppe convexe des $(k,v(a_k))_k$. Ses pentes sont les opposés des valuations des racines de $f$ (comptées avec multiplicité). 
\end{defi}

\begin{lemm}
$$
\text{Newt}^{\,*} (s) = \text{inf} \underbrace{\{v(f(x))\;|\; x\in \overline{K}\;\;\; v(x)\geq s
  \; \}}
_{\text{un intervalle de } v(\overline{K}) }
$$
\end{lemm}
\dem
On vérifie facilement que 
$$
\inf \{v(f(x))\;|\;x\in\overline{K}\;\; v(x)\geq s\} = \sum_{i\in I} \inf \{ s,v(\a_i) \}
$$
\qed

On a donc en termes de géométrie rigide $$f (\mathbb{B} (0,s)) =
\mathbb{B} (0,\text{Newt}^* (s))$$

\subsection{Fonction de Herbrand}

Supposons que $f(0)=0$.

\begin{defi}
Soit  $X^\e$ la composante connexe géométrique de $0$ dans
$Y^\e=f^{-1}(\mathbb{B} (0,\e))$. 
\end{defi}

\begin{prop}\label{Ktdvf}
Si l'on pose $\eta (\bullet) = \text{Newt}^{\,*}(\bullet)$ et $ \psi= \eta^{-1}$
on a les égalité 
$$
X^{\eta(\e)} = \mathbb{B} (0,\e) \; \text{et } \; X^\e= \mathbb{B}(0,\psi (\e))
$$
\end{prop}
\dem
On vérifie aisément que la fonction $\eta : [0,+\infty [\ldrt
[0,+\infty [$ est strictement croissante et que donc $\psi$ est bien
définie.  La seconde égalité résulte de la première. On a d'après le
lemme précédent 
$$
X^{\eta (\e)} = f^{-1} ( f (\mathbb{B} (0,\e)))^0
$$
où $(-)^0$ désigne la composante connexe de $0\in \mathbb{B}^1$. Bien sûr 
$\mathbb{B} (0,\e) \subset  f^{-1} ( f (\mathbb{B} (0,\e)))$ et donc 
$\mathbb{B} (0,\e) \subset  f^{-1} ( f (\mathbb{B} (0,\e)))^0$. De
plus on sait qu'il existe $\e'>0$ tel que $f^{-1} ( f (\mathbb{B}
(0,\e)))^0=\mathbb{B} (0,\e')$ avec nécessairement $\e'\geq \e$. 
 Mais $f (\mathbb{B} (0,\e'))=
\mathbb{B} (0,\eta (\e'))$ or si $\e'>\e \; \eta (\e')>\eta
(\e)$. Donc $\e'=\e$.
\qed
\\

La fonction $\eta$ se calcule de façon usuelle par une intégrale :

\begin{lemm}
$$
\eta (s)= \int^s_0 \text{Card} \left ( f^{-1} (0) \cap \mathbb{B} (0,s) \right ) ds
$$
\end{lemm}

\begin{lemm}\label{korTo}
 Soient $f$ et $g$ deux polynômes comme précédemment. Alors 
$$
\text{Newt}_{f\circ g}^{\,*} = \text{Newt}_f^{\,*} \circ (\text{Newt}_g^{\,*})
$$
En particulier $\text{Newt}_f$ et $\text{Newt}_g$ déterminent
complètement $\text{Newt}_{f\circ g}$. La fonction $\eta$ associée à
$f\circ g$ est la composée de celle associée à $f$ avec celle associée
à $g$.  
\end{lemm}
\dem
On a 
$$
(f\circ g) (\mathbb{B} ( 0,\e)) = f ( \mathbb{B} ( 0, \text{Newt}_g^*
(\e))) = \mathbb{B}\left ( 0, \text{Newt}_f^* ( \text{Newt}_g^*
  (\e))\right )
$$
qui est égal à $\mathbb{B} ( 0, \text{Newt}_{f\circ g}^* (\e))$. 
\qed

\subsection{La filtration de ramification des schémas en groupes monogènes}

Soit $G=\spec ( \O_K [T] /(f))$ où $f(T)=T\prod_i (T-\a_i)$ un schéma
en groupe fini localement libre sur $\O_K$ dont l'algèbre est monogène.

\begin{rema}
Soit $H$ un groupe formel $p$-divisible sur $\O_K$ et $G\subset H$ un sous-groupe fini localement libre. Alors l'algèbre de $G$ est monogène.  En effet, après un bon choix de coordonnées formelles $T$ à la source et au but l'isogénie de groupes formels $H\ldrt H/G$ s'écrit $T\longmapsto f(T)$ pour un polynôme unitaire $f$.  
\end{rema}

\begin{defi}
On pose $G_a= \{ \a_i \; |\; v(\a_i)\geq a\; \}$ vu comme sous-groupe
étale de $G_\eta$. On note également $G_a= \spec (\O_K [T]/(T\prod_{i,
  v(\a_i)\geq a} (T-\a_i)))$ son adhérence schématique dans $G$. On
pose $G^a = G_{\psi (a)}$ où $\psi$ est la fonction de Herbrand
associée à $f$. 
\end{defi}

D'après la proposition \ref{Ktdvf}
le groupe $G^a$ est celui défini par Abbes et Saito, et étudié par
Abbes-Mokrane (\cite{AbbesMokrane}).

\begin{prop}\label{zuytok}
Soit $\ph : H_1\ldrt H_2$ une isogénie de groupes formels $p$-divisibles de dimension
1 sur $\O_K$. Soit $G\subset H_1$ plat fini et $\ph (G)$ l'adhérence schématique de $\ph (G_\eta)$. Alors
$$
\forall a\;\; \ph (G^a)  = \ph (G)^a
$$
\end{prop}
\dem
Posons $G=\spec ( \O_K [T]/ (f\circ g))$, $\ph (G)= \spec (\O_K [T]/
(f))$ et $\ph : G \ldrt \ph (G)$ donné au niveau des algèbres par
$T\mapsto g(T)$. Alors, 
$G^a$ est la composante connexe de $0$ dans $U=\{y\in \mathbb{B}^1
(\overline{K}) \;|\; v (f\circ g (y))\geq a \}$ et $\ph (G)^a$ la
composante connexe de $0$ dans $V= \{y\in \mathbb{B}^1
(\overline{K}) \;|\; v (f (y))\geq a \}$. Il est clair que $\ph (V)
\subset U$. De plus l'image par $\ph$ de la composante connexe contenant $0$ de
$V$ est
encore connexe, contient zéro, et donc est contenue dans $\ph (G)^a$. Donc 
$\ph (G^a)  \subset \ph (G)^a$. Pour l'autre inclusion il suffit de
voir qu'étant donné que $\ph$ induit un morphisme étale fini d'espaces
rigides le morphisme $\ph_{|V} : V\ldrt U$ est étale fini et qu'il
l'est donc encore en restriction à chacune des composantes
connexes. Or un morphisme étale fini entre espaces rigides connexes
est surjectif. 
\qed
\vspace{7mm}

Les groupes de ramifications supérieures se comportent donc bien par
isogénies, mais sont moins concrets que ceux de ramification
inférieure.
Certains des calculs de cet article se ramènent donc en fait à jouer
sur les deux plans en passant de la ramification inférieure à la
ramification
supérieure. 

\bibliographystyle{plain}
\bibliography{biblio}
\end{document}